\documentclass{amsart}
\usepackage[margin=1.0in]{geometry}
\usepackage{amsthm, amsmath, amssymb, mathtools, mathrsfs}
\usepackage{aliascnt}
\usepackage{hyperref}
\usepackage{breakurl}
\usepackage{array, longtable}
\usepackage{makecell}
\usepackage[normalem]{ulem} % for crossing out sentences

\makeatletter
\def\LT@start{%
\let\LT@start\endgraf
\endgraf\penalty\z@\vskip\LTpre
\dimen@\pagetotal
\advance\dimen@ \ht\ifvoid\LT@firsthead\LT@head\else\LT@firsthead\fi
\advance\dimen@ \dp\ifvoid\LT@firsthead\LT@head\else\LT@firsthead\fi
\advance\dimen@ \ht\LT@foot
\dimen@ii\vfuzz
\vfuzz\maxdimen
\setbox\tw@\copy\z@
\setbox\tw@\vsplit\tw@ to \ht\@arstrutbox
\setbox\tw@\vbox{\unvbox\tw@}%
\vfuzz\dimen@ii
\advance\dimen@ \ht
\ifdim\ht\@arstrutbox>\ht\tw@\@arstrutbox\else\tw@\fi
\advance\dimen@\dp
\ifdim\dp\@arstrutbox>\dp\tw@\@arstrutbox\else\tw@\fi
\advance\dimen@ -\pagegoal
\ifdim \dimen@>\z@\unskip\vfil\break\fi
\global\@colroom\@colht
\ifvoid\LT@foot\else
\advance\vsize-\ht\LT@foot
\global\advance\@colroom-\ht\LT@foot
\dimen@\pagegoal\advance\dimen@-\ht\LT@foot\pagegoal\dimen@
\maxdepth\z@
\fi
\ifvoid\LT@firsthead\copy\LT@head\else\box\LT@firsthead\fi
\output{\LT@output}}
\makeatother

\emergencystretch=\maxdimen

\usepackage[textsize=footnotesize]{todonotes}

\usepackage{xcolor}
\definecolor{dblue}{rgb}{0,0,0.70}
\hypersetup{
	unicode=true,
	colorlinks=true,
	citecolor=dblue,
	linkcolor=dblue,
	anchorcolor=dblue
}

\makeatletter
\expandafter\g@addto@macro\csname th@plain\endcsname{%
	\thm@notefont{\bfseries}
}%
\expandafter\g@addto@macro\csname th@remark\endcsname{%
	\thm@headfont{\bfseries}
}%
\makeatother

\newtheorem{theorem}{Theorem}[section]	
\newtheorem*{theorem*}{Theorem}

\newaliascnt{lemma}{theorem}
\newtheorem{lemma}[lemma]{Lemma}
\aliascntresetthe{lemma}
\newtheorem*{lemma*}{Lemma}

\newaliascnt{proposition}{theorem}
\newtheorem{proposition}[proposition]{Proposition}
\aliascntresetthe{proposition}

\newaliascnt{corollary}{theorem}
\newtheorem{corollary}[corollary]{Corollary}
\aliascntresetthe{corollary}

\theoremstyle{remark}

\newaliascnt{remark}{theorem}
\newtheorem{remark}[remark]{Remark}
\aliascntresetthe{remark}
\newaliascnt{question}{theorem}
\newtheorem{question}[question]{Question}
\aliascntresetthe{question}

\newtheorem*{question*}{Question}

\newaliascnt{definition}{theorem}
\newtheorem{definition}[definition]{Definition}
\aliascntresetthe{definition}

\newaliascnt{example}{theorem}

\aliascntresetthe{example}

\newaliascnt{convention}{theorem}
\newtheorem{convention}[convention]{Convention}
\aliascntresetthe{convention}

\makeatletter
\def\l@subsection{\@tocline{2}{0pt}{1pc}{5pc}{}} \def\l@subsection{\@tocline{2}{0pt}{2pc}{6pc}{}}
\makeatother
\usepackage[style = numeric, sorting = nyt]{biblatex}
\newcommand{\lr}{\leftrightarrow}
\newcommand{\Ord}{\mathrm{Ord}}

\newcommand{\dom}{\operatorname{dom}}
\newcommand{\ran}{\operatorname{ran}}
\newcommand{\rank}{\operatorname{rank}}
\newcommand{\crit}{\operatorname{crit}}
\DeclareMathOperator{\restricts}{\upharpoonright}
\DeclareMathOperator{\downwards}{\downarrow}
\newcommand{\mySigma}{\Sigma}

\newcommand{\TC}{\operatorname{TC}}

\newcommand{\mv}{\operatorname{mv}}
\newcommand{\mhyphen}{\text{-}}

\newcommand{\cf}{\operatorname{cf}}
\newcommand{\Def}{\operatorname{Def}}

\newcommand{\up}{\mathsf{up}}
\newcommand{\op}{\mathsf{op}}

\newcommand{\cminus}{\raisebox{.33\height}{\scalebox{0.75}{\ensuremath{-}}}}
\newcommand{\ZFCminusrep}{\ensuremath{\mathsf{ZFC}\cminus}}
\newcommand{\ZFminusrep}{\ensuremath{\mathsf{ZF}\cminus}}

\newcommand{\ZF}{\mathsf{ZF}}
\newcommand{\ZFC}{\mathsf{ZFC}}
\newcommand{\CZF}{\mathsf{CZF}}
\newcommand{\IZF}{\mathsf{IZF}}
\newcommand{\KP}{\mathsf{KP}}
\newcommand{\IKP}{\mathsf{IKP}}
\newcommand{\BTEE}{\mathsf{BTEE}}
\newcommand{\WA}{\mathsf{WA}}
\newcommand{\CGB}{\mathsf{CGB}}
\newcommand{\IGB}{\mathsf{IGB}}
\newcommand{\TR}{\mathsf{TR}}

\newcommand{\ZFCminus}{\mathsf{ZFC}^-}
\newcommand{\ZFminus}{\mathsf{ZF}^-}
\newcommand{\CZFminus}{\mathsf{CZF}^-}
\newcommand{\TIj}{\text{Set Induction\textsubscript{$j$}}}

\newcommand{\dnot}{{\lnot\lnot}}
\newcommand{\Low}{\operatorname{Low}}

\newcommand{\llbr}{[\mkern-2mu[}
\newcommand{\rrbr}{]\mkern-2mu]}
\newcommand{\bigllbr}{\big[\mkern-4mu\big[}
\newcommand{\bigrrbr}{\big]\mkern-4mu\big]}
\newcommand{\lag}{\langle}
\newcommand{\rag}{\rangle}

\newcommand{\rrarrows}{\rightrightarrows}
\newcommand{\lrlrarrows}{\mathrel{\mathrlap{\leftleftarrows}{\rightrightarrows}}}
\newcommand{\Ind}{\mathsf{Ind}}

\makeatletter
\def\oversortoftilde#1{\mathop{\vbox{\m@th\ialign{##\crcr\noalign{\kern3\p@}%
				\sortoftildefill\crcr\noalign{\kern3\p@\nointerlineskip}%
				$\hfil\displaystyle{#1}\hfil$\crcr}}}\limits}

\title{Very large set axioms over constructive set theories}

\author{Hanul Jeon}
\email{ \href{mailto:hj344@cornell.edu}{hj344@cornell.edu}}
\urladdr{ \href{https://hanuljeon95.github.io}{https://hanuljeon95.github.io} }
\address{Department of Mathematics, Cornell University, Ithaca, NY 14853} 
\subjclass[2010]{Primary 03E70; Secondary 03E55}

\author{Richard Matthews}
\email{ \href{richard.matthews@u-pec.fr}{richard.matthews@u-pec.fr}}
\urladdr{ \href{https://richardmatthewslogic.github.io/}{https://richardmatthewslogic.github.io/} }
\address{Universit\'{e} Paris-Est Cr\'{e}teil, LACL, F-94010}
\thanks{The second author was supported by the UK Engineering and Physical Sciences Research Council during this research and is grateful for their support.}

\begin{document}
	\maketitle
	
	\begin{abstract}
		We investigate large set axioms defined in terms of elementary embeddings over constructive set theories, focusing on $\mathsf{IKP}$ and $\mathsf{CZF}$. Most previously studied large set axioms, notably, the constructive analogues of large cardinals below $0^\sharp$, have proof-theoretic strength weaker than full Second-order Arithmetic.
    	On the other hand, the situation is dramatically different for those defined via elementary embeddings. We show that by adding to $\mathsf{IKP}$ the basic properties of an elementary embedding $j\colon V\to M$ for $\Delta_0$-formulas, which we will denote by $\Delta_0\text{-}\mathsf{BTEE}_M$, we obtain the consistency of $\mathsf{ZFC}$ and more. We will also see that the consistency strength of a Reinhardt set exceeds that of $\mathsf{ZF+WA}$. Furthermore, we will define super Reinhardt sets and $\mathsf{TR}$, which is a constructive analogue of $V$ being totally Reinhardt, and prove that their proof-theoretic strength exceeds that of $\mathsf{ZF}$ with choiceless large cardinals.
	\end{abstract}
	
	\tableofcontents
	
	\section{Introduction}
	Large cardinals have played a pivotal role in set theory, and many of them can be defined in terms of elementary embeddings.
	Associating large cardinals with elementary embeddings appeared first in Scott's pioneering paper \cite{Scott1961}, and was further systematically developed throughout the 1960s and 1970s. Many of these results were collected by Reinhardt and Solovay around the early 1970s, which was published later with Kanamori in the expository paper \cite{SolovayReinhardtKanamori1978}. 
	
	The attempt to find ever stronger notions of large cardinal axioms culminated in the principle now known as a \emph{Reinhardt cardinal}, which was first mentioned in Reinhardt's doctoral thesis \cite{ReinhardtPhD}. 
	A Reinhardt cardinal is a critical point of a non-trivial elementary embedding $j\colon V\to V$. 
	Unfortunately, the fate of a Reinhardt cardinal in $\ZFC$ is that of inconsistency, as if Icarus falls into the sea as he flew too close to the sun. 
	A famous result by Kunen \cite{Kunen1971} proves that Reinhardt cardinals are incompatible with the Axiom of Choice and it is still unknown if $\ZF$ with a Reinhardt cardinal is consistent. 
	However, there has been little study of the consistency of Reinhardt cardinals in the choiceless context and 
	{few results about the implications of such axioms have appeared in the literature before 2010. T}he only exception the authors know of is a result of Apter and Sargsyan \cite{ApterSargsyan2004}, and other relevant results about Reinhardt embeddings focused on its inconsistency. 
	Notable examples include Suzuki's non-definability of embeddings $j\colon V\to V$ over $\ZF$ \cite{Suzuki1999} or Zapletal's PCF-theoretic proof of Kunen's inconsistency theorem \cite{Zapletal1996}. \\
	
	Choiceless large cardinals are large cardinal notions that extend a Reinhardt cardinal and are therefore incompatible with the Axiom of Choice.
	A super Reinhardt cardinal was employed by Hugh Woodin in 1983 to prove the consistency of $\mathsf{ZFC+I_0}$, which had a focal role in establishing the consistency of $\mathsf{ZF+AD}^{L(\mathbb{R})}$. 
	A Berkeley cardinal appeared around 1992 by Woodin at his set theory seminar as an attempt to provide a large cardinal notion that was refutable from $\ZF$ alone. While no such inconsistency has been found so far, it has since become an interesting principle in itself.
	Current research on choiceless large cardinals emerged in the mid-2010s as part of a project to explore Woodin's $\mathsf{HOD}$ dichotomy (see \cite[Section 7.1]{Woodin2010} or \cite{BagariaKoellnerWoodin2019} for details), under the thesis that such cardinals would indicate the $V$ was `far' from $\mathsf{HOD}$ in some sense.
	Bagaria, Koellner, and Woodin collected and analyzed notions of choiceless large cardinals in \cite{BagariaKoellnerWoodin2019}, and the theory of choiceless large cardinals was further developed in various papers by authors including Cutolo, Goldberg, and Schlutzenberg. 
	
	One of the most striking results along this line is a result by Goldberg \cite{Goldberg2021EvenOrdinals}, which establishes the consistency of $\ZF + j\colon V_{\lambda+2}\to V_{\lambda+2}$ modulo large cardinals over $\ZFC$:
	
	\begin{theorem}[Goldberg {\cite[Theorem 6.20]{Goldberg2021EvenOrdinals}}] \pushQED{\qed}
	    These two theories are equiconsistent over $\mathsf{ZF+DC}$:
	    \begin{enumerate}
	        \item For some ordinal $\lambda$, there is an elementary embedding $j\colon V_{\lambda+2}\to V_{\lambda+2}$.\footnote{Karagila proposed the term \emph{Kunen cardinal} for a critical point of an elementary embedding $j\colon V_{\lambda+2}\to V_{\lambda+2}$ since Kunen's result \cite{Kunen1971} shows no such cardinal can exist in $\ZFC$.}
	        \item $\mathsf{AC+I_0}$. \qedhere
	    \end{enumerate}
	\end{theorem}
	
	Furthermore, Goldberg proved that the existence of an elementary embedding $j\colon V_{\lambda+3}\to V_{\lambda+3}$ exceeds almost all of the traditional large cardinal hierarchy over $\ZFC$.
	
	\begin{theorem}[Goldberg {\cite[Theorem 6.16]{Goldberg2021EvenOrdinals}}] \label{theorem:Goldberg}
	    Working over $\mathsf{ZF+DC}$, the existence of a $\Sigma_1$-elementary embedding $j\colon V_{\lambda+3}\to V_{\lambda+3}$ implies the consistency of $\ZFC+\mathsf{I_0}$.
	\end{theorem}
	
	% Elementary embeddings over weak set theories
	In the other direction, we can try to `salvage' an elementary embedding $j\colon V\to V$ by weakening the background set theory. One direction of research in this manner was conducted by Corazza. Corazza \cite{Corazza2000} introduced the \emph{Wholeness axiom}, $\WA$, by dropping Replacement for $j$-formulas. $\WA$ is known to be weaker than $\mathsf{I}_3$. Corazza further weakened $\WA$ to the \emph{Basic Theory of Elementary Embeddings}, $\BTEE$, in \cite{Corazza2006} which is the weakest setting one needs to express the existence of an elementary embedding by dropping all axioms for $j$-formulas in the extended language. The resulting axiom is known to be weaker than the existence of $0^\sharp$.
	
	Another direction to weaken the assumptions is by dropping Powerset. However, it should be noted that $\ZFCminusrep$, the theory obtained by ejecting Powerset from $\ZFC$ and only assuming Replacement, is ill-behaved with elementary embeddings. For example, in \cite{GitmanHamkinsJohnstone2016} it is shown that \L o\'s's theorem can fail over $\ZFCminusrep$, and a cofinal $\Sigma_1$-elementary embedding need not be fully elementary. On the other hand, they also show that these issues can be avoided by strengthening $\ZFCminusrep$ to $\ZFCminus$, which is obtained by additionally assuming Collection.
	
	Further research along this line is characterizing large cardinal notions in terms of models of $\ZFCminus$ with an ultrafilter predicate (For example, \cite{HolyLucke2021} or \cite{GitmanSchlichtUnpublished}). 
	Large cardinals defined in this way refine the large cardinal hierarchy between a Ramsey cardinal and a measurable cardinal, and provide bounds for the consistency strength of $\ZFCminus$ with an elementary embedding. As one such example of this characterization, we have the following theorem.
	
	\begin{theorem}[{\cite[Theorem 10.5.7]{MatthwesPhD}}] \label{theorem:ZFC-CriticalStrength}
	    $\ZFC$ with a locally measurable cardinal proves the consistency of $\ZFCminus + \mathsf{DC}_{\mathrm{<Ord}}$ plus the existence of a non-trivial elementary embedding $j\colon V\to M$. 
	\end{theorem}
	
	Finally, it can be shown that an elementary embedding $j \colon V \rightarrow V$ is possible in $\ZFCminus$ under the large cardinal assumption of $\mathsf{ZFC + I_1}$: 

    \begin{theorem}[{\cite[Theorem 9.3.2]{MatthwesPhD}} or {\cite[Theorem 2.2]{Matthews2020}}]
	    $\ZFC$ proves the following: there is an elementary embedding $k\colon V_{\lambda+1}\to V_{\lambda+1}$ if and only if there is an elementary embedding $j\colon H_{\lambda^+}\to H_{\lambda^+}$.
	    
	    As a corollary, $\ZFCminus_j$\footnote{Where $\ZFCminus_j$ is $\ZFCminus$ in the langauge expanded to include $j$. See \autoref{Convention:T_j}.} is compatible (modulo large cardinals over $\ZFC$) with a non-trivial elementary embedding $j\colon V\to V$ which additionally satisfies that $V_{\crit j}$ exists. 
	\end{theorem}
	
	However, \cite{Matthews2020} also showed there is a limitation on the properties of such an elementary embedding \mbox{$j\colon V\to V$} over $\ZFCminus$. One of these restrictions is that $j\colon V\to V$ cannot be cofinal\footnote{{An embedding $j \colon M \rightarrow N$ is said to be \emph{cofinal} if for every $y \in N$ there is an $x \in M$ such that $y \in j(x)$.}}:
	
	\begin{theorem}[{\cite[Theorem 10.2.3]{MatthwesPhD}} or {\cite[Theorem 5.4]{Matthews2020}}]
	    Working over $\ZFCminus_j$, if $j\colon V\to V$ is a non-trivial $\Sigma_0$-elementary embedding such that $V_{\crit j}$ exists, then it cannot be cofinal.
	\end{theorem}
	
	% Large set axioms over constructive set theories
	We may also weaken the background set theory by dropping the law of excluded middle, that is, moving into a constructive setting. Since some statements are no longer equivalent over a constructive background, we need to carefully formulate constructive counterparts of classical notions, including the axiom systems of constructive set theories.
	
	The first form of a constructive set theory was defined by H. Friedman \cite{Friedman1973}, and is now known as \emph{Intuitionistic $\ZF$}, $\IZF$. In \cite{Friedman1973}, Friedman showed that $\IZF$ and $\ZF$ are mutually interpretable by using the combination of double-negation translation and non-extensional set theory.
	Another flavor of constructive set theory appeared as an attempt to establish a formalization of Bishop's Constructive analysis. Myhill \cite{Myhill1975} gave his own formulation of a constructive set theory $\mathsf{CST}$. However, its language is different from the standard one -- $\mathsf{CST}$ includes natural numbers as primitive objects, while $\ZF$ does not.
	The other formulation of a constructive set theory, which is closer to standard $\ZF$, was given by Aczel \cite{Aczel1978} via a type-theoretic interpretation. Aczel's constructive set theory is called \emph{Constructive $\ZF$,} $\CZF$. Aczel further developed a theory on $\CZF$ and its relationship with Martin-L\"of type theory in the consequent works \cite{Aczel1982} and \cite{Aczel1986}. In particular, Aczel's last paper in the sequel, \cite{Aczel1986}, defined \emph{regular sets}, which begins the program to define large cardinals over $\CZF$.
	
	The first research on constructive analogues of large cardinal axioms was done by Friedman and \v{S}\v{c}edrov \cite{FriedmanScedrov1984}. Unlike classical set theories, ordinals over constructive set theories are not well-behaved. This motivated defining large cardinal notions over constructive set theories in a structural manner, resulting in \emph{large set axioms}.
	They defined and analyzed inaccessible sets, Mahlo sets, and various elementary embeddings over $\IZF$, and proved that their consistency strength is no different from their classical counterparts.
	
	Large set axioms over constructive set theories appeared first in various papers of Rathjen (for example, \cite{Rathjen1993}, \cite{Rathjen1998}, \cite{RathjenGrifforPalmgren1998}, \cite{Rathjen1999Realm}). The first appearance of large set axioms over $\CZF$ was in a proof-theoretic context, and their relationship with better known theories, like extensions of Martin-L\"of type theories or $\KP$, were emphasized. 
	An initial analysis of large set axioms over $\CZF$ can be found in \cite{AczelRathjen2001} or \cite{AczelRathjen2010}, and this has been further extended by Gibbons \cite{GibbonsPhD} and Ziegler \cite{ZieglerPhD} in each of their doctoral theses.
	Gibbons \cite{GibbonsPhD} extended Rathjen's analysis on the proof-theoretic strength of $\CZF$ in \cite{Rathjen1993} to $\CZF$ with Mahlo sets. Furthermore, Gibbons' thesis is the first publication that shows the definition of a critical set\footnote{He called it a \emph{measurable set}. We will address this terminology in \autoref{Section:LargeLargeSet}.}
	
	Unlike other results around large set axioms over $\CZF$, Ziegler's thesis focused on what we can derive about large set axioms from $\CZF$ alone. For example, he observed that the number of inaccessible sets may not affect the proof-theoretic strength:
	
	\begin{theorem}[{\cite[Chapter 6]{ZieglerPhD}}] \label{theorem:ZieglerInaccessibleStrength}
	    The following theories are equiconsistent:
	    \begin{enumerate}
	        \item $\CZF$ with an inaccessible set,
	        \item $\CZF$ with two inaccessible sets, and
	        \item $\CZF$ with $\omega$ inaccessible sets.
	    \end{enumerate}
	\end{theorem}
	
	Ziegler also examined elementary embeddings over $\CZF$ in detail, and one of his striking results is that every Reinhardt embedding $j\colon V\to V$ must be cofinal (see \autoref{Proposition:CZFReinhardtEmbeddingCofinal} for the formal statement of the theorem). 
	
	Ziegler's thesis ends with the following result that elementary embeddings are incompatible with the principle of \emph{subcountability}, which asserts that for every set there is a partial surjection from $\omega$ onto it. This can be seen as a constructive analogue to Scott's early result in \cite{Scott1961} that measurable cardinals are incompatible with the Axiom of Constructibility.
	
	\begin{theorem}[Ziegler {\cite[Theorem 9.93]{ZieglerPhD}}] \pushQED{\qed}
	    Over $\mathsf{{C}ZF}$, the combination of the Axiom of Subcountability and the existence of a critical set results in a contradiction. \qedhere
	\end{theorem}
	
	Rathjen analyzed the proof-theoretic strength of \emph{small large sets}, large set notions whose classical counterpart is weaker than the existence of $0^\sharp$, over $\CZF$ (See \autoref{Section:SmallLargeSet}, especially \autoref{Proposition:Prooftheoreticstrength-smalllargesets}), and the proof-theoretic strength of all currently known small large set axioms over $\CZF$ is weaker than that of Second-order Arithmetic. However, the proof-theoretic strength of \emph{large large sets}, large set notions defined in terms of elementary embeddings, is yet to receive a formal rigorous treatment.
	
	We end this introduction by noting some history of the development of this work. A first version of this paper can be found in \cite{Jeon2021CZF} where the first author studied the consistency strength of a Reinhardt set over $\CZF$ with Full Separation. The second author's PhD thesis, \cite{MatthwesPhD}, included an initial investigation into elementary embeddings of $\KP$ and $\IKP$. This version can be seen as a combination of the previous work by the individual authors and extensively extends the results found in either source.
	
	\subsection*{Main results}
	
	It turns out that the proof-theoretic strength of large large set vastly exceeds that of $\ZFC$. In fact, we just need a small fragment of the properties of an elementary embedding{, that we shall denote by} $\Delta_0\mhyphen\BTEE_M$, which is the minimal theory needed to claim that $j\colon V\to M$ is a $\Delta_0$-elementary embedding, to exceed the proof-theoretic strength of $\ZFC$.
	
	\begin{theorem*}
	\phantom{a}
	    \begin{enumerate}
	        \item (\autoref{Theorem:CriticalPointModelsIZF}) Working over $\IKP$, let $K$ be a transitive set such that $K\models \Delta_0\mhyphen\mathsf{Sep}$, $\omega\in K$ and let $j \colon V\to M$ be a $\Delta_0$-elementary embedding whose critical point is $K$. Then $K\models \IZF$.
	        
	        \item (\autoref{Corollary:LambdaModelsIZFBTEEInd}) Furthermore, if we additionally allow Set Induction and Collection for $\Sigma^{j,M}$-formulas and add Separation for $\Delta_0^{j,M}$-formulas, then we can define $j^\omega(K):=\bigcup_{n\in\omega} j^n(K)$ and prove that $j^\omega(K)$ satisfies $\mathsf{IZF+BTEE}$ plus Set Induction for $j$-formulas.
	        
	        \item (\autoref{Theorem:IKPSigmaOrdimpliesIZF+BTEE} and \autoref{Theorem:InterpretationResults}) As a consequence, the following two theories prove the consistency of $\mathsf{ZFC+BTEE}$ {plus Set Induction for $j$-formulas}: $\CZF_{j, M}$ with a critical set, and $\IKP_{j,M}$ plus a \mbox{$\Sigma$-$\Ord$-inary} elementary embedding with a critical point $\kappa\in\Ord$.
	    \end{enumerate}
	\end{theorem*}
	
	The next natural question would be how strong a Reinhardt set is. It turns out that Reinhardt sets are very strong over $\CZF$:
	
	\begin{theorem*}(\autoref{Theorem:CZF+ReinhardtinterpretsZF+WA})
	    $\CZF$ with a Reinhardt set proves the consistency of $\mathsf{ZF+WA}$.
	\end{theorem*}
	
	These two results motivate an idea that the proof-theoretic strength of stronger large set notions may go beyond that of $\ZF$ with choiceless large cardinals.
	{This idea turns out to hold}, and a constructive formulation of super Reinhardt cardinals witnesses this.
	We can push this idea further, so that we can expect that we may reach an `equilibrium' by strengthening large set axioms once more, in the sense that adding some large set notions to $\CZF$ has the same proof-theoretic strength as that of $\IZF$ plus the same axiom. We can see that the assertion that $V$ resembles $V_\kappa$ for a total Reinhardt cardinal $\kappa$, which we will call $\mathsf{TR}$, witnesses this claim.
	
	However, both our analogues of super Reinhardts and total Reinhardts require a second-order formulation. We will resolve this issue by defining $\CGB$, which is a constructive version of $\mathsf{GB}$. Moreover, formulating an elementary embedding in a constructive context requires infinite connectives because there is no obvious way to cast the elementarity into a single formula. (It is possible in a classical context because every $\Sigma_1$-elementary embedding $j\colon V\to M$ is fully elementary.)
	This motivates $\CGB_\infty$, $\CGB$ with infinite connectives. Extending $\CGB$ with infinite connectives will turn out to be `harmless' in the sense that $\CGB_\infty$ is conservative over $\CGB$. Under this setting, we have the following results:
	
	\begin{theorem*}
	    \phantom{a}
	    \begin{enumerate}
	        \item (\autoref{Theorem:LowerBound-superReinhardt}) $\CGB_\infty$ with a super Reinhardt set proves the consistency of $\ZF$ with a Reinhardt cardinal,
	        \item (\autoref{Corollary:CGBTRinterpretsGBTR}) $\CGB_\infty + \mathsf{TR}$ interprets $\mathsf{ZF+TR}$.
	    \end{enumerate}
	\end{theorem*}
	
    \subsection*{The structure of the paper} 
    This paper is largely divided into two parts: the `internal' analysis of large set axioms over constructive set theories, and deriving a lower bound for large large set axioms in terms of extensions of classical set theories.
    Defining some of these axioms will require largely unexplored notions including second-order constructive set theories, so we define the necessary preliminary notions in \autoref{Section:Prelim}.
    
    Next, we review the basic properties of large set axioms over constructive set theories. This is also divided into two sections: in \autoref{Section:SmallLargeSet}, we review small large sets over constructive set theories and their connection with classical set theories. 
    Then in \autoref{Section:LargeLargeSet}, we define large large sets over constructive set theories, including analogues of choiceless large cardinals.
    Having laid down the necessary framework, in \autoref{Section:StructuralLowerbdd} we provide an internal analysis of large set axioms. The main consequence of \autoref{Section:StructuralLowerbdd} is that we provide lower bounds for the consistency strength of large large set axioms in terms of extensions of $\IZF$, from which it will be easier to interpret classical theories by applying a double-negation translation.
    
    In any case, we want to derive the consistency strengths in terms of extensions of classical set theories. Hence we need to develop one of the methods to transform intuitionistic theories into classical ones.
    This is what \autoref{Section:HeytingInterpretation} mainly focuses on. In \autoref{Section:HeytingInterpretation}, we review Gambino's Heyting-valued interpretation defined in \cite{Gambino2006} and investigate this interpretation under the double-negation topology. We will see that the interpretation translates $\IZF$ into the classical theory $\ZF$. Furthermore, we will also reduce extensions of $\IZF$ to those of $\ZF$ by using the double-negation topology.
    \autoref{Section:DNT-SOST} is devoted to the double-negation translation of second-order set theories and the concept of the universe being totally Reinhardt.
    We summarize the lower bound of the consistency strength of large large set axioms over constructive set theories in terms of extensions of $\ZFC$ in \autoref{Section:ConsistencyStrength:Final} before ending by posing some questions for future investigation in \autoref{Section:RemarkQuestions}.

	\section{Preliminaries}\label{Section:Prelim}
	In this section, we will briefly review $\ZFC$ without Power Set, $\ZFCminus$, and constructive set theory. There are various formulations of constructive set theories, but we will focus on $\CZF$.
	In addition, we will define the second-order variant $\CGB$ and $\IGB$ of $\CZF$ and $\IZF$ respectively.
	
	\subsection{\texorpdfstring{$\ZFC$}{ZFC} without Power Set}
		We will frequently mention $\ZFC$ without Power Set, denoted $\ZFCminus$. However, $\ZFCminus$ is not obtained by just dropping Power Set from $\ZFC$: 
	
	{
	\begin{definition}
		$\ZFminus$ is the theory obtained from $\ZF$ by dropping Power Set and using Collection instead of Replacement. $\ZFCminus$ is obtained by adding the Well-Ordering Principle to $\ZFminus$.
	\end{definition}
	}
	
	Note that using Collection instead of Replacement is necessary to avoid pathologies. See \cite{GitmanHamkinsJohnstone2016} for the details. It is also known by \cite{FriedmanGitmanKanovei2019} that $\ZFCminus$ does not prove the reflection principle.
	
	\subsection{Intuitionistic set theory \texorpdfstring{$\IZF$}{IZF} and Constructive set theory \texorpdfstring{$\CZF$}{CZF}}
	There are two possible constructive formulations of $\ZF$, namely $\IZF$ and $\CZF$, although we will focus on the latter.
	
	% IZF
	$\IZF$ appeared first in H. Friedman's paper \cite{Friedman1973} on the double-negation of set theory. Friedman introduced $\IZF$ as an intuitionistic counterpart of $\ZF$ and showed that there is a double-negation translation from $\ZF$ to $\IZF$, analogous to that from $\mathsf{PA}$ to $\mathsf{HA}$.
	
	\begin{definition}
		$\IZF$ is the theory that comprises the following axioms: Extensionality, Pairing, Union, Infinity, Set Induction, Separation, Collection, and Power Set.
	\end{definition}

    \begin{remark} \label{InfinityDefinition}
    We take the axiom of Infinity to be the statement $\exists a (\exists x (x \in a) \land \forall x \in a \, \exists y \in a \, (x \in y))$. See Remark \ref{IKPformulation} for alternative, equivalent, ways to define Infinity.    
    \end{remark}
    
	% Axioms of CZF
	Constructive Zermelo-Fraenkel set theory, $\CZF$, is introduced by Aczel \cite{Aczel1978} with his type-theoretic interpretation of $\CZF$. We will introduce subtheories called \emph{Basic Constructive Set Theory}, $\mathsf{BCST}$, and $\CZFminus$ before defining the full $\CZF$.
	
	\begin{definition}
		$\mathsf{BCST}$ is the theory that consists of Extensionality, Pairing, Union, Emptyset, Replacement, and $\Delta_0$-Separation.		
		$\CZFminus$ is obtained by adding the following axioms to $\mathsf{BCST}$: Infinity, Set Induction, and \emph{Strong Collection} that states the following: if $\phi(x,y)$ is a formula such that for given $a$, if $\forall x\in a\exists y \phi(x,y)$, then we can find $b$ such that
		\begin{equation*}
			\forall x\in a\exists y \in b \phi(x,y) \land \forall y\in b\exists x\in a \phi(x,y).
		\end{equation*}
	\end{definition}
	
	We also provide notation for frequently-mentioned axioms:
	\begin{definition}
		We will use $\mathsf{Sep}$, $\Delta_0\mhyphen\mathsf{Sep}$, $\Delta_0\mhyphen\mathsf{LEM}$ for denoting Full Separation (i.e., Separation for all formulas), $\Delta_0$-Separation and the law of excluded middle for $\Delta_0$-formulas respectively.
	\end{definition}
	
	% Binary Intersection
	The combination of Full Separation and Collection proves Strong Collection, but the implication does not hold if we weaken Full Separation to $\Delta_0$-Separation. It is also known that $\Delta_0$-Separation is equivalent to the existence of the intersection of two sets. See Section 9.5 of \cite{AczelRathjen2010} for its proof.
	
	\begin{proposition}\pushQED{\qed} 
		Working over $\mathsf{BCST}$ without $\Delta_0$-Separation, $\Delta_0$-Separation is equivalent to the \emph{Axiom of Binary Intersection}, which asserts that $a\cap b$ exists if $a$ and $b$ are sets. \qedhere
		\popQED
	\end{proposition}
	
	% Multi-valued function
	It is convenient to introduce the notion of \emph{multi-valued function} to describe the Strong Collection and Subset Collection axioms that we will discuss shortly. Let $A$ and $B$ be classes. A relation $R\subseteq A\times B$ is a \emph{multi-valued function from $A$ to $B$} if $\dom R=A$. In this case, we write $R\colon A\rrarrows B$. We use the notation $R \colon A\lrlrarrows B$ if both $R \colon A\rrarrows B$ and $R \colon B\rrarrows A$ hold. The reader is kept in mind that the previous definition must be rephrased in an appropriate first-order form if one of $A$, $B$, or $R$ is a (definable) proper class in the same way as how we translate classes over $\ZF$.
	        
	Then we can rephrase Strong Collection as follows: for every set $a$ and a class multi-valued function \mbox{$R \colon a\rrarrows V$}, there is an `image' $b$ of $a$ under $R$, that is, a set $b$ such that $R \colon a\lrlrarrows b$.
	
	% Subset Collection and Fullness
	Now we can state the Axiom of Subset Collection:
	\begin{definition}
		The Axiom of Subset Collection states the following: Assume that $\phi(x,y,u)$ is a formula that defines a collection of multi-valued functions from $a$ to $b$ parametrized by $u\in V$: that is, $\phi(x,y,u)$ satisfies $\forall u \forall x\in a\exists y\in b \phi(x,y,u)$.
		Then we can find a set $c$ such that
		\begin{equation*}
			\forall u \exists d\in c [\forall x\in a\exists y\in d \phi(x,y,u) \land \forall y\in d\exists x\in a \phi(x,y,u)].
		\end{equation*}
		
		$\CZF$ is the theory obtained by adding Subset Collection to $\CZFminus$.
	\end{definition}
	We may state Subset Collection informally as follows: for every first-order definable collection of multi-valued class functions $\langle R_u \colon a\rrarrows b \mid u\in V\rangle$ from $a$ to $b$, we can find a set $c$ of all `images' of $a$ under some $R_u$. That is, for every $u\in V$ there is $d\in c$ such that $R_u \colon a \lrlrarrows d$.
	
	There is a simpler axiom equivalent to Subset Collection, known as \emph{Fullness}, which is a bit easier to understand.
	
	\begin{definition}
		The Axiom of Fullness states the following: Let $\mv(a,b)$ be the class of all multi-valued functions from $a$ to $b$. Then there is a subset $c\subseteq \mv(a,b)$ such that if $r\in \mv(a,b)$, then there is $s\in c$ such that $s\subseteq r$.
		Such a $c$ is said to be \emph{full} in $\mv(a,b)$.
	\end{definition}
	
	Then the following hold:
	\begin{proposition}[\cite{AczelRathjen2001}, \cite{AczelRathjen2010}]\label{Proposition:SubsetCollection}
		\leavevmode
		\begin{enumerate}
			\item\label{Item:Fullness} \normalfont{($\CZFminus$)} Subset Collection is equivalent to Fullness.
			\item \normalfont{($\CZFminus$)} Power Set implies Subset Collection.
			\item \normalfont{($\CZFminus$)} Subset Collection proves the function set ${^a}b$ exists for all $a$ and $b$.
			\item \normalfont{($\CZFminus$)} If $\Delta_0\mhyphen \mathsf{LEM}$ holds, then Subset Collection implies Power Set.
		\end{enumerate}
	\end{proposition}
	
	We will not provide a proof for the above proposition, but the reader may consult with \cite{AczelRathjen2001} or \cite{AczelRathjen2010} for its proof. We also note here that \cite{Rathjen1993} showed that Subset Collection does not increase the proof-theoretic strength of $\CZFminus$ while \cite{Rathjen2012Power} showed that the Axiom of Power Set does.
	
	The following lemma is useful to establish \eqref{Item:Fullness} of \autoref{Proposition:SubsetCollection}, and is also useful to treat multi-valued functions in general:
	
    \begin{lemma}\label{Lemma:PrelimAdjuectmentFtn}
        Let $R \colon A\rightrightarrows B$ be a multi-valued function. Define 
		$\mathcal{A}(R) \colon A\rightrightarrows A\times B$ by
		\begin{equation*}
			\mathcal{A}(R) = \{\lag a,\lag a,b\rag\rag \mid 
			\lag a,b\rag \in R\}.
		\end{equation*}
		%For $S \subseteq A \times B$, we have the following:
        For $S\subseteq A\times B$,
        let $\mathcal{A}^S(R) = \{ \langle a, \langle a, b \rangle \rangle | \langle a, b \rangle \in R \cap S\}$.\footnote{After this lemma we will often abuse notation by referring to $\mathcal{A}^S(R)$ simply as $\mathcal{A}(R)$.}  Then
		\begin{enumerate}
			\item ${\mathcal{A}^S(R)} \colon A\rightrightarrows S \iff R\cap S \colon A\rightrightarrows B$,
			\item ${\mathcal{A}^S(R)} \colon A\leftleftarrows S\iff S\subseteq R$.
		\end{enumerate}
	\end{lemma}
	
	\begin{proof}
		For the first statement, observe that $\mathcal{A}^S(R) \colon A\rightrightarrows S$ is equivalent to
		\begin{equation*}
			\forall a\in A \exists s\in S [\lag a,s\rag \in \mathcal{A}^S(R)].
		\end{equation*}
		By the definition of $\mathcal{A}^S$, this is equivalent to
		\begin{equation*}\label{Formula:Mvaluedftn_eq00}
			\forall a\in A \exists s\in S [\exists b\in B ( s=\lag a,b\rag \land \lag a,b\rag\in R\cap S)].
		\end{equation*}
		We can see that the above statement is equivalent to $\forall a\in A\exists b\in B [\lag a,b\rag\in R\cap S]$, which is the definition of $R\cap S \colon A\rightrightarrows B$.
		For the second claim, observe that $\mathcal{A}^S(R) \colon A\leftleftarrows S$ is equivalent to
		\begin{equation*}
			\forall s\in S \exists a\in A [\lag a,s\rag \in \mathcal{A}^S(R)].
		\end{equation*}
		By rewriting $\mathcal{A}^S$ to its definition, we have
		\begin{equation*}
			\forall s\in S \exists a\in A [\exists b\in B ( s=\lag a,b\rag \in R\cap S)].
		\end{equation*}
		We can see that it is equivalent to $S\subseteq R$.
	\end{proof}

	The following lemma is useful when we work with multi-valued functions because it allows us to replace class multi-valued functions over $A$ with set multi-valued functions in $A$:
	
	\begin{lemma}\label{Lemma:SetMV}
	Assume that $A$ satisfies second-order Strong Collection, that is, for every $a\in A$ and \mbox{$R \colon a\rrarrows A$,} we have $b\in A$ such that $R \colon a\lrlrarrows b$.\footnote{If $A$ is also transitive then {$A$ shall be called} \emph{regular}, and this will be formally defined in \autoref{Definition:RegularityandInaccessibility}.}
	If $a\in A$ and $R \colon a\rrarrows A$, then there is a set $c\in A$ such that $c\subseteq R$ and $c \colon a\rrarrows A$. 
	\end{lemma}
	
	\begin{proof}
	Consider $\mathcal{A}(R) \colon a\rrarrows a\times A$. By second-order Strong Collection over $A$, there is $c\in A$ such that $\mathcal{A}(R) \colon a\lrlrarrows c$. Hence by \autoref{Lemma:PrelimAdjuectmentFtn}, we have $c\subseteq R$ and $c \colon a\rrarrows A$.
	\end{proof}

	It is known that every theorem of $\CZF$ is also provable in $\IZF$. Moreover, $\IZF$ is quite strong in the sense that its proof-theoretic strength is the same as that of $\ZF$. On the other hand, it is known that the proof-theoretic strength of $\CZF$ is equal to that of Kripke-Platek set theory $\KP$.
	$\IZF$ is deemed to be \emph{impredicative} due to the presence of Full Separation and Power Set.\footnote{There is no consensus on the definition of predicativity. The usual informal description of predicativity is rejecting self-referencing definitions.} On the other hand, $\CZF$ is viewed as predicative since it allows for a \emph{type-theoretic interpretation} such as the one given by Aczel, \cite{Aczel1978}. However, adding the full law of excluded middle to $\IZF$ or $\CZF$ results in the same theory, namely $\ZF$.
	
	\subsection{Kripke-Platek set theory}
	Kripke-Platek set theory is the natural intermediate theory between arithmetic and stronger set theories like $\ZF$. $\KP$ has a natural intuitionistic counterpart called Intuitionistic Kripke-Platek, which we denote by $\IKP$.
	
	\begin{definition}
		$\IKP$ is the theory consisting of Extensionality, Pairing, Union, Infinity, Set Induction, $\Delta_0$-Collection, and $\Delta_0$-Separation. 
	\end{definition}
	
	\begin{remark}\label{IKPformulation}
	    The reader is reminded that there are different formulations of $\KP$ and $\IKP$.
	   
	    \begin{enumerate}
	        \item Some authors such as in \cite{Avigad2000} restricts Set Induction in $\KP$ and $\IKP$ to $\Pi_1$-formulas. We include Full Set Induction in $\KP$ and $\IKP$.
	        
	        \item The formulation of Infinity over $\IKP$ is more subtle. Some authors such as \cite{Avigad2000} exclude Infinity from $\KP$ and $\IKP$, and denote $\KP$ with Infinity $\mathsf{KP\omega}$.
	        We also have an apparently stronger formulation, namely Strong Infinity, which is defined as follows: let $\Ind(a)$ be the formula $\varnothing \in a\land \forall x\in a (x\cup\{x\}\in a)$.
        	Then Strong Infinity is the statement
        	\begin{equation*}
                \exists a [\Ind (a) \land \forall b[\Ind(b)\to a\subseteq b]].
        	\end{equation*}
	        Moreover, Lubarsky \cite{Lubarsky2002IKP} uses another alternative formulation of Infinity stated as follows:
	        \begin{equation*}
	            \exists a [\Ind(a)\land \forall x\in a[x=0\lor \exists y\in a (x=y\cup\{y\})]].
	        \end{equation*}
	        However, these formulations are all equivalent over $\IKP$. The equivalence of Strong Infinity and Lubarsky's Infinity is easy to prove. The harder part is proving Strong Infinity from Infinity. This is done over $\CZF$ in Proposition 4.7 of \cite{AczelRathjen2001}, and one can verify that the proof also works over $\IKP$.
        \end{enumerate}
	\end{remark}

    Finally, let us observe that $\IKP$ proves Collection for a broader class of formulas named \emph{$\Sigma$-formulas}:
    \begin{definition} \label{SigmaFormulasDefinition}
        The collection of $\Sigma$-formulas is the least collection which contains the $\Delta_0$-formulas and is closed under conjunction, disjunction, bounded quantifications, and unbounded $\exists$. 
    \end{definition}
    
    \begin{theorem} \label{Theorem:IKPSigmaCollection}
        For every $\Sigma$-formula $\varphi(x, y, u)$ the following is a theorem of $\IKP$: For every set $u$, if $\forall x \in a \exists y \varphi(x, y, u)$ then there is a set $b$ such that
        \[
        \forall x \in a \exists y \in b \varphi(x, y, u) \land \forall y \in b \exists x \in a \varphi(x, y, u).
        \]
    \end{theorem}
    
    We refer the reader to Section 19 of \cite{AczelRathjen2010} or Section 11 of \cite{AczelRathjen2001} for some of the basic axiomatic consequences of $\IKP$ and their proofs.
	
	\subsection{Inductive definition}
	% Inductive definition
	Various recursive constructions on $\CZF$ are given by inductive definitions. The reader might refer to \cite{AczelRathjen2001} or \cite{AczelRathjen2010} to see general information about the inductive definition, but we will review some of the details for the reader who are not familiar with it.
	
	\begin{definition} \label{InductiveDefinition}
		An \emph{inductive definition} $\Phi$ is a class of pairs $\lag X,a\rag$.
		To any inductive definition $\Phi$, associate the operator $\Gamma_\Phi(C)=\{a\mid \exists X\subseteq C \lag X,a\rag\in\Phi\}$. A class $C$ is $\Phi$-closed if $\Gamma_\Phi(C)\subseteq C$.
	\end{definition}
	
	We may think of $\Phi$ as a generalization of a deductive system, and $\Gamma_\Phi(C)$ as a class of theorems derivable from the class of axioms $C$. Some authors use the notation $X\vdash_\Phi a$ or $X/a\in\Phi$ instead of $\lag X,a\rag\in \Phi$.
    The following theorem says each inductive definition induces a least fixed point:
	
	\begin{theorem}[Class Inductive Definition Theorem, $\CZFminus$]
		Let $\Phi$ be an inductive definition. Then there is the smallest $\Phi$-closed class $I(\Phi)$.
	\end{theorem}
	
	The following lemma is the essential tool for the proof of the Class Inductive Definition Theorem.
	See Lemma 12.1.2 of \cite{AczelRathjen2010} for its proof:
	\begin{lemma}[$\CZFminus$]\label{Lemma:ItreationClass}
		Every inductive definition $\Phi$ has a corresponding \emph{iteration class} $J$, which satisfies \mbox{$J^a=\Gamma\left(\bigcup_{x\in a}J^x\right)$} for all $a$, where $J^a=\{x\mid \lag a,x\rag\in J\}$.
	\end{lemma}
	
	\subsection{Constructive \texorpdfstring{$L$}{L}}
	In this subsection, we will define the constructible universe $L$ and discuss its properties over $\IKP$.
	
	Constructing $L$ under a constructive manner was first studied by Lubarsky \cite{Lubarsky1993L}. Lubarsky developed the properties of $L$ over $\IZF$, and showed that $\IZF$ proves $L$ satisfies $\IZF$ plus $V=L$.
	Crosilla \cite{CrosillaPhD2000} showed that the construction of $L$ carries over to $\IKP$\footnote{Crosilla uses $\mathsf{CZF^r}$ to denote what we called $\IKP$}. 
	
	There are at least two ways of defining $L$: the first is using definability over a set model, which Lubarsky \cite{Lubarsky1993L} and Crosilla \cite{CrosillaPhD2000} had taken. Another one is using \emph{fundamental operations}, also called G\"odel operations in the classical context, which was taken by the second author in \cite{MatthwesPhD}. We will follow the second method.
	
	\begin{definition}[Fundamental operations]
	    Define
	    \begin{itemize}
	        \item $1^{st}(x)=a$ iff $\exists u\in x\exists b\in u(x=\langle a,b\rangle)$,
	        \item $2^{nd}(x)=b$ iff $\exists u\in x\exists a\in u(x=\langle a,b\rangle)$,
	        \item $y^"\{z\} := \{u\mid \langle z,u\rangle \in y\}$.
	    \end{itemize}
	    
	    Then we define the fundamental operations as follows:
	    \begin{itemize}
	        \item $\mathcal{F}_p(x,y):=\{x,y\}$,
	        \item $\mathcal{F}_\cap(x,y):=x\cap \bigcap y$,
	        \item $\mathcal{F}_\cup(x):=\bigcup x$,
	        \item $\mathcal{F}_\setminus(x,y) : =x\setminus y$,
	        \item $\mathcal{F}_\times(x,y) :=x\times y$,
	        \item $\mathcal{F}_\to(x, y) = x \cap \{ z \mid \text{$y$ is an ordered pair and } (z \in 1^{st}(y) \rightarrow z \in 2^{nd}(y) ) \}$,
	        \item $\mathcal{F}_\forall(x,y):=\{x^"\{z\}\mid z\in y\}$,
	        \item $\mathcal{F}_d(x,y):=\dom x$
	        \item $\mathcal{F}_r(x,y):=\ran x$
	        \item $\mathcal{F}_{123}(x,y):=\{\langle u,v,w\rangle\mid \langle u,v\rangle\in x\land w\in y\}$,
	        \item $\mathcal{F}_{132}(x,y):=\{\langle u,w,v\rangle\mid \langle u,v\rangle\in x\land w\in y\}$,
	        \item $\mathcal{F}_=(x,y):=\{\langle v,u\rangle\in y\times x\mid u=v\}$,
	        \item $\mathcal{F}_\in(x,y):=\{\langle v,u\rangle\in y\times x\mid u\in v\}$.
	    \end{itemize}
	    For simplicity, we shall let $\mathcal{I}$ be the set of all indices $i$ of $\mathcal{F}_i$ presented in the above definition.
	\end{definition}
	
	The following lemma says that we can represent every $\Delta_0$-formula in terms of fundamental operations:
	\begin{lemma}[\cite{MatthwesPhD}, Lemma 5.2.4, $\IKP$]
	    \pushQED{\qed}
	    Let $\phi(x_1,\cdots, x_n)$ be a bounded formula whose free variables are all expressed. Then there is a term $\mathcal{F}_\phi$ built up from fundamental operations such that
	    \begin{equation*}
	        \IKP\vdash \mathcal{F}_\phi(a_1,\cdots,a_n) = \{\langle x_n,\cdots,x_1\rangle\in a_n\times\cdots\times a_1\mid\phi(x_1,\cdots,x_n)\}. \qedhere 
	    \end{equation*}
	\end{lemma}
	
	\begin{definition}
	    For a set $a$, define
	    \begin{itemize}
	        \item $\mathcal{E}(a) := a\cup \{\mathcal{F}_i(\vec{x})\mid \vec{x}\in a\land i\in\mathcal{I}\}$,
	        \item $\mathcal{D}(a):=\mathcal{E}(a\cup\{a\})$, and
	        \item $\Def(a) := \bigcup_{n\in\omega} \mathcal{D}^n(a)$.
	    \end{itemize}
	\end{definition}
	
    \begin{definition}
	    For an ordinal $\alpha$, define $L_\alpha :=\bigcup_{\beta\in\alpha}\Def(L_\beta)$ and $L:=\bigcup_{\alpha\in\Ord} L_\alpha$.
	\end{definition}
	
	Then we have the following properties of the constructible hierarchy: 
	\begin{proposition}[{\cite[Proposition 5.3.12]{MatthwesPhD}}, $\IKP$] \label{Proposition:propertiesofL}
	\pushQED{\qed}
	    For all ordinals $\alpha,\beta$,
	    \begin{enumerate}
	        \item If $\beta\in\alpha$, then $L_\beta\subseteq L_\alpha$,
	        \item $L_\alpha\in L_{\alpha+1}$,
	        \item $L_\alpha$ is transitive, and
	        \item\label{item:LalphamodelsboundedSep} $L_\alpha$ is a model of Bounded Separation. \qedhere
	    \end{enumerate}
	\end{proposition}
	
	Moreover we can see that $\IKP$ proves $L$ is a model of $\IKP$:
	\begin{theorem}[{\cite[Theorem 5.3.6 and 5.3.7]{MatthwesPhD}}, $\IKP$] \pushQED{\qed}
	   $\IKP^L$ holds. That is, if $\sigma$ is a theorem of $\IKP$, then $\IKP$ proves $\sigma^L$. Furthermore, $L$ thinks $V=L$ holds.
	   \qedhere 
	\end{theorem}
	
	We will not examine its proof in detail, but it is still worthwhile to mention relevant notions that are necessary for the proof. 
	One of these is \emph{hereditary addition}, which was first formulated by Lubarsky \cite{Lubarsky1993L}. This is necessary because $\alpha\in\beta$ does not entail $\alpha+1\in\beta+1$ constructively.
	
	\begin{definition} \label{Definition:HereditaryAddition}
	    For ordinals $\alpha$ and $\gamma$, define $\alpha+_H\gamma$ by recursion on $\alpha$:
	    \begin{equation*}
	        \alpha +_H \gamma := \left(\bigcup \{\beta +_H\gamma \mid\beta\in\alpha\}\cup\{\alpha\}\right) + \gamma.
	    \end{equation*}
	\end{definition}
	
	Another relevant notion is \emph{augmented ordinals}. Lubarsky \cite{Lubarsky1993L} introduced this notion to develop properties of the constructible hierarchy over $\IZF$. Augmented ordinals are not needed to verify that $L$ is a model of $\IKP$, but are used to prove the Axiom of Constructibility, and we will use them when working with elementary embeddings.
	
	\begin{definition} \label{Definition:AugmentedOrdinal}
	    Let $\alpha$ be an ordinal. Then $\alpha$ \emph{augmented}, $\alpha^\#$, is defined recursively on $\alpha$ as
	    \begin{equation*}
	        \alpha^\# := \bigcup\{\beta^\#\mid\beta\in\alpha\}\cup (\omega+1).
	    \end{equation*}
	\end{definition}
	
	\begin{remark}
	    In \cite{MatthwesPhD}, the second author distinguished Strong Infinity and Infinity, and take a more cautious way to define $L$: first, one defines a subsidiary hierarchy $\mathbb{L}_\alpha$ by $\mathbb{L}_\alpha:=\bigcup_{\beta\in\alpha} \mathcal{D}(\mathbb{L}_\beta)$. Then define $\mathbb{L}:=\bigcup_{\alpha\in\Ord}\mathbb{L}_\alpha$.
	    Finally, it is shown that if Strong Infinity holds, then $\mathbb{L}$ satisfies $\IKP$ and $\mathbb{L}=L$.
	    
	    The reason for taking this method is that the equivalence of Infinity and Strong Infinity over $\IKP$ without Infinity is quite non-trivial. This way of defining $L$ also has the benefit that it works even in the absence of any form of Infinity.
	\end{remark}
	
	\begin{remark}
	Unlike either $\IKP$ or $\IZF$, $\CZF$ does not prove that $L$ satisfies $\CZF$, in particular $\CZF$ cannot prove that Exponentiation holds in $L$. Crosilla \cite{CrosillaPhD2000} showed that $\CZF$ proves $L$ validates $\IKP$ with full Collection. Full details of the above construction using fundamental operations, alongside a full investigation of which axioms of $\CZF$ hold in $L$ was undertaken by the second author and Michael Rathjen in \cite{MatthewsRathjen22}.
	\end{remark}
	
	\subsection{Second-order set theories}
	We need a second-order formulation of constructive set theories in order to define large set axioms corresponding to large cardinals beyond choice. We will formulate constructive analogues of G\"odel-Bernays set theory $\mathsf{GB}$ whose first-order counterparts are $\CZF$ and $\IZF$, respectively.
	We will use $\forall^0$ and $\exists^0$ for quantifications over sets, and $\forall^1$ and $\exists^1$ for quantifications over classes. We omit the superscript if the context is clear. Following standard conventions, we will use uppercase letters for classes and lowercase letters for sets, unless specified otherwise.
	The reader could consult with \cite{WilliamsPhD} if they are interested in classical second-order set theory.

	\begin{definition}\label{Definition:CGB}
		\emph{Constructive G\"odel-Bernays set theory} ($\CGB$) is defined as follows: the language of $\CGB$ is two-sorted, that is, $\CGB$ has sets and classes as its objects. $\CGB$ comprises the following axioms:
		\begin{itemize}
			\item Axioms of $\CZF$ for sets.
			\item Every set is a class, and every element of a class is a set. 
			\item Class Extensionality: two classes are equal if they have the same set of members. Formally,
			\begin{equation*}
				\forall^1 X,Y [X=Y\lr \forall^0x (x\in X \lr x\in Y)].
			\end{equation*}
			
			\item Elementary Comprehension: if $\phi(x,p,C)$ is a first-order formula with a class parameter $C$, then there is a class $A$ such that $A=\{x\mid \phi(x,p,C)\}$. Formally,
			\begin{equation*}
				\forall^0p\forall^1C\exists^1 A \forall^0x [x\in A\lr \phi(x,p,C)].
			\end{equation*}
			
			\item Class Set Induction: if $A$ is a class, and if we know a set $x$ is a member of $A$ if every element of $x$ belongs to $A$, then $A$ is the class of all sets. Formally,
			\begin{equation*}
				\forall^1A \big[ [\forall^0 x(\forall^0y\in x (y\in A)\to x\in A)]\to \forall^0x (x\in A) \big].
			\end{equation*}
			
			\item Class Strong Collection: if $R$ is a class multi-valued function from a set $a$ to the class of all sets, then there is a set $b$ which is an `image' of $a$ under $R$. Formally,
			\begin{equation*}
				\forall^1 R\forall^0a[R \colon a\rrarrows V \to \exists^0b(R \colon a\lrlrarrows b)].
			\end{equation*}
		\end{itemize}
	\end{definition}
	
	\begin{remark}
	    We do not need to add the Class Subset Collection axiom
    	\begin{equation*}
    		\forall^1 R\forall^0a\forall^0b[[R\subseteq V\times V\times V\land \forall^0u (R\restricts u \colon a\rrarrows b)]\to \exists^0c\forall^0u\exists^0d\in c(R\restricts u \colon a \lrlrarrows d) ] \footnotemark
    	\end{equation*} 
    	\footnotetext{Here $R\restricts u$ is the class $\{\langle x,y\rangle \mid \langle u,x,y\rangle \in R\}$.} 
    	as an axiom of $\CGB$ because it is derivable from the current formulation of $\CGB$: we can prove it from the first-order Fullness and the Class Strong Collection by mimicking the proof of Set Collection from Fullness and Strong Collection.
    	
    	The reader should also note that, unlike classical $\mathsf{GB}$, there is no additional Separation axiom for sets: that is, we do not assume $a\cap A=\{x\in a \mid x\in A\}$ is a set for a given class $A$. In fact, the assumption that $a\cap A$ is always a set implies Full Separation due to Elementary Comprehension.
	\end{remark}
	
	Thus we introduce a new terminology for classes such that $A\cap a$ is always a set:
	\begin{definition}
		A class $A$ is \emph{amenable} if $A\cap a$ is a set for any set $a$.
	\end{definition}
	
	The following lemma shows every class function is amenable:
	\begin{lemma}[$\CGB$]\label{Lemma:FunctionalClass-amenable}
	    Let $F$ be a class function, then $F$ is amenable.
	\end{lemma}
	\begin{proof}
	    It suffices to show that if $F$ is a class function then $F\restricts a$ is a set for each $a\in V$.
	    Consider the class function $\mathcal{A}(F\restricts a):a\to a\times V$, where $\mathcal{A}$ is the operation we defined in \autoref{Lemma:PrelimAdjuectmentFtn}. By Class Strong Collection, we can find a set $b$ such that $\mathcal{A}(F\restricts a) : a\lrlrarrows b$.
	    
	    We claim that $b=F\restricts a$. For $F\restricts a\subseteq b$, if $x\in a$ then we can find $y\in b$ such that $\lag x,y\rag \in \mathcal{A}(F\restricts a)$. Hence $y=\lag x,z\rag \in F\restricts a$ for some $z$. By functionality of $F$, we have $y=\lag x,F(x)\rag$.
	    For the remaining inclusion, if $y\in b$, then there is $x\in a$ such that $\lag x,y\rag \in\mathcal{A}(F\restricts a)$. Now we can prove that $y=\lag x,F(x)\rag$ since $F$ is a class function.
	\end{proof}
	
	Next, we define the second-order variant of $\IZF$ which we denote by $\IGB$. $\IGB$ allows Separation for arbitrary classes as classical $\mathsf{GB}$ does.
	\begin{definition} \label{Definition:IGB}
		\emph{Intuitionistic G\"odel-Bernays set theory} ($\IGB$) is obtained by adding the following axioms to $\CGB$:
		\begin{enumerate}
			\item Axioms of $\IZF$ for sets.
			\item Class Separation: every class is amenable.
		\end{enumerate}
	\end{definition}

	We know that classical $\mathsf{GB}$ is a conservative extension of $\ZF$. We expect the same results to hold for constructive set theories. The following theorem shows it actually holds, however, its proof will require a small amount of proof theory. We borrow ideas of the proof from \cite{2010MOanswer}. 
	
	\begin{proposition}\label{Proposition:Secondordersettheory-conservativity}
		$\CGB$ is conservative over $\CZF$. $\IGB$ is conservative over $\IZF$.
	\end{proposition}
	
	\begin{proof}
		We only provide the proof for the conservativity of $\CGB$ over $\CZF$ since the same argument applies to the conservativity of $\IGB$ over $\IZF$. We will rely on the cut-elimination theorem of intuitionistic predicate logic. The reader who might be unfamiliar with this can consult with \cite{NegrivonPlato2008} or \cite{Arai2020}.
		
		Assume that $\sigma$ is a first-order sentence that is deducible from $\CGB$. Then we have a cut-free derivation of $\sigma$ from a finite set $\Gamma$ of axioms of $\CGB$. It is known that every cut-free derivation satisfies the \emph{subformula property}, that is, every formula appearing in the deduction is a subformula of $\sigma$ or $\Gamma$.
		
		Thus we have a finite set $\{X_0,\cdots X_n\}$ of class variables that appear in the deduction. Now we divide into the following cases:
		\begin{enumerate}
			\item If $X_i$ appears in an instance of Class Comprehension in the deduction, then the deduction contains $x\in X_i\lr \phi(x,p,C)$ for some class-quantifier free $\phi$. Now replace every $x\in X_i$ with $\phi(x,p,C)$, and $X_i=Y$ with $\forall x (\phi(x,p,C)\lr x\in Y)$ in the deduction.
			\item Otherwise, replace $x\in X_i$ with $x=x$, and $X_i=Y$ with $\forall x (x=x\lr x\in Y)$. 
		\end{enumerate}
		Then we can see that the resulting deduction is a first-order proof of $\sigma$. Hence $\CGB\vdash \sigma$.
	\end{proof}

	\subsection{Second-order set theories with infinite connectives}
	We will use infinite conjunctions to circumvent technical issues about defining elementary embeddings. Hence we define an appropriate infinitary logic and the corresponding second-order set theories. The following definition appears in \cite{Negri2021}:
	\begin{definition}[$\mathsf{G3i_\omega}$, \cite{Negri2021}]
		The intuitionistic cut-free first-order sequent calculus $\mathsf{G3i_\omega}$ is defined by the following rules: initial sequents are of the form
		$P,\Gamma\implies \Delta,P$
		and its deduction rules are:
		\renewcommand*{\arraystretch}{2.7}
		\setlength{\tabcolsep}{12pt}
		\begin{longtable}{*2{>{\(\displaystyle} c <{\)}}}
			\frac{A,B,\Gamma\implies \Delta}{A\land B,\Gamma\implies\Delta} \ L\land &
			\frac{\Gamma\implies \Delta, A \qquad \Gamma\implies\Delta, B}{\Gamma\implies \Delta,A\land B} \ R\land \\
			
			\frac{A_k,\bigwedge_{n\in\omega}A_n,\Gamma\implies \Delta}{\bigwedge_{n\in\omega}A_n,\Gamma\implies\Delta} \ \sideset{}{_k}{L\bigwedge} &
			\frac{\{\Gamma\implies \Delta, A_n\mid n>0\}}{\Gamma\implies \Delta,\bigwedge_{n\in\omega}A_n} \ R\bigwedge \\
			
			\frac{A,\Gamma\implies \Delta\qquad B,\Gamma\implies\Delta}{A\lor B,\Gamma\implies\Delta} \ L\lor &
			\frac{\Gamma\implies \Delta,A,B}{\Gamma\implies \Delta,A\lor B} \ R\lor \\
			
			\frac{\{\Gamma,A_n\implies \Delta\mid n<\omega\}}{\Gamma,\bigvee_{n\in\omega}A_n\implies \Delta} \ L\bigvee & 
			\frac{\Gamma\implies\Delta,\bigvee_{n\in\omega}A_n,A_k}{\Gamma\implies \Delta, \bigvee_{n\in\omega}A_n} \ \sideset{}{_k}{R\bigvee} \\
			
			\frac{\Gamma\implies \Delta,A \qquad B,\Gamma\implies\Delta}{A\to B,\Gamma\implies\Delta}\ L\to & 
			\frac{A,\Gamma\implies \Delta, B}{\Gamma\implies \Delta,A\to B}\ R\to \\
			
			\frac{}{\bot,\Gamma\implies\Delta}\ L\bot &\\
			
			\frac{A[t/x], \Gamma \implies \Delta}{\exists x A, \Gamma \implies \Delta}\ L\exists & 
			\frac{\Gamma\implies \Delta, \exists x A, A[t/x]}{\Gamma\implies \Delta, \exists x A}\ R\exists \\
			
			\frac{\forall x A, A[t/x], \Gamma\implies \Delta}{\forall x A, \Gamma \implies \Delta}\ L\forall &
			\frac{\Gamma\implies \Delta, A[y/x]}{\Gamma\implies\Delta,\forall x A}\ R\forall\ \text{($y$ fresh)}
		\end{longtable}
	\end{definition}

    The difference between finitary logic and $\mathsf{G3i}_\omega$ is that the latter allows for infinite conjunctions and disjunctions over countable sets of formulas. This will only be needed in order to formalise the concepts of Super Reinhardt sets and the total Reinhardtness of $V$ in \autoref{Subsection: Largeset-BeyondChoice}. Otherwise, throughout this paper, we will stick to finitary logic.
 
	It is well-known that the usual classical and intuitionistic predicate calculus enjoys the cut-elimination rule. This is the same for $\mathsf{G3i_\omega}$. We state the following proposition without proof. The reader might refer to \cite{Negri2021} for its proof.
	
	% CUT-ELIMINATION
	\begin{proposition}[Cut elimination for $\mathsf{G3i_\omega}$]\label{Proposition:CutEliminationG3iomega}
		$\mathsf{G3i_\omega}$ with the cut rule proves the same sequents as $\mathsf{G3i_\omega}$ itself (that is, without the cut rule.)
	\end{proposition}
	
	{Finally, w}e are ready to define set theories over $\mathsf{G3i_\omega}$:
	
	\begin{definition} \label{Definition:CGBIGBinfty}
		$\CGB_\infty$ and $\IGB_\infty$ are obtained by replacing the underlying logic of $\CGB$ and $\IGB$ with $\mathsf{G3i_\omega}$, respectively.
	\end{definition}
	
	Note that the infinitary logic $\mathsf{G3i_\omega}$ has no role in defining new sets from the given sets: for example, we do not have Separation schemes for infinite conjunctions and disjunctions.
	The only reason we are using $\mathsf{G3i}$ is to enhance the expressive power of set theories as we introduced classes into first-order set theories. The following theorem shows that infinite connectives do not provide new theorems into $\CGB$ or $\IGB$:	
	
	\begin{proposition}
		$\CGB_\infty$ is conservative over $\CGB$ for formulas with no infinite connectives. The analogous result also holds between $\IGB_\infty$ and $\IGB$.
	\end{proposition}
	
	\begin{proof}
		We only prove the first claim since the proof for the second claim is identical.
		We can see from \autoref{Proposition:CutEliminationG3iomega} that the \emph{subformula property for $\mathsf{G3i_\omega}$} holds: that is, every formula in a derivation of $\Gamma\Rightarrow\Delta$ is a subformula of $\Gamma$ or $\Delta$.
		
		Now assume that we have a derivation $\Gamma\implies\sigma$ from a finite set $\Gamma$ of axioms of $\CGB$ and an infinite-connective free formula $\sigma$. By the subformula property, the derivation has no infinite connectives. This means that the entire derivation can be done without rules for infinite connectives.
		Hence $\CGB$ proves $\sigma$.
	\end{proof}

	\section{Small Large set axioms}\label{Section:SmallLargeSet}
	In this and the next section, we will discuss large set axioms, which is an analogue of large cardinal axioms over constructive set theories. Since ordinals over $\CZF$ could be badly behaved (for example, they need not be well-ordered), we focus on the structural properties of given sets to obtain higher infinities over $\CZF$.
	This approach is not unusual even under classical context since many large cardinals over $\ZF$ are characterized and defined by structural properties of $V_\kappa$ or $H_\kappa$. (One example would be the definition of indescribable cardinals.)
	We also compare the relation between large cardinal axioms over well-known theories like $\ZF$ and its large set axiom counterparts.
	
	% Regular sets and Inaccessible sets: explanation
	The first large set notions over $\CZF$ would be \emph{regular sets}. Regular sets appear first in Aczel's paper \cite{Aczel1986} about inductive definitions over $\CZF$. As we will see later, regular sets can `internalize' most inductive constructions, which turns out to be useful in many practical cases.
	We will follow definitions given by \cite{ZieglerPhD}, and will briefly discuss differences in the terminology between different references.
	
	\vbox{
	\begin{definition} \label{Definition:RegularityandInaccessibility}
		A transitive set $A$ is
		\begin{enumerate}
			\item \emph{Regular} if it satisfies second-order Strong Collection: 
			\begin{equation*}
				\forall a\in A\forall R [R \colon a\rrarrows A\to \exists b\in A (R \colon a\lrlrarrows b)].
			\end{equation*}
			\item \emph{Weakly regular} if it satisfies second-order Collection,
			\item \emph{Functionally regular} if it satisfies second-order Replacement,
		\item \emph{$\bigcup$-regular} if it is regular and $\bigcup a\in A$ for all $a\in A$.
			\item \emph{Strongly regular} if it is $\bigcup$-regular and $\prescript{a}{}b\in A$ for all $a, b \in A$, and
			\item \emph{Inaccessible} if $(A,\in)$ is a model of second-order $\CZF$.
		\end{enumerate}
		
		The \emph{Regular Extension Axiom} $\mathsf{REA}$ asserts that every set is contained in some regular set. The \emph{Inaccessible Extension Axiom} $\mathsf{IEA}$ asserts that every set is contained in an inaccessible set.
	\end{definition}
	}
	
	Note that there is an alternative characterization of inaccessible sets:
	\begin{lemma}[{\cite[Corollary 10.27]{AczelRathjen2001}}, $\CZF$] \pushQED{\qed} A regular set $A$ is inaccessible if and only if $A$ satisfies the following conditions:
	    \begin{enumerate}
	        \item $\omega\in A$,
	        \item $\bigcup a\in A$ if $a\in A$,
	        \item $\bigcap a\in A$ if $a\in A$ and $a$ is inhabited, and
	        \item for $a,b\in A$ we can find $c\in A$ such that $c$ is full in $\mv(a,b)$. \qedhere
	    \end{enumerate}
	\end{lemma}
	\begin{corollary}\label{Corollary:Inaccessiblesets-ExpClosed}
	    Inaccessible sets are strongly regular. Especially, if $a,b\in A$ and $A$ is inaccessible, then ${}^{a}b\in A$.
	\end{corollary}
	
	\begin{proof}
	    Let $A$ be an inaccessible set and $a,b\in A$. If $r \colon a\to b$ is a function and $c\in A$ is full in $\mv(a,b)$, then we can find a multi-valued function $s\in c$ of the domain $a$ such that $s\subseteq r$. Hence $r=s\in c\in A$, which proves $r\in A$, and thus we have ${{}^a}b\subseteq A$. 
        To see ${{}^a}b\in A$, since $A$ satisfies fullness, we can find $c\in A$ which is full in $\mv(a,b)\cap A$. Then ${{}^a}b\subseteq \mv(a,b)\cap A$, and so ${{}^a}b=\{f\in c\mid f\colon a\to b\}\in A$ by $\Delta_0$-Separation over $A$.
	\end{proof}
	
	There is no need for the notion `pair-closed regular sets' since every regular set is closed under pairings if it contains 2:
	\begin{lemma}[{\cite[Lemma 11.1.5]{AczelRathjen2010}}, $\CZFminus$]\label{Lemma:Regular2-Pairing}\pushQED{\qed}
		If $A$ is regular and $2\in A$, then $\lag a,b\rag \in A$ for all $a,b\in A$. \qedhere
	\end{lemma}
	$\mathsf{REA}$ has various consequences. For example, $\mathsf{\CZFminus+REA}$ proves Subset Collection. Moreover, it also proves that every \emph{bounded} inductive definition $\Phi$ has a set-sized fixed point $I(\Phi)$. (See \cite{AczelRathjen2001} or \cite{AczelRathjen2010} for details.)
	
	The notion of regular set is quite restrictive, as it does not have Separation axioms, not even for $\Delta_0$-formulas. Thus we have no way to do any internal construction over a regular set. The following notion is a strengthening of regular set, which resolves the issue of internal construction:
	
	\begin{definition}\label{Definition: BCST-regular}
		A regular set $A$ is \emph{BCST-regular} if $A\models \mathsf{BCST}$. Equivalently, $A$ is a regular set satisfying Union, Pairing, Empty set, and Binary Intersection.
	\end{definition}
	
	We do not know if $\CZF$ proves every regular set is BCST-regular, although Lubarsky and Rathjen \cite{LubarskyRathjen2003} proved that the set of all hereditarily countable sets in the Feferman-Levy model is functionally regular but not $\bigcup$-regular.
	It is not even known that the existence of a regular set implies that of a BCST-regular set. However, every inaccessible set is BCST-regular, and every BCST-regular set we work with in this paper will also be inaccessible.
	
	What are regular sets and inaccessible sets in the classical context? The following result illustrates what these sets look like under some well-known classical set theories: 
	
	\begin{proposition} \label{Prop:ClassicalSmallLargeSets} \phantom{a} 
		\begin{enumerate}
			\item \normalfont{($\ZFminus$)} Every $\bigcup$-regular set containing $2$ is a transitive model of second-order $\ZFminus$, $\ZFminus_2$.
			
			\item \normalfont{($\ZFCminus$)} Every $\bigcup$-regular set containing $2$ is of the form $H_\kappa$ for some regular cardinal $\kappa$.
			
			\item \normalfont{($\ZFminus$)} Every inaccessible set is of the form $V_\kappa$ for some inaccessible cardinal $\kappa$.\footnote{Without choice, the various definitions of inaccessibility are no longer equivalent. Therefore, following \cite{HayutKaragila2020}, we define an ordinal $\kappa$ to be inaccessible if $V_\kappa$ is a model of second-order $\ZF$. Equivalently, $\kappa$ is inaccessible if $\kappa$ is a regular cardinal, and for every $\alpha<\kappa$ there is no surjection from $V_\alpha$ to $\kappa$.}
        \end{enumerate}
	\end{proposition}
	\begin{proof} \phantom{a}
		\begin{enumerate}
			\item Let $A$ be a {$\bigcup$}-regular set containing $2$. We know that $A$ satisfies Extensionality, Set Induction, Union, and the second-order Collection. Hence it remains to show that second-order Separation holds.
			
			Let $X\subseteq A$, $a\in A$ and suppose that $X\cap a$ is inhabited. Fix $c$ in this intersection. Now consider the function $f \colon a\to A$ defined by
			\begin{equation*}
				f(x)=\begin{cases}
					x&\text{if }x\in X,\text{ and}\\
					c&\text{otherwise}.
				\end{cases}
			\end{equation*}
			By second-order Strong Collection over $A$, we have $b\in A$ such that $f \colon a\lrlrarrows
			b$. It is easy to see that $b=X\cap a$ holds.
			
			\item Let $A$ be a regular set. Let $\kappa$ be the least ordinal that is not a member of $A$. Then $\kappa$ must be a regular cardinal: if not, there is $\alpha<\kappa$ and a cofinal map $f \colon \alpha\to \kappa$. By transitivity of $A$ and the definition of $\kappa$, we have $\alpha\in A$, so $\kappa\in A$ by the second-order Replacement and Union, a contradiction.
			
			We can see that $\ZFCminus$ proves $H_\kappa=\{x \mid |\TC(x)|<\kappa\}$ is a class model of $\ZFCminus$, and $A$ satisfies the Well-ordering Principle. 
			We can also show that $A\subseteq H_\kappa$ holds: if not, there is a set $x\in A\setminus H_\kappa$. Since $A$ is closed under transitive closures, $A$ contains a set whose cardinality is at least $\kappa$. Now derive a contradiction from second-order Replacement.
			
			We know that $A\cap\Ord= H_\kappa\cap \Ord=\kappa$. By second-order Separation over $A$, $\mathcal{P}(\Ord)\cap A=\mathcal{P}(\Ord)\cap H_\kappa$. Hence we have $A=H_\kappa$: for each $x\in H_\kappa$, we can find $\theta<\kappa$, $R\subseteq \theta\times\theta$, and $X\subseteq\theta$ such that $(\TC(x),\in,x)\cong (\theta,R,X)$. (Here we treat $x$ as a unary relation.) Then $(\theta,R,X)\in A$, so $x\in A$ by Mostowski Collapsing Lemma.
			
			% For inaccessible sets	
			\item If $A$ is inaccessible, then $A$ is closed under the true power set of its elements since second-order Subset Collection implies $A$ is closed under exponentiation $a,b\mapsto {}^ab$. Hence $A$ must be of the form $V_\kappa$ for some $\kappa$. Moreover, $\kappa$ is inaccessible because $V_\kappa=A\models \mathsf{ZF_2}$.
			\qedhere
		\end{enumerate}
	\end{proof}
	
	It is easy to see from results in \cite{LubarskyRathjen2003} that $\ZF$ proves every rank of a regular set is a regular cardinal. However, the complete characterization of regular sets in a classical context is open:
	
	\begin{question}
		Is there a characterization of regular sets over $\ZFC$? How about $\bigcup$-regular sets over $\ZFminus$?
	\end{question}
	
	{While w}e can see that every inaccessible set over $\ZF$ is closed under power sets, there is no reason to believe that inaccessible sets over constructive set theories are closed under power sets, even if the Axiom of Power Set holds.
	Hence we introduce the following notion:
	
	\begin{definition}\label{Definition:PowerInaccessible}
		An inaccessible set $K$ is \emph{power inaccessible} if, for every $a\in K$, every subset of $a$ belongs to $K$. $\mathsf{pIEA}$ is the assertion that every set is an element of some power inaccessible set.
	\end{definition}
	
	Notice that the existence of a power inaccessible set implies the Axiom of Power Set over $\CZF$. This is because if $K$ is power inaccessible, then $\mathcal{P}(1)\subseteq K$, and from this one can see that every set has a power set since there is a bijection between $\mathcal{P}(a)$ and ${^a}\mathcal{P}(1)$ through characteristic functions, see Proposition 10.1.1 of \cite{AczelRathjen2010} for further details.
    We will mention power inaccessible sets only in the context of $\IZF$. Also, the reader should note that \cite{MatthwesPhD} and \cite{FriedmanScedrov1984} use the word `inaccessible sets' to denote power inaccessible sets.
	
	The reader is reminded that the meaning of an inaccessible set varies over the references. On the one hand, \cite{ZieglerPhD} and \cite{AczelRathjen2010} follows our definition of inaccessibility. On the other hand, other references like \cite{Rathjen1998}, \cite{RathjenGrifforPalmgren1998}, \cite{AczelRathjen2001} and \cite{Rathjen2017} use the following definition, which we will call \emph{REA-inaccessibility}:
	
	\begin{definition}\label{Definition:REA-inaccessible}
		A set $I$ is \emph{REA-inaccessible} if $I$ is inaccessible and $I\models\mathsf{REA}$.
	\end{definition}
	
	The disagreement in the terminology may come from the difference between their set-theoretic and proof-theoretic properties. 	
	On the one hand, \autoref{Prop:ClassicalSmallLargeSets} shows that inaccessible sets over $\ZFminus$ are exactly of the form $V_\kappa$ for some inaccessible cardinal $\kappa$. However, an REA-inaccessible set over $\ZFminus$ is not only of the form $V_\kappa$ for some inaccessible $\kappa$, but also satisfies $\mathsf{REA}$.
	Gitik \cite{Gitik1980} proved that $\ZF$ with no regular cardinal other than $\omega$ is consistent if there is a proper class of strongly compact cardinals, and Gitik's construction can carry over $V_\kappa$ for an inaccessible cardinal $\kappa$ which is a limit of strongly compact cardinals while preserving the inaccessibility of $\kappa$. Thus $\ZF$ with the existence of an inaccessible set does not imply there is an REA-inaccessible set.
	
	On the other hand, \cite{Rathjen2003AFA} showed that the theory $\mathsf{CZF+IEA}$ is equiconsistent with $\mathsf{CZF+REA}$.\footnote{\cite{Rathjen2003AFA} defined inaccessible sets as transitive models of $\CZFminus$ and second-order Strong Collection. These two definitions are equivalent since $\mathsf{REA}$ implies Subset Collection.}
	However, $\CZF + \forall x\exists y(x\in y\land \text{$y$ is REA-inaccessible})$ has a stronger consistency strength than that of $\mathsf{CZF+REA}$.
	
	Before finishing this section, let us remark that the consistency strength of small large sets over $\CZF$ is quite weak compared to their counterparts over $\ZF$. The reader should refer to \cite{Rathjen1993}, \cite{Rathjen1998} or \cite{Rathjen1999Realm} for additional accounts for the following results:
	
	\begin{theorem}[Rathjen]\pushQED{\qed} \label{Proposition:Prooftheoreticstrength-smalllargesets}
		Each pair of theories have an equal consistency strength and proves the same $\Pi^0_2$-sentences.
		\begin{enumerate}
			\item $\mathsf{CZF+REA}$ and $\mathsf{KPi}$, the theory $\KP$ with a proper class of admissible ordinals.
			\item $\CZF+\forall x\exists y (x\in y\land\text{$y$ is REA-inaccessible})$ and $\mathsf{KPI}$, the theory $\KP$ with a proper class of recursively inaccessible ordinals. \qedhere \popQED
		\end{enumerate}
	\end{theorem}

	\section{Large Large set axioms}\label{Section:LargeLargeSet}
    There is no reason to refrain from defining larger large sets. Hence we define stronger large set axioms. We have two ways of defining large cardinals above measurable cardinals: elementary embeddings and ultrafilters. On a fundamental level, ultrafilters make use of the law of excluded middle because either a set or its complement must be in the ultrafilter. This means that it is possible to prove that a set is in the ultrafilter by showing that a different set is not in the ultrafilter, an inherently nonconstructive principle. On a more structural note, the ultrafilter characterization does not immediately entail the existence of smaller large set notions. For example, it is well-known that it is possible to have a model of $\ZF$ in which there is a non-principal, countably complete ultrafilter over $\omega_1$. To obtain the consistency of inaccessible cardinals from this, one is then obliged to work in an inner model, which adds additional complexity to the arguments. 

    Unlike ultrafilters, the embedding characterization will give direct, positive results on large sets even when working in a weak constructive setting. For example, we will see in \autoref{Corollary:CriticalSetsModelsIZFIEA} that, over $\IKP$, if $K$ is a transitive set (which satisfies some minor additional assumptions) which is a `critical point' of an elementary embedding $j \colon V\to M$, then $K$ is a model of $\mathsf{IZF+pIEA}$ and more.
    Therefore, we will use the elementary embedding characterizations to access the notions of measurable cardinals and stronger principles.

	% critical set, Reinhardt set
	\subsection{Critical sets and Reinhardt sets}
	
	In this section, we will give the basic tools that we shall need to work with elementary embeddings. When working with some embedding $j \colon V \rightarrow M$ we will not in general be assumed that either $j$ or $M$ is definable in a first-order manner. This means that we cannot state Collection or Separation for formulas involving $j$ or $M$ in a purely first-order way. Therefore, to deal with them, we will need to expand our base theory to accommodate class parameters. This leads to the following convention, whose main purpose is to simplify later notation.
	
	\begin{convention}\label{Convention:T_j}
	    Let $T$ be a `reasonable' first-order set theory over the language $\{\in\}$ such as $\CZF$, $\IKP$, or $\IZF$ and let $A$ be a class predicate. The theory $T_{A}$ is the theory with the same axiom schemes as $T$ in the language expanded to include the predicate $A$. For example, $\CZF_A$ has the axiom schemes of $A$ plus $\Delta_0$-Separation, Set Induction, Strong Collection, and Subset Collection in the expanded language $\{ \in, A \}$.
	\end{convention}
	
	\begin{definition}\label{Definition:ElementaryEmbedding}
		We work over the language $\in$ extended by a unary functional symbol $j$ and a unary predicate symbol $M$. Let $T$ be a `reasonable' set theory such as $\CZF$, $\IKP$ or $\IZF$ and suppose that $V \models T$. We say that $j \colon V \rightarrow M$ is a \emph{(fully) elementary embedding} if it satisfies the following:
		
		\begin{enumerate}
		    \item $j$ is a map from $V$ into $M$; $\forall x M(j(x))$,
			\item $M$ is transitive; $\forall x (M(x) \rightarrow \forall y \in x M(y))$,
			\item $M \models T$, 
			\item\label{Item:Elementarity} (Elementarity) $\forall \vec{x} [\phi(\vec{x})\lr \phi^M(j(\vec{x}))]$ for every $\in$-formula $\phi$, where $\phi^M$ is the relativization of $\phi$ over $M$.
		\end{enumerate}
	
	For a class of formulas $Q$ (usually $Q=\Delta_0$ or $\mySigma$), we call $j$ \emph{$Q$-elementary} if $j$ satisfies the above definition with \eqref{Item:Elementarity} restricted to $\phi$ a $Q$-formula.	
	Finally, we say that $K$ is a \emph{critical point} of $j$ if it is a transitive set that satisfies $K \in j(K)$ and $j(x)=x$ for all $x\in K$.
	\end{definition}
	
	We begin with the weak theory of the Basic Theory of Elementary Embeddings, originally introduced by Corazza in \cite{Corazza2006}. This is the minimal theory required to state that $j \colon V \rightarrow V$ is an elementary embedding with a critical point but without assuming anything about the relevant Separation, Collection, or Set Induction schemes over the extended language. For our purposes, it will in fact be beneficial to weaken this further to $\BTEE_M$ which is the minimal theory stating that $j \colon V \rightarrow M$ is an elementary embedding with a critical point for some $M \subseteq V$.
	
	We shall see in \autoref{Section:StructuralLowerbdd} that $\IKP_j$ with a critical point implies the consistency of $\mathsf{IZF + BTEE}$ and in \autoref{Section:HeytingInterpretation} that the latter is equiconsistent with $\ZF + \BTEE$, which already has a quite high consistency strength.
	
	\begin{definition}\label{Definition:BTEE}
	    We work over the language $\in$ extended by a unary functional symbol $j$ and a unary predicate symbol $M$. Let $\BTEE_M$ be the assumption that there exists some transitive class $M \subseteq V$ and elementary embedding $j\colon V\rightarrow M$ which has a critical point. We drop $M$ and simply use $\BTEE$ when $V=M$ holds, that is, $\forall x M(x)$.
	    In addition, continuing our previous notation, given a class of formulas $Q$ we let $Q$-$\BTEE_M$ denote $\BTEE_M$ with elementarity replaced by $Q$-elementarity.
    \end{definition}
	
	\begin{remark}
	    If $V \models T + \BTEE_M$ (where $T$ is a reasonable set theory) we will always assume that $V \models T_M$, that is we will assume that $M$ is allowed to appear in the axiom schemes of $T$.
	\end{remark}
	
	Classically, $\BTEE$ is a relatively weak principle which is implied by the large cardinal axiom $0^\#$, which is the statement that there is a non-trivial elementary embedding $j \colon L \rightarrow L$. The reason this is weak is because we do not have Separation in the extended language which means that we cannot produce sets of the form $j \restricts a$ for arbitrary sets $a$. 
	
	\begin{convention}
	    For $T$ a reasonable theory, we will slightly abuse the notation $T_{j, M}$ to extend it to be the theory $T$ plus a fully elementary embedding $j \colon V \rightarrow M$ which has a critical point and which satisfies the axiom schemes of $T$ in the language expanded to include predicates for $j$ and $M$. In particular, we will allow $j$ and $M$ to appear in the Separation Schemes of $T$ (for example, $\Delta_0$-Separation if $T = \CZF$), Set Induction and appropriate Collection axioms of $T$ (for example, $\Sigma$-Collection if $T=\IKP$ or Strong Collection and Subset Collection if $T=\CZF$).
	\end{convention}
	
	Over $\IKP_j$, being a critical point of an elementary embedding is sufficient to prove some of the basic properties one would want from the theory of elementary embeddings. For example, one can show that if there is a critical point then there is (at least one) critical point which is an ordinal (see Chapter 9 of \cite{ZieglerPhD} or Chapter 7 of \cite{MatthwesPhD} for details) or that there is a model of (first-order) $\IZF$. 
	
	However, such an assumption is not necessarily the most natural counterpart to the classical notion of non-trivial elementary embeddings. For example, consider $\ZFC$ plus an elementary embedding $j \colon V \rightarrow M$ with critical point $\kappa$. Then $L_\kappa$ can easily be seen to be a critical point as defined above however it is the `wrong' set to work within this context because it is not a model of second-order $\ZFC$. In particular, one can show that $\mathcal{P}(\omega)^L$ is merely a countable set. Instead, what one wants to work with is a critical point which is an inaccessible set, and this is what we shall define next. Note that in weak theories there is no reason that having a critical point will imply the existence of such a set. For example, it is proven in Chapter 10 of \cite{MatthwesPhD} that it is possible to produce a model of $\ZFCminus$ with a non-trivial elementary embedding, while also having that $\mathcal{P}(\omega)$ is a proper class, and in such a model there will be no set sized models of second-order $\ZFCminus$.
	
	Now let us define critical sets and Reinhardt sets in terms of elementary embeddings:
	
	\begin{definition}[$T_{j, M}$]\label{Definition:CriticalReinhardt}
		Let $K$ be a critical point of an elementary embedding $j \colon V\to M$.
		We call $K$ a \emph{critical set} if $K$ is also inaccessible. If $M=V$ and $K$ is transitive and inaccessible, then we call $K$ a \emph{Reinhardt set}.
	\end{definition}
	
	We introduce one final definition which is intermediate between $\BTEE$ and Reinhardt sets; the \emph{Wholeness Axiom}. This was first formulated by Corazza \cite{Corazza2000} and later stratified by Hamkins \cite{Hamkins2001WA} in order to investigate what theory was necessary to produce the Kunen Inconsistency under $\ZFC$.
	
	\begin{definition}\label{Definition:WA}
	    The \emph{Wholeness axiom} $\WA$ is the combination of $\BTEE$ and Separation\textsubscript{$j$}. $\WA_0$ is the combination of $\BTEE$ and $\Delta_0$-Separation\textsubscript{$j$}.
	\end{definition}
	
	We will use the term critical point and critical set simultaneously, so the reader should distinguish the difference between these two terms. (For example, a critical point need not be a critical set unless it is inaccessible.) 
	Note that the definition of critical sets is different from that suggested by Hayut and Karagila \cite{HayutKaragila2020}: they defined a critical cardinal as a critical point of an elementary embedding $j \colon V_{\kappa+1}\to M$ for some transitive set $M$. This is done to ensure that the embedding $j$ is a set and therefore first-order definable.
    We will instead take the more natural, but not obviously first-order definable, definition where the domain is taken to be the universe, which is what Schlutzenberg calls \emph{V-critical} in \cite{Schlutzenberg2020Extenders}.
	Finding and analyzing the $\CZF$-definition of a critical set that is classically equivalent to a critical set in the style of Hayut and Karagila would be good future work.
	Also, \cite{ZieglerPhD} uses the term `measurable sets' to denote critical sets, but we will avoid this term for the following reasons: it does not reflect that the definition is given by an elementary embedding, and it could be confusing with measurable sets in measure theory.
	
	We do not know if every elementary embedding $j \colon V\to M$ over $\CZF$ enjoys being cofinal. Surprisingly, the following lemma shows that $j$ becomes a cofinal map if $M=V$.
	Note that the following lemma heavily uses Subset Collection. See Theorem 9.37 of \cite{ZieglerPhD} for its proof.
	
	\begin{proposition}[Ziegler \cite{ZieglerPhD}, $\CZF_j$]\label{Proposition:CZFReinhardtEmbeddingCofinal}\pushQED{\qed}
		Let $j \colon V\to V$ be a non-trivial elementary embedding. Then $j$ is \emph{cofinal}, that is, we can find $y$ such that $x\in j(y)$ for each $x$. \qedhere
	\end{proposition}
	
	Note that \cite{ZieglerPhD} uses the term \emph{set cofinality} to denote our notion of cofinality. However, we will use the term cofinality to harmonize the terminology with that of \cite{Matthews2020}. Ziegler's theorem on the cofinality of Reinhardt embedding was an early indication of the high consistency strength of $\CZF$ with a Reinhardt set. This is because \cite{Matthews2020} shows that a related weak set theory, namely $\ZFCminus$ with the $\mathsf{DC}_\mu$-scheme for all cardinals $\mu$, already proves there is no cofinal non-trivial elementary embedding $j \colon V\to V$.
	
	The following lemma, due to Ziegler \cite{ZieglerPhD}, has an essential role in the proof of \autoref{Proposition:CZFReinhardtEmbeddingCofinal} and developing properties of large large sets. We note here that the proof will not require $V$ to satisfy any of our axioms in the language expanded with $j$ or for the embedding to be fully elementary and therefore the Lemma is also provable over weaker theories such as $\mathsf{IKP + BTEE}$. Note that, when working without Power Set, $\mathcal{P}(\mathcal{P}(a))$ should be seen as an abbreviation for the class consisting of those sets whose elements are all subsets of $a$.
	
	\begin{lemma}[$\IKP + \Delta_0\mhyphen\BTEE_M$] \label{Lemma:PowerSetPreserving}
		Assume that $j \colon V\to M$ is a $\Delta_0$-elementary embedding and $K$ is a transitive set such that $j(x)=x$ for all $x\in K$. Then we have the following results for each $a\in K$:
		\begin{enumerate}
			\item If $b\subseteq a$, then $j(b)=b$.
			\item If $b\subseteq \mathcal{P}(a)$, then $j(b)=b$.
			\item If $b\subseteq \mathcal{P}(\mathcal{P}(a))$, then $j(b)=b$.
		\end{enumerate}
		Furthermore, if we can apply the induction to $j$-formulas, then we can show that $j(b)=b$ for all $b\subseteq \mathcal{P}^n(a)$ for each $n\in\omega$, where $\mathcal{P}^n(a)$ is the $n\textsuperscript{th}$ application of the power set operator to $a$.
	\end{lemma}
	
	\begin{proof}
		\leavevmode
		\begin{enumerate}
			\item If $x\in b$, then $x\in K$ by transitivity of $K$, so $x=j(x)\in j(b)$. Hence $b\subseteq j(b)$.
			
			Conversely, $x\in j(b)\subseteq j(a)=a$ implies $x\in K$, so $j(x)=x$. Hence $j(x)\in j(b)$, and we have $x\in b$. This shows $j(b)\subseteq b$. 
			
			\item If $x\in b$, then $x\in\mathcal{P}(a)$. Hence $j(x)=x$ by the previous claim, and we have $x=j(x)\in j(b)$. In sum, $b\subseteq j(b)$.
			
			Conversely, suppose $x\in j(b)$. Observe that $\forall t \in b ( t \subseteq a)$ and this is $\Delta_0$. Thus, by elementarity, $x \subseteq j(a) = a$.
   
			Hence by the previous claim, $x=j(x)$. This shows $j(x)\in j(b)$, so $x\in b$. Hence $j(b)\subseteq b$.
		\end{enumerate}		
	The proof of the remaining case is identical, so we omit it.
	\end{proof}
	
	\subsection{Large set axioms beyond choice}
	\label{Subsection: Largeset-BeyondChoice}
    As before, there is no reason why we should stop at Reinhardt sets. While Reinhardt cardinals are incompatible with $\ZFC$, it is known that all proofs of this use the Axiom of Choice in an essential way. This has led to the study of stronger large cardinals in the $\ZF$ context, the most notable example being \cite{BagariaKoellnerWoodin2019}. We will focus on the concept of super Reinhardt sets and totally Reinhardt sets. However, before we do so, we note a technical difficult{ly regarding} truth predicates and first-order definability, which is why we will need to work in an infinitary second-order theory. 
	 
	Using \cite[Theorem 1]{Gaifman1974} or \cite[Proposition 5.1]{Kanamori2008}, in $\ZF$, an elementary map $j \colon V\to M$ is fully elementary if it is $\Sigma_1$-elementary. Unfortunately, in our context, we do not know whether $\Sigma_n$-elementary embeddings are fully elementary even if $n$ is sufficiently large.
	 
	During this subsection, we will fix an enumeration $\langle \phi_{n}(X,\vec{x})\mid n\in\mathbb{N}\rangle$ of all first-order formulas over the language of set theory with one class variable $X$.
	
	\begin{definition}[$\CGB_\infty$]
		Let $A$ be a class.
		A class $j$ is an $A$-\emph{(fully) elementary embedding} from $V$ to $M$ if $j$ is a class function and $\bigwedge_{n\in\mathbb{N}}\forall\vec{x}[\phi_n(A,\vec{x})\lr \phi^M_{n}(A\cap M,j(\vec{x}))]$. 
		We simply call $j$ \emph{elementary} if $A=V$.
	\end{definition}

	\begin{definition}[$\CGB_\infty$]\label{Definition:SuperReinhardt}
		An inaccessible set $K$ is \emph{super Reinhardt} if for every set $a$ there is an elementary embedding $j \colon V\to V$ such that $K$ is a critical point of $j$ and $a\in j(K)$.
		
		More generally, an inaccessible set $K$ is \emph{$A$-super Reinhardt} for an amenable $A$ if for each set $a$ we can find an $A$-elementary embedding $j \colon V\to V$ whose critical point is $K$ such that $a\in j(K)$.
	\end{definition}
	
    The reader should notice that we have formulated critical sets and Reinhardt sets over 
    $T_{j,M}$ where $T$ was some general, unspecified (first-order) theory, whereas the definition of super Reinhardt sets in formulated in the specified second-order theory of $\CGB_\infty$.
    %$\CGB_\infty$ instead of $\CZF_{j,M}$. 
    Some readers could think that formulating criticalness and Reinhardtness over $\CGB_\infty$ would be more natural, as Bagaria-Koellner-Woodin formulated Reinhardt cardinals over $\mathsf{GB}$ instead of $\ZF_j$.
	It will turn out that sticking to the more first-order formulation of criticalness and Reinhardtness is better for the following technical reason. We will take a double-negation interpretation based on Gambino's Heyting-valued interpretation \cite{Gambino2006} over first-order constructive set theories in \autoref{Section:HeytingInterpretation}. It is harder to extend Gambino's interpretation to a second-order set theory than extending it to $\CZF_{j,M}$ because the former requires extending Gambino's forcing to all classes and defining the interpretation of second-order quantifiers, while the latter does not. However, we cannot avoid the full second-order formulation of super Reinhardtness, unlike we did for Reinhardtness, because its definition asks for class many elementary embeddings.

	Our definition of $A$-super Reinhardtness is different from that of Bagaria-Koellner-Woodin \cite{BagariaKoellnerWoodin2019} because they required $j$ to satisfy $j^+[A] :=\bigcup_{x\in V} j(A\cap x)$ is equal to $A$ instead of $A$-elementarity. It turns out by \autoref{Lemma:DeltaA0elementarityandfullelementarity} that our definition is stronger than that of Bagaria-Koellner-Woodin in our constructive context, but they are equivalent in the classical context.
	
	\begin{lemma}[$\CGB$] \label{Lemma:RestrictingAClassOntoaset}
        Let $\phi(X, x_0,\cdots, x_n)$ be a $\Delta^X_0$-formula whose free variables are all expressed. If $A$ is an amenable class and if we take $a = \TC(\{x_0,\cdots,x_n\})$, then $$\phi(A,x_0,\cdots,x_n) \lr \phi(A\cap a,x_0,\cdots, x_n).$$
	\end{lemma}
	
	\begin{proof}
	    The proof proceeds by induction on $\phi$.
	    \begin{itemize}
	        \item Assume that $\phi$ is an atomic formula. The only non-trivial case is $x\in A$, and it is easy to see that $x\in A\lr x\in \TC(\{x\})\cap A$.
	        \item The cases for logical connectives are easy to check.
	        \item Assume that $\phi(A,x_0,\cdots, x_n)$ is \(\forall y\in x_0 \psi(A,y,x_0,\cdots, x_n)\), and assume that
	        \begin{equation*}
	            \psi(A,y,x_0,\cdots, x_n)\lr \psi(A\cap f(y),x_0,\cdots, x_n)
	        \end{equation*}
	        for all $y$, where $f(y)=\TC(\{y,x_0,\cdots, x_n\})$. If $y\in x_0$, then $f(y)=a=\TC(\{x_0,\cdots, x_n\})$. Hence
	        \begin{align*}
	            \forall y\in x_0 \psi(A,y,x_0,\cdots, x_n) \lr \forall y\in x_0 \psi(A\cap a,y,x_0,\cdots, x_n).
	        \end{align*}
	        The case for $\exists y\in x_0\psi(y,x_0,\cdots, x_n)$ is also similar. \qedhere
	    \end{itemize}
	\end{proof}
	
	\begin{remark} \label{Remark:SuperReinhardtsareReinhardts}
	    It is worthwhile to note that if $K$ is a super Reinhardt set then it is also a Reinhardt set, in the sense that it is a critical point for some elementary embedding $j \colon V \rightarrow V$ for which $\lag V, j \rag$ satisfies every axiom of $\CZF_j$. This is not obvious because, while $\CGB_\infty$ satisfies Class Set Induction and Class Strong Collection, it does not include a Class Separation axiom. Therefore, for an arbitrary class $A$, there is no reason why $\Delta_0^A$-Separation should hold. To circumvent this issue let $j \colon V \rightarrow V$ be an elementary embedding such that $K \in j(K)$ (and everything in $K$ is fixed by $j$). We observe that, by \autoref{Lemma:FunctionalClass-amenable}, $j$ is an amenable class and therefore we can apply the above Lemma from which one can conclude that $\Delta_0^j$-Separation holds.
	\end{remark}

	\begin{lemma}[$\CGB_{\infty}$] \label{Lemma:DeltaA0elementarityandfullelementarity}
		Let $A$ be an amenable class and $j \colon V \rightarrow V$ be a class function. 
		\begin{enumerate}
			\item If $j$ is fully elementary, then $j$ is $\Delta^A_0$-elementary if and only if $j^+[A]=A$.
			\item If $\mathsf{LEM}$ holds, then every $\Sigma_1\cup \Delta^A_0$-elementary embedding $j$ is elementary for all $A$-formulas.
		\end{enumerate}
	\end{lemma}
	
	\begin{proof} \phantom{a}
		\begin{enumerate}
			\item Assume that $j$ is elementary for $\Delta^A_0$-formulas. We first show that $j(A\cap x)=A\cap j(x)$: we know that $\forall y [y\in A\cap x \lr y\in A \land y\in x]$. We can view this sentence as a conjunction of two $\Delta^A_0$-formulas, so we have
			\begin{equation*}
				\forall y [y\in j(A\cap x) \lr y\in A \land y\in j(x)].
			\end{equation*}
			This shows the claim we desired. Hence $j^+[A]=\bigcup_{x\in V} j(A\cap x) = \bigcup_{x\in V} A \cap j(x) = A$. (The last equality holds by the cofinality of $j$, \autoref{Proposition:CZFReinhardtEmbeddingCofinal}.)
			
			Conversely, assume that $j^+[A]=A$ holds and let $\phi(X,x_0,\cdots,x_n)$ be a $\Delta^X_0$-formula with a unique class variable $X$, all of whose free variables are displayed.
			By \autoref{Lemma:RestrictingAClassOntoaset}, $a=\TC(\{x_0,\cdots, x_n\})$ satisfies $\phi(A,x_0,\cdots, x_n)\lr \phi(A\cap a,x_0,\cdots, x_n)$.
			Since $j$ is elementary and $A \cap a$ is a set by amenability, we have
			\begin{equation*}
			    \phi(A\cap a,x_0,\cdots,x_n)\lr \phi(j(A\cap a),j(x_0),\cdots, j(x_n)).
			\end{equation*} 
			Note that $\forall y\in j(a)[y\in j^+[A] \lr y\in j(A\cap a)]$ holds. Furthermore, every bounded variable $y$ appearing in $\phi$ is bounded by $a$.
			Hence 
			\begin{equation*}
			    \phi(j(A\cap a),j(x_0),\cdots, j(x_n)) \lr \phi(j^+[A],j(x_0),\cdots, j(x_n)).
			\end{equation*}
			From the assumption $j^+[A]=A$, we finally have $\phi(A,x_0,\cdots, x_n)\lr \phi(A, j(x_0),\cdots, j(x_n))$.
			
			\item $\mathsf{LEM}$ implies our background theory is $\mathsf{GB}$. We follow Kanamori's proof \cite[Proposition 5.1(c)]{Kanamori2008} that every $\Sigma_1$-elementary cofinal embeddings over $\ZF$ is fully elementary.
			
			Assume that we have the $\Sigma_n^A$-elementarity of $j$. Consider the $A$-formula $\exists x\phi(A,x,a)$, where $\phi$ is a $\Pi^A_n$-formula. Then $\phi(A,x,a)$ implies $\phi(A,j(x),j(a))$, and we have $\exists y \phi(A,y,j(a))$.
			
			Conversely, assume that we have $\exists x \phi(A,x,j(a))$. Since $j$ is cofinal by \autoref{Proposition:CZFReinhardtEmbeddingCofinal}, we can find $\alpha$ such that $x\in V_{j(\alpha)}$. Since $j$ is elementary, $V_{j(\alpha)} = j(V_\alpha)$ and hence $\exists x\in j(V_\alpha)\phi(A,x,j(a))$ holds. Note that, by using second-order Collection for formulas using the class $A$, if $\theta$ is $\Sigma_n^A$ then so is $\forall x \in a \theta$ and, since we have a classical background theory, by taking negations, if $\theta$ is $\Pi_n^A$ then so is $\exists x \in a \theta$ and thus the formula $\exists x\in j(V_\alpha)\phi(A,x,j(a))$ is $\Pi^A_n$. So, by $\Sigma^A_n$-elementarity, we have $\exists x\in V_\alpha \phi(A,x,a)$. \qedhere
		\end{enumerate}
	\end{proof}
    
	The upshot of the second claim in the above lemma is that $\mathsf{GB}$ proves 
	that given a fully elementary embedding, it is $A$-elementary if it is $\Delta^A_0$-elementary. Furthermore, the condition $j^+[A]=A$ is equivalent to $\Delta^A_0$-elementarity. Hence \autoref{Lemma:DeltaA0elementarityandfullelementarity} justifies that our definition of super Reinhardtness is a reasonable constructive formulation of that of Bagaria-Koellner-Woodin.
	
	We also note here that we can weaken the assumption in the second clause of \autoref{Lemma:DeltaA0elementarityandfullelementarity}: namely, $\mathsf{GB}$ proves every \emph{cofinal} $\Delta^A_0$-elementary embedding $j\colon V\to V$ is fully elementary for $A$-formulas.
	
	We end this section with one final large set axiom which we will see suffices to bring the constructive theory into proof-theoretic equilibrium with its classical counterpart.
	
	\begin{definition}\label{Definition:TR}
	\emph{$V$ is totally Reinhardt} ($\TR$) is the following statement: for every class $A$, there is an $A$-super Reinhardt set.
	\end{definition}

	\section{How strong are large large sets over \texorpdfstring{$\CZF$}{CZF}?}\label{Section:StructuralLowerbdd}
	
	In this section, we perform the preparatory work needed to derive a lower bound for the consistency strength of various large large set axioms over $\CZF$. We will show that over $\IKP$, a critical point with a moderate property is a model of $\mathsf{IZF+pIEA}$. We will see that we can find a model of a stronger extension of $\IZF$ under a large large set axioms over $\CZF$.
	
	One might ask why we focus on deriving the models of $\IZF$ with large set properties. The reason is that deriving consistency strength from $\CZF$ is difficult: double-negation translation does not behave well over $\CZF$ with large set axioms. On the contrary, $\IZF$ or its extensions go well with double-negation translations.
	
	\subsection{The strength of elementary embeddings over \texorpdfstring{$\BTEE_M$}{BTEE-M}}
	The aim of this subsection is to derive a lower bound for the consistency strength of an elementary embedding $j\colon V\to M$ with a critical point. 
	One amazing consequence of \autoref{Lemma:PowerSetPreserving} is that the existence of a critical set is already quite strong, compared to small large set axioms over $\CZF$. As in the Lemma, we remark here that we do not require $V$ to satisfy any axiom in the expanded language and therefore the proof also works over $\IKP + \Delta_0\mhyphen \BTEE_M$.
	
	\begin{theorem}[$\IKP + \Delta_0\mhyphen \BTEE_M$] \label{Theorem:CriticalPointModelsIZF}
		Let $K$ be a transitive set such that $K\models \Delta_0\mhyphen\mathsf{ Sep}$, $\omega\in K$ and let $j \colon V\to M$ be a $\Delta_0$-elementary embedding whose critical point is $K$.
		Then $K\models \IZF$.
	\end{theorem}
	
	\begin{proof}
		$K$ satisfies Extensionality and Set Induction since $K$ is transitive. We only prove that Power Set, Full Separation, and Collection are valid over $K$ since the validity of other axioms is not hard to prove. The focal fact of the proof is that $j(K)$ is also a transitive model of $\Delta_0$-Separation and $K\in j(K)$.
		\begin{enumerate}
			\item Power Set: let $a\in K$, and define $b$ as $b=\{c\in K \mid c\subseteq a\}$. We can see that $b\in j(K)$ and $b=j(b)$ by \autoref{Lemma:PowerSetPreserving}. Hence $j(b)\in j(K)$, and thus we have $b\in K$.
			
			It remains to show that $K$ thinks $b$ is a power set of $a$. We claim that $K\models \forall x (x\in b\lr x\subseteq a)$, which is equivalent to $\forall x\in K (x\in b\lr x\subseteq a)$, but this is obvious from the definition of $b$.
			
			\item Full Separation: let $a\in K$ and $\phi(x,p)$ be a first-order formula with a parameter $p\in K$. Observe that the relativization $\phi^K(x,p)$ is $\Delta_0$, so $b=\{x\in a\mid \phi^K(x,p)\}\in j(K)$. Since $b\subseteq a$, we have $j(b)=b$. Hence $b\in K$.
			
			It suffices to show that $K$ thinks $b$ witnesses this instance of Separation for $\phi$ and $a$. Formally, it means $\forall x\in K [x\in b\lr x\in a\land \phi^K(x,p)]$ holds, which trivially holds by the definition of $b$.
			
			\item Collection: Since $K$ models Full Separation, by $\Delta_0$-elementarity of $j$, we also have $j(K)$ satisfies Full Separation.
			Now, let $\phi(x,y,p)$ be a formula and $a,p\in K$, and suppose that $\forall x\in a \exists y\in K \phi^K(x,y,p)$ holds. Now, for each $x\in a$ and $y\in K$, since everything is fixed by $j$, we have
			\begin{equation*}
				V\models \phi^K(x,y,p) ~~~ \implies ~~~ M\models \phi^{j(K)} (x,y,p).
			\end{equation*}
			Therefore, $M\models \forall x\in a\exists y\in K \phi^{j(K)}(x,y,p)$. Define $b=\{y\in K \mid \exists x\in a\ \phi^{j(K)}(x,y,p)\}\in j(K)$, then $b$ witnesses $M\models \exists b\in j(K) \forall x\in a \exists y\in b\  \phi^{j(K)}(x,y,p)$. Thus, by elementarity, we have $\exists b\in K\forall x\in a\exists y\in b \ \phi^K(x,y,p)$. Namely, $K$ satisfies this instance of Collection.
			\qedhere 
		\end{enumerate}
	\end{proof}
	
	Hence the existence of an elementary embedding with a critical point is already tremendously strong compared to $\CZF$.
	We will see in \autoref{Section:HeytingInterpretation} that $\IZF$ interprets classical $\ZF$, so the existence of a critical set over $\CZF$ is stronger than $\ZF$. However, it turns out that the hypotheses given in \autoref{Theorem:CriticalPointModelsIZF} proves the critical point $K$ satisfies not only $\IZF$, but also large set axioms.
	For example, we next show that $K$ is a model of $\mathsf{IZF+pIEA}$. The proof begins with the following lemmas:
	
	\begin{lemma}[$\IKP + \Delta_0\mhyphen \BTEE_M$] \label{Lemma:KandjK-samepowerset}
	    Assume that $K$ satisfies the conditions given in the hypotheses of \autoref{Theorem:CriticalPointModelsIZF}. If $a\in K$, then $\mathcal{P}(a)\cap K = \mathcal{P}(a)\cap j(K)$.
	\end{lemma}
	
	\begin{proof}
	    Let $a\in K$, $b\in j(K)$ and $b\subseteq a$. By \autoref{Lemma:PowerSetPreserving}, we have $b\in K$.
	\end{proof}
	
	\begin{lemma}[$\IKP + \Delta_0\mhyphen \BTEE_M$] \label{Lemma:jKThinksKInacc}
	    Under the hypotheses for $K$ given in \autoref{Theorem:CriticalPointModelsIZF}, $j(K)$ thinks $K$ is power inaccessible.
	\end{lemma}
	
	\begin{proof}
	    We proved that $K$ is a model of $\IZF$ and $\mathcal{P}(a) \cap K = \mathcal{P}(a) \cap j(K)$ for any $a \in K$. Thus it suffices to show that $j(K)$ believes $K$ is regular. That is,
	    \begin{equation*}
	        \forall a \in K \forall R \in j(K) [ R\colon a\rrarrows K\to \exists b \in K (R \colon a \lrlrarrows b)].
	    \end{equation*}
	    The proof employs basically the same argument as the proof that $K$ satisfies Collection given in \autoref{Theorem:CriticalPointModelsIZF}:
	    let $a\in K$, $R\in j(K)$, and suppose that $R\colon a\rrarrows K$. Then 
	    \begin{equation*}
	        V\models R: a\rrarrows K ~~~ \implies ~~~ M\models j(R): a\rrarrows j(K).
	    \end{equation*}

		Now let $b=\{y\in K \mid \exists x\in a \ \langle x,y\rangle\in j(R)\}\in M$. Since $b$ is a subset of $K\in j(K)$, we have $b\in j(K)$. Furthermore, $b$ witnesses $M \models \exists b\in j(K) [j(R)\colon a\lrlrarrows b]$.
		Hence, by elementarity, there exists $b \in K$ such that $\forall x\in a\exists y\in b (R\colon a\lrlrarrows b)$.
	\end{proof}
	
	\begin{corollary}[$\IKP + \Delta_0\mhyphen\BTEE_M$] \label{Corollary:CriticalSetsModelsIZFIEA}
		Under the hypotheses for $K$ given in \autoref{Theorem:CriticalPointModelsIZF}, 
		$K$ satisfies $\mathsf{pIEA}$.
	\end{corollary}
	
	\begin{proof}
	    By \autoref{Lemma:KandjK-samepowerset}, $j(K)$ believes $K$ is power inaccessible.
	    Fix $a\in K$, then $b=K$ witnesses the following statement:
	    \begin{equation*}
	        M \models \big[ \exists b\in j(K) \big( a\in b \land j(K)\models \text{$b$ is power inaccessible} \big) \big].
	    \end{equation*}
	    Hence by elementarity, there is $b\in K$ such that $a\in b$ and $K$ believes $b$ is power inaccessible. Since $a$ is arbitrary, we have $K \models \mathsf{IZF + pIEA}$.
	\end{proof}
	
	In fact, we can derive more: we can show that not only does $K$ satisfy $\mathsf{pIEA}$, but $K$ also satisfies \mbox{$\forall x\exists y[x\in y\land \Phi(y)]$}, where $\Phi$ is any large set properties such that $j(K) \models \Phi(K)$. Especially, $\Phi$ can be a combination of power inaccessibility with \emph{Mahloness} or \emph{2-strongness}. (See \cite{Rathjen1998} or \cite{ZieglerPhD} for the details of Mahlo sets or 2-strong sets.)
	
	\subsection{Iterating \texorpdfstring{$j$}{j} to a critical point}
    In the classical context, one usually works with $\lambda=j^\omega(\kappa)=\sup_{n<\omega}j^n(\kappa)$ for $\kappa=\crit j$ when they study large cardinals at the level of rank-into-rank embeddings and beyond.
    We can see in the classical context that if $j\colon V\to M$ is elementary embedding with $\kappa=\crit j$, then $V_\lambda$ is still a model of $\ZFC$ with some large cardinal properties.
    
    We may expect the same in a constructive manner, but defining the analogue of $\lambda=j^\omega(\kappa)$ cannot be done over $\IKP+\mySigma\mhyphen \BTEE_M$. The problem is that the definition of $\lambda$ requires defining $\langle j^n(\kappa) \mid n\in\omega\rangle$, which requires recursion for $j$-formulas. 
    It turns out that we can define this sequence if we have Set Induction for $\mySigma^{j, M}$-formulas. Moreover, even if we have this sequence, we require an instance of $\mySigma^{j, M}$-Replacement to turn the sequence into a set, which is when we will have to start working in $\IKP_{j,M}$.
    
    \begin{definition}[$\IKP+\mySigma\mhyphen \BTEE_M$]\label{Definition:Iterationofj}
        Let $j\colon V\to M$ be an elementary embedding with a critical point $K$. Following Corazza \cite[(2.7) to (2.9)]{Corazza2006}, let $\Theta$ and $\Phi$ be the following formulas:
        \begin{align*}
            \Theta(f,n,x,y) & \equiv \text{$f$ is a function}\land \dom f = n+1 \land f(0) = x \, \land \\ &\forall i \big( 0<i\le n\to f(i) = j(f(i-1)) \land f(n) = y\big)\\
            \Phi(n,x,y) & \equiv n\in\omega \to \exists f\ \Theta (f,n,x,y).
        \end{align*}
    \end{definition}
    
    Informally speaking, $\Theta(f,n,x,y)$ states that $f$ is a function with domain $n + 1$ computing $j^n(x)=y$, and $\Phi(n,x,y)$ asserts that $j^n(x)=y$.
    We can see that $\Theta$ is $\Delta^{j,M}_0$ and $\Phi$ is $\Sigma^{j,M}$. A careful analysis will show that $\Phi$ allows a $\Pi^{j,M}$-formulation, so $\Phi$ is actually $\Delta^{j,M}$. However, this fact is irrelevant in our context.
    
    \begin{lemma}[{Corazza,} $\IKP+\mySigma\mhyphen \BTEE_M+\mySigma^{j,M}\mhyphen\text{Set Induction}$] \label{Lemma:DefinabilityIterationOfj}
        \phantom{a} 
        
        \begin{enumerate}
            \item For all $n\in\omega$ and $x$, $y$, there is at most one $f$ such that $\Theta(f,n,x,y)$ holds. That is,
            \begin{equation*}
                \forall n\in\omega \forall x,y,f,g\ [\Theta(f,n,x,y)\land \Theta(g,n,x,y)\to f=g].
            \end{equation*}
            \item $\Phi$ defines a function. That is, $\forall n\in\omega \forall x \exists! y \Phi(n,x,y)$.
        \end{enumerate}
    \end{lemma}
    
    \begin{proof} \phantom{a}
        \begin{enumerate}
            \item Fix $x,y,f,g$ and assume that both $\Theta(f,n,x,y)$ and $\Theta(g,n,x,y)$ hold. To be precise, we apply Set Induction to the following formula:
            \begin{equation*}
                \theta(i) \equiv [i\le n\to f(i)=g(i)].
            \end{equation*}
            Assume that $\forall j\in i \ \theta(j)$ holds. 
            Then we can see that $f(i) = j(f(i-1))=j(g(i-1))=g(i)$ holds. In other words, we have $\theta(i)$.
            By Set Induction for $\Delta_0$-formulas, we have $\forall i \ \theta(i)$, so $f=g$.
            
            \item For uniqueness, assume that we have $\Phi(n,x,y_0)$ and $\Phi(n,x,y_1)$. By the previous claim, we can see that for the same $f$, $\Theta(f,n,x,y_0)$ and $\Theta(f,n,x,y_1)$ hold. Hence $y_0=f(n)=y_1$.
            
            For existence, we claim that $\forall n\in\omega\forall x \exists y,f\ \Theta(f,n,x,y)$. This uses Set Induction for $\mySigma^{j,M}$-formulas{. F}ix $x$ and let
            \begin{equation*}
                \phi(n, x) \equiv n\in\omega \to \exists y\exists f\ \Theta(f,n,x,y).
            \end{equation*}
            Assume that $\forall i\in n \ \phi(i, x)$ holds.
            Suppose that $y$ and $f$ witnesses $\Theta(f,n-1,x,y)$, then we can extend $f$ to $f'$ by setting $f':=f\cup \{\langle n, j(y) \rangle\}$. It is immediate that $\Theta(f',n,x,y)$. Hence, by Set Induction for $\mySigma^{j,M}$-formulas, we have $\forall n \phi(n,x)$. \qedhere 
        \end{enumerate}
    \end{proof}
    
    \begin{remark}
	    When working with elementary embeddings derived from ultrafilters in $\ZFC$ one can show that the ultrapower construction can be iterated. That is, one starts with $j \colon V \rightarrow M_0$ and shows that for every $n$ there is a transitive class $M_{n+1}$ and an embedding $i_n \colon M_n \rightarrow M_{n + 1}$. Then one can show that $j^2 = i_0 \circ j \colon V \rightarrow M_1$ is an elementary embedding with a critical set. However, this argument does not go through in our weaker context. In particular, it is unclear what the codomain of $j^2$ should be and therefore why $j^2$ should be an elementary embedding. Therefore, for our purposes, in general, $j^n(x)$ should be a formal object which is the result of applying $j$ to $x$ $n$ times. On the other hand, one should note that if $j$ is restricted to some set $A$ then it is possible to iterate the (set) embedding $j \restricts A \colon A \rightarrow j(A)$. 
	\end{remark}
    
    Hence we can define the constructive analogue of $\lambda=j^\omega(\kappa)$, by using an instance of $\mySigma^{j,M}$-Replacement:
    \begin{definition}[$\IKP_{j,M}$]
        Let $j\colon V\to M$ be a $\Delta_0$-elementary embedding with a critical point $K$. By \autoref{Lemma:DefinabilityIterationOfj} and Replacement for $\mySigma^{j,M}$-formulas, $\langle j^n(K)\mid K\in\omega\rangle$ is a set. Define $\Lambda := \bigcup_{n\in\omega} j^n(K)$.
    \end{definition}
    
	\begin{remark}
	    In order to define $\Lambda$ in the above definition it is important that $V$ satisfies at least $\mySigma^{j,M}$-Replacement and $\mySigma^{j,M}$-Induction for $\omega$. One can see that the class function sending $n$ to $j^n(K)$ is $\mySigma^{j,M}$-definable and that $\mySigma^{j,M}$-Induction then implies that this function is total. Finally, an instance of $\mySigma^{j,M}$-Replacement gives us that $\{ j^n(K) \mid n \in \omega \}$ is a set, and therefore so is its union. We refer to the end of Section 2 of \cite{Corazza2006} or Section 6.3 of \cite{MatthwesPhD} for details of this and what issues can potentially arise without assuming any induction in the extended language. 
	\end{remark}
	
	\begin{remark}
	    In the above definition, we keep talking about $\Delta_0$-elementary embedding while we work over $\IKP_{j,M}$, which technically includes the full elementarity of $j$ as an axiom. Hence we actually work with a weaker subtheory of $\IKP_{j,M}$, which is obtained by weakening the elementarity scheme to $\Delta_0$-formulas. This fact becomes important in \autoref{Subsection:Sigma-Ord-inary-IKP} where we work in $\IKP$ with a $\Sigma$-elementary embedding $j\colon V\to M$.
	    
	    However, we will continue to refer to the background theory $\IKP_{j,M}$ when this does not cause confusion. To emphasize, we work with the Set Induction scheme over the extended language as well as the $\Delta_0^{j,M}$-Separation and the $\mySigma^{j,M}$-Collection schemes. 
	\end{remark}
	
	Now we can see that $K$ is an elementary submodel of $\Lambda$:
	\begin{lemma}[$\IKP_{j,M}$]\label{Lemma:KelementarysubmodelLambda}
		Let $j \colon V\to M$ be an $\Delta_0$-elementary embedding with a critical point $K$, and assume that $K\models \Delta_0\mhyphen \mathsf{Sep}$. Then $K$ is an elementary submodel of $\Lambda$. That is, the following holds: for every formula $\phi(\vec{x})$ (without $j$ or $M$) all of whose free variable are expressed, 
		\begin{equation}\label{Formula:PhiKLambdaAbsoluteness}
			V\models \forall \vec{x}\in K [\phi^K(\vec{x})\lr \phi^\Lambda(\vec{x})].
		\end{equation}
	\end{lemma}
	
	\begin{proof}
		We proceed with the proof by induction on formulas. Note that we can formulate our proof over $\IKP$ since the truth predicate for transitive models is $\Sigma$-definable.
		Atomic cases and the cases for $\land$, $\lor$, and $\to$ are trivial, and so we concentrate on cases for quantifiers.
		
		\begin{itemize}
			\item Case $\forall$: assume that $\vec{a}\in K$ and $\phi(x,\vec{a})$ is absolute between $K$ and $\Lambda$. That is, assume that \eqref{Formula:PhiKLambdaAbsoluteness} holds for $\phi$. Obviously \mbox{$V\models (\forall x\phi(x,\vec{a}))^\Lambda\to (\forall x\phi(x,\vec{a}))^K$.}
			
			Conversely, assume that $V\models \forall x\in K \phi^K(x,\vec{a})$. Since this is $\Delta_0$-expressible in $V$, by elementarity we have $M\models \forall x\in j(K) \phi^{j(K)}(x,\vec{a})$.
			Since $\forall x\in j(K) \phi^{j(K)}(x,\vec{a})$ is bounded and $M$ is a transitive subclass of $V$, $\forall x\in j(K) \phi^{j(K)}(x,\vec{a})$ is absolute between $V$ and $M$, so $V\models \forall x\in j(K) \phi^{j(K)}(x,\vec{a})$.
			We can iterate $j$ in a similar fashion, so we have $V\models \forall x\in j^n(K) \phi^{j^n(K)}(x,\vec{a})$ for all $n\in\omega$. Here $\vec{a}$ is unchanged because it is in $K$. 
			Similarly, by elementarity and absoluteness of bounded formulas, we have from \eqref{Formula:PhiKLambdaAbsoluteness} that
			\begin{equation*}
			    V\models \forall n\in\omega \forall x,\vec{a}\in j^n(K) [\phi^{j^n(K)}(x,\vec{a})\lr \phi^\Lambda(x,\vec{a})].
			\end{equation*}
			
			Hence $V$ thinks $\phi^\Lambda(x,\vec{a})$ for all $n\in\omega$ and $\vec{a}\in K$. Since $\Lambda=\bigcup_{n\in\omega} j^n(K)$, we have \mbox{$\forall x\in\Lambda \phi^\Lambda(x,\vec{a})$.}
            
            \item Case $\exists$: The proof is similar to the case for $\forall$. 
			Assume the same conditions to $\vec{a}$ and $\phi(x,\vec{a})$ as we did before. Showing \mbox{$V\models (\exists x\phi(x,\vec{a}))^K\to (\exists x\phi(x,\vec{a}))^\Lambda$} is trivial.
			For the converse, assume that $V\models \exists x\in\Lambda\phi^\Lambda(x,\vec{a})$. Then we can find $n\in\omega$ such that $V\models \exists x\in j^n(K)\phi^\Lambda(x,\vec{a})$. We can prove $\phi^\Lambda(x,\vec{a})$ is equivalent to $\phi^{j^n(K)}(x,\vec{a})$ for $x,\vec{a}\in j^n(K)$, so $V\models \exists x\in j^n(K) \phi^{j^n(K)}(x,\vec{a})$.
			
			Then we have the desired result if $n=0$. If $n>0$, 
			then by absoluteness of bounded formulas, $M\models \exists x\in j^n(K) \phi^{j^n(K)}(x,\vec{a})$. By applying elementarity, we have $V\models \exists x\in j^{n-1}(K) \phi^{j^{n-1}(K)}(x,\vec{a})$. Now we can see by repeating this argument that $V\models \exists x\in K \phi(x,\vec{a})$. \qedhere
			\end{itemize}
	\end{proof}
	
    An upshot of \autoref{Lemma:KelementarysubmodelLambda} is that $j^n(K)$ is an elementary submodel of $\Lambda$.\footnote{In fact, the proof of \autoref{Lemma:KelementarysubmodelLambda} implicitly implies $K$ is an elementary submodel of $j^n(K)$.} Especially, by applying this fact to \autoref{Lemma:jKThinksKInacc}, we have
	
	\begin{lemma}[$\IKP_{j,M}$]\pushQED{\qed}
	\label{Lemma:LambdaThinksKInacc}
		Assume that $K$ satisfies the conditions given in the hypotheses of \autoref{Theorem:CriticalPointModelsIZF}. Then $\Lambda$ thinks $K$ is a power inaccessible set. \qedhere 
	\end{lemma}
	
	Especially, we have an analogue of \autoref{Lemma:KandjK-samepowerset} between $K$ and $\Lambda$:
	\begin{lemma}[$\IKP_{j,M}$]\label{Lemma:KandLambda-samepowerset} \pushQED{\qed}
		Assume that $K$ satisfies the conditions given in the hypotheses of \autoref{Theorem:CriticalPointModelsIZF}. For $n\in\omega$ and $a\in j^n(K)$, we have $\mathcal{P}(a)\cap\Lambda=\mathcal{P}(a)\cap j^n(K)$. \qedhere 
	\end{lemma}
	
	As a corollary of the previous results, we have
	\begin{corollary}[$\IKP_{j,M}$]\label{Corollary:LambdaModelsIZFBTEEInd}
		Assume that $K$ satisfies the conditions given in the hypotheses of \autoref{Theorem:CriticalPointModelsIZF}.\
	    Then $\Lambda$ satisfies $\lag \Lambda,j\restricts\Lambda \rag \models\mathsf{IZF+BTEE}+\text{Set Induction}_j$.
	\end{corollary}
	
	\begin{proof}
	    $\Lambda$ is a model of $\IZF$ by \autoref{Lemma:KelementarysubmodelLambda} and \autoref{Theorem:CriticalPointModelsIZF}. Moreover, $j\restricts \Lambda:\Lambda\to\Lambda$ has a critical point $K\in\Lambda$, so $ \lag \Lambda, j\restricts \Lambda \rag$ is a model of $\BTEE$.
	    Since $\Lambda$ is transitive and Set Induction for $\Delta_0^{j,M}$ formulas holds, we have $\lag \Lambda, j\restricts\Lambda \rag$ believes Set Induction for $j$-formulas holds.
	\end{proof}
	
	\begin{remark}
	    The reader might notice that we always try to reveal where the given formula holds over, like $V\models\phi$ or $M\models\phi$, and explicitly state how to transfer between statements holding over $V$ and over $M$, for example, by relying on absoluteness of bounded formulas.
	    The main reason for it is that while $V\models\phi(K)$ implies $M\models \phi(j(K))$, we cannot say anything about $V\models \phi(j(K))$ and $M\models\phi(K)$.
	    
	    For example, work over $\ZFC$ with a measurable cardinal $\kappa$, and consider an ultrapower map \mbox{$j\colon V\to\operatorname{Ult}(V,U)\cong M$.} Then $V\models (V_\kappa\models \ZFC_2)$ and $M\models (V_{j(\kappa)}\models \ZFC_2)$. However, $j(\kappa)$ is not even a cardinal over $V$, so $V\not\models (V_{j(\kappa)}\models \ZFC_2)$. Similarly, while $V$ thinks $\kappa$ is measurable, $M$ does not in general think $\kappa$ is measurable.
	\end{remark}
	
	\subsection{Playing with a \texorpdfstring{$\Sigma$-$\Ord$-inary}{Sigma-Ord-inary} elementary embedding over \texorpdfstring{$\IKP$}{IKP}}
	\label{Subsection:Sigma-Ord-inary-IKP}
	
	 Let us examine what we can derive from a $\mySigma$-elementary embedding with an ordinal critical point. An illuminating result along this line is that of Ziegler \cite{ZieglerPhD}, which is proved by considering the rank of any fixed critical point $K$:
	
	\begin{proposition}[Ziegler {\cite[Section 9.1]{ZieglerPhD}} or {\cite[Lemma 7.2.4]{MatthwesPhD}}, $\IKP+ \mathrm{\mySigma}\mhyphen \BTEE_M$] \pushQED{\qed}
	
	    Let $j\colon V\to M$ be a $\mySigma$-elementary map.
	    Then the following statements are all equivalent:
	    \begin{enumerate}
	        \item There is $K$ such that $K\in j(K)$ and $\forall x\in \TC(K) \ (j(x)=x)$,
	        \item There is a transitive $K$ such that $K\in j(K)$ and $\forall x\in K\ (j(x)=x)$,
	        \item There is an ordinal $\kappa$ such that $\kappa\in j(\kappa)$ and $\forall\alpha\in\kappa \ (j(\alpha)=\alpha)$. \qedhere
	    \end{enumerate}
	\end{proposition}
	
	The main idea of the above proposition is extracting the rank of a given $K$, then observing that $j$ respects the rank of a set since the rank function is $\mySigma$-definable.
	
	However, it is not likely that we can extend this result further. Ziegler provided a way to give a model of $\mathsf{CZF^-_{Rep}+Exp}$, $\CZF^-$ with Exponentiation and Replacement in place of Strong Collection, from a $\Delta_0$-elementary embedding $j\colon V\to M$ with an ordinal critical point. However, Ziegler's result requires an additional assumption on $j$ called $j\mathsf{IEA}$. 
	
    Despite that, $\mySigma$-elementary embedding over $\IKP$ still shows quite a strong consistency strength provided if it has an ordinal critical point. The following result is from \cite{MatthwesPhD}:
	
	\begin{theorem}[{\cite[Theorem 7.3.2]{MatthwesPhD}}, $\IKP+ \mathrm{\mySigma}\mhyphen \BTEE_M$] \label{Theorem:SigmaOrdEmbedding-LmodelsIZF}
	   Let $j\colon V\to M$ be a $\mySigma$-$\Ord$-inary elementary embedding with witnessing ordinal $\kappa$, that is, a $\mySigma$-elementary embedding with a critical point $\kappa$ that is an ordinal. Furthermore, let $\kappa^\#$ be defined as in \autoref{Definition:AugmentedOrdinal}.
	   Then $L_{\kappa^\#}\models\IZF$.
	\end{theorem}
	
    Using the results we previously developed we can extend this to a model of $\mathsf{IZF +pIEA}.$
    
	\begin{proposition}[$\IKP+\mySigma\mhyphen \BTEE_M$]
	    Let {$j \colon V \rightarrow M$} be a $\mySigma$-$\Ord$-inary elementary embedding with witnessing ordinal $\kappa$.
	    Then $j\restricts L^V : L^V\to L^M$ is a $\mySigma$-elementary embedding over $L$ and $L_{\kappa^\#}\models \Delta_0\mhyphen \mathsf{Sep}$. 
	\end{proposition}
	
	\begin{proof}
	    $\Delta_0$-elementarity of $j\restricts L^V : L^V\to L^M$ follows from the fact that the formula $x\in L$ is $\mySigma$-definable and $j$ is $\mySigma$-elementary.
	    Secondly, by \eqref{item:LalphamodelsboundedSep} of \autoref{Proposition:propertiesofL},
	    $L_{\kappa^\#}\models \Delta_0\mhyphen \mathsf{Sep}$.)
	\end{proof}
	
	By applying \autoref{Theorem:CriticalPointModelsIZF} and \autoref{Corollary:CriticalSetsModelsIZFIEA}, we have
	
	\begin{corollary}[$\IKP+\mySigma\mhyphen \BTEE_M$]
	    Let $j$ be a $\mySigma$-$\Ord$-inary elementary embedding with witnessing ordinal $\kappa$. Then $L_\kappa^\#$ is a model of $\mathsf{IZF+pIEA}$. 
	\end{corollary}
	
	\begin{remark}
	    Classically, $L^V=L^M$ if $M$ is a proper class transitive model of $\KP$. Furthermore, if $M$ is a transitive (set or class) model of $\KP$, then $L^M=L\cap M$ only depends on $\Ord^M=\Ord\cap M$.
	    
	    This is not constructively valid even if $M$ is a proper class since there is no reason to believe \mbox{$\Ord\cap V=\Ord\cap M$.} In fact, we can construct a Kripke model of $\IZF$ that satisfies $\Ord\cap V\neq \Ord\cap L$. See \cite[Section 5.5]{MatthwesPhD} for the details.
	\end{remark}
	
	Of course, we can derive more: \cite{MatthwesPhD} proved that $L_{\kappa^\#}$ also thinks every set is included in a \emph{totally indescribable set}. Also, under the presence of Set Induction and Collection for $\mySigma^{j, M}$-formulas, we have
	
	\begin{theorem}[$\IKP_{j,M}$]\label{Theorem:IKPSigmaOrdimpliesIZF+BTEE}
	   Let $\lambda=\bigcup_{n\in\omega} j^n(\kappa^\#)$. Then $L_\lambda = \bigcup_{n\in\omega} L_{j^n(\kappa^\#)}$ and $ \lag L_\lambda, j\restricts L_\lambda \rag$ is a model of $\mathsf{IZF+BTEE} + \TIj$.
	\end{theorem}
	
	\begin{proof}
	    Set Induction and Collection for $\mySigma^{j,M}$-formulas are necessary to ensure the existence of $\lambda$.
	    
	    By definition of $L_\alpha$, $L_\lambda = \bigcup_{\alpha\in\lambda}\Def(L_\alpha) = \bigcup_{n\in\omega}\bigcup_{\alpha\in j^n(\kappa^\#)} \Def(L_\alpha) = \bigcup_{n\in\omega} L_{j^n(\kappa^\#)}$.
	    Now, the desired result follows from \autoref{Corollary:LambdaModelsIZFBTEEInd}.
	\end{proof}
	
	\subsection{Reinhardt sets}
	Let $j \colon V\to M$ be an elementary embedding and $K$ a critical point of $j$. It is not generally true that $j(K)$ is also a regular set, although $M$ believes it is. However, for a Reinhardt embedding $j \colon V\to V$, $j(K)$ is regular. This yields a better lower bound for the consistency strength of Reinhardt sets. Observe that the proof only requires $j$ to be $\Delta_0$-elementary because all of the formulas can be bound by $\Lambda$.
	
	\begin{theorem}[{$\IKP_j$}]\label{Theorem:ReinhardtCriticalPtModelsWA}
		If $K$ is Reinhardt, then $\Lambda$ satisfies $\mathsf{IZF+WA}$.
	\end{theorem}
	
	Since we already know that $\Lambda \models \IZF$, the above proposition follows immediately from the following lemma:
	
	\begin{lemma}
		For any $j$-formula $\phi$, $t\in\omega$ and $a,p\in j^t(K)$, $\{x\in a\mid \langle \Lambda, j\restricts\Lambda\rangle \models \phi(x, p)\}\in j^t(K)$.
	\end{lemma}
	
	\begin{proof}
		Let $F_\theta(a,p):=\{x\in a\mid \langle \Lambda, j\restricts\Lambda\rangle \models \theta(x,p)\}$ for a formula $\theta$. We will prove that{, for every formula $\theta$,} $F_\theta(a,p)\in j^t(K)$ for all $a,p\in j^t(K)$ by induction on the complexity of $\theta$. 
		
		Atomic cases in which $j$ does not appear follow immediately from the inaccessibility of $j^t(K)$. We consider the case for equality where $j$ appears, the cases for element-hood being similar. To do this, we need to show that for $a \in j^t(K)$, 
		\begin{equation*}
		    \{ \langle x, y \rangle \in a \mid x = j(y) \} \in j^t(K).
		\end{equation*}
		First, since $j^{t+1}(K)$ is inaccessible in $V$, $j^{t+1}(K)$ is Exp-closed by \autoref{Corollary:Inaccessiblesets-ExpClosed}, and therefore \mbox{$j \restricts a \in \mv(\prescript{a}{}j(a))$} is in $j^{t+1}(K)$. Next, since $j^{t+1}(K)$ is closed under intersections, $j \restricts a \cap a \in j^{t+1}(K)$. However, this is a subset of $a \in j^t(K)$ so, by \autoref{Lemma:KandLambda-samepowerset}, $j \restricts a \cap a = \{ \langle x, y \rangle \in a \mid x = j(y) \} \in j^t(K)$.
		
		Conjunctions, disjunctions and implications follow from the fact that $j^t(K)$ satisfies Union and \mbox{$\Delta_0$-Separation}: let us examine the proof for implications. 
		Suppose that $F_\phi(a,p),F_\psi(a,p)\in j^t(K)$ for any $a,p\in j^t(K)$. Then{, by $\Delta_0$-Separation,}
		\begin{equation*}
			F_{\phi\to\psi}(a,p) = \{x\in a \mid x\in F_\phi(a,p) \to x\in F_\psi(a,p)\}\in j^t(K).
		\end{equation*}
		Next, suppose that $\phi(x,p)$ is $\forall y \psi(x,y,p)$. For each $n\in\omega$, let
		\begin{equation*}
			S_n := \{x\in a \mid \forall y\in j^{t+n}(K) \langle \Lambda, j\restricts\Lambda\rangle \models\psi(x,y,p) \}.
		\end{equation*}
		We claim that $S_n\in j^t(K)$ for every $n\in\omega$. By the inductive hypothesis, for every $y\in j^{t+n}(K)$, 
		\begin{equation*}
			F_\psi(a,\langle y,p\rangle ):= \{x\in a \mid \langle \Lambda, j\restricts\Lambda\rangle \models\psi(x,y,p)\} \in j^{t+n}(K).
		\end{equation*}
		
		Then we can define a function $R \colon j^{t+n}(K)\to j^{t+n+1}(K)$ given by $R(y)=F_\psi(a,\langle y,p\rangle)$. Furthermore, $R(y)\in \mathcal{P}(a)\cap j^{t+n+1}(K)=\mathcal{P}(a)\cap j^{t+n}(K)$. 
		Hence the codomain of $R$ is $\mathcal{P}(a)\cap j^{t+n}(K)$, which is an element of $j^{t+n+1}(K)$. Since $K$ is regular, so inaccessible by \autoref{Theorem:CriticalPointModelsIZF}, so is $j^{t+n+1}(K)$.
		Hence by \autoref{Corollary:Inaccessiblesets-ExpClosed}, $R\in j^{t+n+1}(K)$. Hence $S_n =\bigcap\ran R \in j^{t+n+1}(K)$. By \autoref{Lemma:KandLambda-samepowerset} and $S_n\subseteq a \in j^t(K)$, we have $S_n \in j^t(K)$.
		
		Finally for this case, since $S_n\in j^t(K)$ for each $n\in\omega$, we can define a function $S \colon \omega\to j^t(K)$ by $S(n):=S_n$. By repeating the previous argument with the inaccessibility of $j^t(K)$ and $\omega\in j^t(K)$, we have $S\in j^t(K)$. Hence
		\begin{equation*}
			{\textstyle \bigcap_{n\in\omega} S_n} = \{x\in a\mid \forall n\in\omega \forall y\in j^{t+n}(K) \langle\Lambda, j\restricts\Lambda\rangle \models\psi(x,y,p)\} = F_\phi(a,p)
		\end{equation*}
		is also a member of $j^t(K)$.
		
		The last case is when $\phi(x,p)$ is $\exists y \psi(x,y,p)$. The proof is quite similar to the previous one. Similar to the previous case, for every $n\in\omega$, let
		\begin{equation*}
			S'_n := \{x\in a \mid \exists y\in j^{t+n}(K) \langle \Lambda, j\restricts\Lambda\rangle \models\psi(x,y,p) \}.
		\end{equation*}
		
		We again prove that for each $n\in\omega$, $S'_n\in j^t(K)$ by first obtaining that $R\in j^{t+n+1}(K)$, where $R \colon y\mapsto F_\psi(a,\langle y,p\rangle)$ is the same function we defined in the proof for the previous case.
		Therefore, $S'_n = \bigcup_{y\in j^{n+1}(K)}R(y) \in j^{n+1}(K)$. Furthermore, $S'_n\subseteq a\in j^t(K)$ shows $S'_n\in j^t(K)$. If we define $S'\colon \omega\to j^t(K)$ by $S'(n):=S'_n$, then $S'\in j^t(K)$. Hence $F_\phi(a,p)=\bigcup_{n\in\omega} S'(n)\in j^t(K)$.
	\end{proof}
	
	\subsection{Super Reinhardt sets}
	The following theorem shows that super Reinhardt sets reflect first-order properties of $V$.
	
	\begin{theorem}[$\CGB_\infty$] \label{Theorem:SuperReinhardtReflection}
		Let $K$ be a super Reinhardt set. Then $K$ is an elementary submodel of $V$. That is, for every formula (without $j$) $\phi(\vec{x})$ all of whose variables are displayed,
		\begin{equation*}
			\forall\vec{x}\in K \phi^K(\vec{x})\lr\phi(\vec{x}).
		\end{equation*}
	\end{theorem}
	
	\begin{proof}
		Atomic cases and the cases for logical connectives are trivial. Hence we focus on quantifications.
		\begin{itemize}
			\item Case $\forall$: assume that $a\in K$ and $\phi(x,a)$ is absolute between $K$ and $V$. Then clearly we have \mbox{$\forall x \phi(x,a)\to\forall x\in K \phi^K(x,a)$}.
			Conversely, assume that $\forall x\in K \phi^K(x,a)$. Fix $b\in V$ and $j$ such that $b\in j(K)$. Then $\forall x\in j(K) \phi^{j(K)}(x,a)$ implies $\phi^{j(K)}(b,a)$.
			Furthermore, we can see that $\phi(x,a)$ is also absolute between $j(K)$ by applying $j$ to our inductive hypothesis; $\forall x, a\in K \phi^K(x,a)\lr\phi(x,a)$.
			Thus $\phi(b,a)$ for all $b\in V$.
			\item Case $\exists$: assume the same conditions on $a$ and $\phi(x,a)$ as we did before. Then obviously we have $\exists x\in K \phi^K(x,a)\to \exists x \phi(x,a)$.
			For the converse, assume that there is $b$ such that $\phi(b,a)$. Find $j$ such that $b\in j(K)$. Since $\phi$ is also absolute between $j(K)$ and $V$, we have $\exists x\in j(K) \phi^{j(K)}(x,a)$. Thus $\exists x\in K \phi^K(x,a)$. \qedhere
		\end{itemize}
	\end{proof}
	
	Note that the above theorem requires the full elementarity of elementary embeddings. 
	Next, we shall see that not only do strong Reinhardt sets reflect all first-order properties of $V$, but also they contain every true subset of a member of themselves:
	
	\begin{proposition}[$\CGB_\infty$ ] \label{Prop:SuperReinhardtPowerInaccessible}
		Let $K$ be a super Reinhardt set and $a\in K$. If $b\subseteq a$, then $b\in K$.
		Especially, $\mathcal{P}^K(a)=\mathcal{P}(a)$ for $a\in K$, so $K$ is power inaccessible.
	\end{proposition}
	
	\begin{proof}
		Find $j \colon V\to V$ such that $b\in j(K)$. By \autoref{Lemma:PowerSetPreserving}, $j(b)=b\in j(K)$, so $b\in K$.
	\end{proof}
	
	\begin{corollary}[$\CGB_\infty$] \label{Corollary:SuperReinhardtmodelsIZFpIEA}
	        Suppose that there is a super Reinhardt set. Then $V$ is a model of $\mathsf{IZF + pIEA}$.
	\end{corollary}
	
	However, it is not, in general, true that we will also satisfy the full second-order theory of $\IGB$. The issue here is that there is no reason why a (not first-order definable) class should be amenable. On the other hand, one should observe that by restricting our attention to amenable classes we will obtain a model of $\IGB$. Also, note that being a super Reinhardt is a second-order property, so its existence does not reflect down to $K$.
	
	Bagaria-Koellner-Woodin \cite{BagariaKoellnerWoodin2019} showed that super Reinhardt cardinals \emph{rank-reflect} Reinhardt cardinals, that is, there is an inaccessible cardinal $\gamma$ such that $(V_\gamma,V_{\gamma+1})$ models $\ZF_2$ with a Reinhardt cardinal.
	We will show in \autoref{Theorem:LowerBound-superReinhardt} that $\CGB_\infty$ with a super Reinhardt set interprets $\ZF$ with a proper class of $\gamma$ such that $(V_\gamma,V_{\gamma+1})\models \mathsf{ZF_2}+\text{`there is a Reinhardt cardinal.'}$
	However, its proof `mixes up' large set arguments over $\CZF$ with a double-negation translation, so the following question is still open:
	
	\begin{question}
		Working over $\CGB_\infty$ with a super Reinhardt set, can we prove there is an inaccessible set $M$ such that $(M,\mathcal{P}(M))$ satisfies $\CGB_\infty$ with the existence of a Reinhardt set?
	\end{question}
	
	However, we can still derive various large set principles from super Reinhardtness. For example, we can see that super Reinhardtness implies the analogue of $j \colon V_{\lambda+n}\prec V_{\lambda+n}$ over $\ZF$:
	
	\begin{proposition}[$\CGB_\infty$]
		Assume that there is a super Reinhardt set. Define $V_\alpha(x)$ recursively as $V_\alpha({ x})=\TC({ x})\cup \bigcup_{\beta\in\alpha}\mathcal{P}(V_\beta(x))$.
		If $j(\xi)=\xi$, then for each set $a$ we can find $\Lambda\ni a$, which is a countable union of power inaccessible sets, with an elementary embedding $j \colon V_\xi(\Lambda)\to V_\xi(\Lambda)$.
		
		Especially, for each $n\in\omega$ and $a\in V$, we can find $\Lambda\ni a$, which is a countable union of power inaccessible sets, such that there is an elementary embedding $j \colon \mathcal{P}^n(\Lambda)\to \mathcal{P}^n(\Lambda)$.
	\end{proposition}
	
	\begin{proof}
		Let $K$ be a super Reinhardt set and $j$ be an elementary embedding with a critical point $K$ such that $a\in j(K)$. Now let $\Lambda = j^\omega(K)$. We can see that $j\restricts V_\xi(\Lambda) \colon V_\xi(\Lambda)\to j(V_\xi(\Lambda)) = V_\xi(\Lambda)$ is the desired elementary embedding. 
		The latter claim follows by letting $\xi=n$.
	\end{proof}
	
	\subsection{Totally Reinhardt sets}
	Then how about the case for totally Reinhardt sets? We examined that if $K$ is super Reinhardt, then $K\prec V$. However, $K$ does not reflect $j$-formulas. We can see that $A$-super Reinhardt sets reflect not only usual set-theoretic formulas, but also $A$-formulas. The proof of the following theorem is identical to that of \autoref{Theorem:SuperReinhardtReflection}, so we omit it. We note here that in the theorem we will not need to assume that $A \cap K$ is a set, and if $A$ is not amenable then in fact this will not be the case.
	
	\begin{theorem}[$\CGB_\infty$] \pushQED{\qed}
		Let $K$ be an $A$-super Reinhardt set. Then $K$ reflects every $A$-formula. That is, for every formula $\phi(X, \vec{x})$ with one class parameter $X$ and all of whose variables are displayed,
		\begin{equation*}
			\forall\vec{x}\in K [\phi^K(A \cap K, \vec{x})\lr\phi(A, \vec{x})]. \qedhere
		\end{equation*}
	\end{theorem}
	
	Note that $A$-elementarity of $j$ is necessary for the above theorem: The reader can see that the proof of \autoref{Theorem:SuperReinhardtReflection} applies $j$ to the inductive hypothesis, namely, $\phi(\vec{x})$ is absolute between $K$ and $V$, and this is where we need the $A$-elementarity. We can also see that the proof breaks down if we do not assume $A$-elementarity: if the proof were to work without $A$-elementarity, then \autoref{Theorem:SuperReinhardtReflection} would hold even for $j$-formulas. This would imply $K$ thinks there is a critical point of $j$, which is invalid.
	
	One consequence of the reflection of $A$-formulas is that $V$ satisfies Full Separation for $A$-formulas when $A$ is amenable. This is because, for amenable classes, it is relatively straightforward to prove that $K \models \mathsf{Sep}_A$.
	The proof of the following lemma is similar to that of \autoref{Theorem:CriticalPointModelsIZF}, so we omit it.
	
	\begin{lemma}[$\CGB_\infty$] \pushQED{\qed}
		Let $A$ be an amenable class and $K$ be an $A$-super Reinhardt set.
		Then $K$ satisfies Full Separation for $A$-formulas. \qedhere
	\end{lemma}
	
	\begin{corollary}[$\CGB_\infty$]\label{Corollary:TRprovesASep} \pushQED{\qed}
		Assume that there is an $A$-super Reinhardt set for an amenable $A$. Then Full Separation for $A$-formulas hold.
		Especially, if $V$ is totally Reinhardt, then Full Separation for $A$-formulas hold for all amenable classes $A$. \qedhere
	\end{corollary}
	
	Note that every $A$-super Reinhardt set is super Reinhardt, so is power inaccessible by \autoref{Prop:SuperReinhardtPowerInaccessible}. 
	In sum, {$\mathsf{CGB_\infty+TR}$} proves {$\mathsf{IGB_\infty+TR}$} without Class Separation. It is unknown if Class Separation follows from the remaining axioms of $\mathsf{IGB_\infty+TR}$. However, we can still see that $\mathsf{IGB_\infty+TR}$ is interpretable within itself without Class Separation:
	
	\begin{theorem}[$\CGB_\infty$] \label{Theorem:CGBTRinterpretsIGBTR}
		$\mathsf{CGB_\infty+TR}$ interprets $\mathsf{IGB_\infty+TR}$.
	\end{theorem}

	\begin{proof}
		We proved that $\mathsf{CGB_\infty+TR}$ proves every axiom of $\mathsf{IGB_\infty+TR}$ except for Class Separation. We claim that $\mathsf{IGB_\infty+TR}$ is interpreted in itself without Class Separation.
		
		Let $\phi$ be a formula of $\mathsf{IGB_\infty+TR}$, and $\phi^\mathfrak{a}$ be a formula obtained by bounding every second-order quantifier to the collection of all amenable classes. That is, we get $\phi^\mathfrak{a}$ by replacing every $\forall^1 X$ and $\exists^1 X$ occurring in $\phi$ to $\forall^1 X (\text{$X$ is amenable})\to\cdots$ and $\exists^1 X (\text{$X$ is amenable})\land\cdots$.
		
		By \autoref{Corollary:TRprovesASep}, Class Separation holds for amenable classes: that is, if $A$ is amenable and $\phi(x,p,C)$ is a second-order quantifier-free formula whose free variables are all expressed, then $\{x\mid \phi(x,p,A)\}$ is also amenable. Thus the $\mathfrak{a}$-interpretation of Class Separation holds. Moreover, it is easy to see that the $\mathfrak{a}$-interpretation of the other axioms of $\mathsf{IGB_\infty+TR}$ are valid.
	\end{proof}

	\section{Heyting-valued interpretation and the double-negation interpretation}\label{Section:HeytingInterpretation}
	In this section, we will develop tools to analyze the consistency strength of large set axioms over constructive set theories. The main tool we will use is the double-negation translation.
	Especially, we will heavily rely on Gambino's Heyting-valued interpretation (\cite[Chapter 5]{GambinoPhDThesis} or \cite{Gambino2006}) with the double-negation topology which has a similar presentation to the classical Boolean-valued interpretations. Grayson \cite{Grayson1979} and Bell \cite{Bell2014} have also introduced a way to interpret classical $\mathsf{ZF}$ from $\IZF$ using an easier method whose connection to Boolean-valued models is even clearer. However, we will use Gambino's method because it works over the much weaker background theory of $\CZFminus$. 
	
	\subsection{Heyting-valued interpretation of \texorpdfstring{$\CZFminus$}{CZF--}}
	Forcing is a powerful tool to construct a model of set theory. Gambino's definition of Heyting-valued model (or alternatively, forcing) opens up a way to produce models of $\CZFminus$. His Heyting-valued model starts from a \emph{formal topology}, which formalizes a poset of open sets with a covering relation:
	
	\begin{definition}
		A \emph{formal topology} is a structure $\mathcal{S}=(S,\le,\vartriangleleft)$ such that $(S,\le)$ is a poset and $\vartriangleleft\subseteq S\times\mathcal{P}(S)$ satisfies the following conditions:
		\begin{enumerate}
			\item if $a\in p$, then $a\vartriangleleft p$,
			\item if $a\le b$ and $b\vartriangleleft p$, then $a\vartriangleleft p$,
			\item if $a\vartriangleleft p$ and $\forall x\in p (x\vartriangleleft q)$, then $a\vartriangleleft q$, and
			\item if $a\vartriangleleft p, q$, then $a\vartriangleleft (\downwards p)\cap (\downwards q)$, where $\downwards p = \{b\in S \mid \exists c \in p (b\le c )\}$. 
		\end{enumerate}
	\end{definition}
	
	Intuitively, $S$ describes a basis of a topology, and $\vartriangleleft$ is a covering relation. Then, for each collection of `open sets' $p$, we have the notion of a \emph{nucleus}, $\jmath p$, which is the set of all open sets that are covered by $p$. We can view $\jmath p$ as a `union' of all open sets in $p$, defined by
	\begin{equation*}
		\jmath p = \{x\in S \mid x\vartriangleleft p\}.
	\end{equation*}
	Then the class $\Low(\mathcal{S})_\jmath$ of all lower subsets\footnote{A subset $p\subseteq S$ is a lower set if $\downwards p = p$.} that are stable under $\jmath$ (i.e., $\jmath p = p$) form a \emph{set-generated frame}:
	
	\begin{definition}
		A structure $\mathcal{A}=(A,\le,\bigvee,\land,\top,g)$ is a \emph{set-generated frame} if $(A,\le,\bigvee,\land,\top)$ is a complete distributive lattice with the generating set $g\subseteq A$, such that the class $g_a=\{x\in g\mid x\le a\}$ is a set, and $a=\bigvee g_a$ for any $a\in A$.
	\end{definition}
	
	The reader is reminded that we can endow a Heyting algebra structure over a set-generated frame. For example, we can define $a\to b$ by $a\to b = \bigvee\{x\in g\mid x\land a\le b\}$, $\bot$ by $\varnothing$, and $\bigwedge p$ by $\bigvee\{x\in g\mid \forall y\in p (x\le y)\}$.
	
	\begin{proposition}\pushQED{\qed}
		For every formal topology $\mathcal{S}$, the class $\Low(\mathcal{S})_\jmath$ has a set-generated frame structure under the following definition of relations and operations:
		\begin{itemize}
			\item $p\land q=p\cap q$,
			\item $\bigvee p=\jmath(\bigcup p)$,
			\item $\top = S$,
			\item $\le$ as the inclusion relation, and
			\item $g=\{\{x\}\mid x\in S\}$.
		\end{itemize}
		Furthermore, we can make $\Low(\mathcal{S})_\jmath$ a Heyting algebra with the following additional operations:
		\begin{itemize}
			\item $p\lor q=\jmath(p\cup q)$,
			\item $p\to q=\{x\in S\mid x\in p\to x\in q\}$,
			\item $\bigwedge p=\bigcap p$. \qedhere
		\end{itemize}
	\end{proposition}
	
	We extend the nucleus $\jmath$ to a lower subclass $P\subseteq S$, which is
	a subclass of $\mathcal{S}$ satisfying \mbox{$P = \downwards P := \{a\in\mathcal{S} \mid \exists b\in P (a\le b)\}$,} by taking 
	\begin{equation*}
		JP := \bigcup\{\jmath p \mid p\subseteq P\}.
	\end{equation*}
	Then we define Heyting operations for classes as follows:
	\begin{itemize}
		\item $P\land Q=P\cap Q$,
		\item $P\lor Q=J(P\cup Q)$,
		\item $P\to Q =\{x\in S\mid x\in P\to x\in Q\}$.
	\end{itemize}
	
	For a set-indexed collection of classes $\{P_x\mid x\in I\}$, take $\bigwedge_{x\in I} P_x=\bigcap_{x\in I} P_x$ and $\bigvee_{x\in I} P_x=J\left(\bigcup_{x\in I} P_x\right)$. \\
    
    The following results can be proven by direct calculation and so we omit their proofs:
	
	\vbox{
	\begin{remark} \label{Remark:lowersubclassproperties} \phantom{a}
	    \begin{itemize}
	        \item If $p$ is a set then $Jp = \jmath p$.
	        \item For any lower subclass $P \subseteq S$, $P \subseteq JP$.
	        \item If $P$, $Q$ and $R$ are subclasses of $S$ which are stable under $J$ then $R \subseteq (P \rightarrow Q)$ if and only if $R \land P \subseteq Q$.
	        \item If $\{ P_x \mid x \in I \}$ is a family of subclasses of $S$ such that for each $x \in I$, $JP_x = P_x$ and $R$ is another subclass of $S$ such that $JR = R$, then $R \subseteq \bigwedge_{x \in I} P_x$ if and only if $R \subseteq P_x$ for each $x \in I$.
	    \end{itemize}
	\end{remark}
	}

	Given $x \in S$ one can consider the class of all $p \subseteq S$ such that $x \in \jmath p$, namely the collection of all covers of $x$. In general, this need not be a set however in many cases it can be sufficiently approximated by a set. This is particularly vital to verify Subset Collection in our eventual model in \autoref{Theorem:CZFPersistence}. 
	
	\begin{definition} \label{Definition:FormalTopology}
	    A formal topology $(\mathcal{S},\le,\vartriangleleft)$ is said to be \emph{set-presentable} if there is a \emph{set-presentation} \mbox{$R:\mathcal{S}\to\mathcal{P(P(S))}$,} which is a set function satisfying
	    \begin{equation*}
	        a\vartriangleleft p \lr \exists u\in R(a) [u\subseteq p]
	    \end{equation*}
	    for all $a\in \mathcal{S}$ and $p\in\mathcal{P(S)}$.
	\end{definition}
	
	Since some readers may not be familiar with the definition of formal topology,
	we give here a brief informal description of it. We then refer the reader to Chapter 4 of \cite{GambinoPhDThesis} or Chapter 15 of \cite{AczelRathjen2010} for a more detailed account. The notion of formal topology stems from an attempt to formulate point-free topology over a predicative system such as Martin-L\"of type theory. Thus we may view $(\mathcal{S}, \le)$ as a collection of open sets.
	
	We usually describe open sets by using a subbasis, and sometimes the full topology is more complex than the subbasis. This issue can particularly arise in the constructive set-theoretic context and, even worse, it could be that while a subbasis is a set, the whole topology generated by the subbasis is a proper class. In that case, we want to have a simpler surrogate for the full topology. This explains why we do not define the formal topology as a $\bigvee$-semilattice.
	Since the covering relation plays a pivotal role in topology and sheaf theory, we also formulate the covering relation $\vartriangleleft$ into the definition of formal topology. Producing $\Low(\mathcal{S})_\jmath$ from the formal topology corresponds to recovering the full topology from a subbasis.
	
	Although we hope that $\mathcal{S}$ will be as simple as possible, the use of a `complex' formal topology is sometimes unavoidable. For example, even the natural double-negation formal topology defined over $\mathcal{P}(1)$, which we will shortly define, is a class unless we have Power Set. Thus we want to define a `small formal topology' separately, and the notion of set-presented formal topology is exactly for such a purpose. Roughly, the set-representation $R$ decomposes $a\in\mathcal{S}$ into some collection of `open sets,' and we can track the covering relation by using $R$. \\ 

	% Heyting-valued universe
	The Heyting universe $V^\mathcal{S}$ over $\mathcal{S}$ is defined inductively as follows: $a\in V^\mathcal{S}$ if and only if $a$ is a function from a set-sized subset of $V^{\mathcal{S}}$ to $\Low(\mathcal{S})_\jmath$. For each set $x$, we have the canonical representation $\check{x}\in V^{\mathcal{S}}$ of $x$ recursively defined by $\dom \check{x}=\{\check{y}\mid y\in x\}$ and $\check{x}(\check{y})=\top$.
	We can now define the Heyting interpretation, $\llbr\phi\rrbr\in \Low_\jmath$, with parameters in $V^\mathcal{S}$ as follows:
	
	\begin{definition} \label{Definition: Heyting interpretation}
		Let $\phi$ be a formula of first-order set theory and $\vec{a}\in V^\mathcal{S}$. Then we define the Heyting-valued interpretation, $\llbr\phi(\vec{a})\rrbr$, as follows:
		\begin{itemize}
			\item $\llbr a= b\rrbr=\left(\bigwedge_{x\in\dom a}a(x)\to\bigvee_{y\in\dom b} b(y)\land\llbr x=y\rrbr\right) \land \left(\bigwedge_{y\in\dom b}b(y)\to\bigvee_{x\in\dom a}a(x)\land\llbr x=y\rrbr\right)$,
			\item $\llbr a\in b\rrbr=\bigvee_{y\in\dom b} b(y)\land\llbr a=y\rrbr$,
			\item $\llbr\bot\rrbr=\bot$, $\llbr \phi\land \psi\rrbr=\llbr\phi\rrbr\land \llbr\psi\rrbr$, $\llbr \phi\lor \psi\rrbr=\llbr\phi\rrbr\lor \llbr\psi\rrbr$, $\llbr \phi\to \psi\rrbr=\llbr\phi\rrbr\to \llbr\psi\rrbr$, and $\llbr \lnot\phi\rrbr=\llbr\phi\to\bot\rrbr$,
			\item $\llbr\forall x\in a\phi(x)\rrbr=\bigwedge_{x\in\dom a}a(x)\to \llbr\phi(x)\rrbr$ and $\llbr\exists x\in a\phi(x)\rrbr=\bigvee_{x\in\dom a}a(x)\land \llbr\phi(x)\rrbr$,
			\item $\llbr\forall x\phi(x)\rrbr=\bigwedge_{x\in V^\mathcal{S}}\llbr\phi(x)\rrbr$ and $\llbr\exists x\phi(x)\rrbr=\bigvee_{x\in V^\mathcal{S}}\llbr\phi(x)\rrbr$.
		\end{itemize}
		
	We write $V^{\mathcal{S}} \models \phi$ when $\llbr \phi \rrbr = \top$ holds, and in such a case we say that $\phi$ is \emph{valid} in $V^{\mathcal{S}}$.
	\end{definition}

    The next result is that the interpretation validates every axiom of $\CZFminus$. Since the proof of this was already done in \cite{Gambino2006}, we will omit most of the proof. However we will replicate the proof of the validity of Strong Collection because it is more involved and the method will be necessary to show the persistence of BCST-regularity in \autoref{Theorem:BCSTset-preserving}.
	
	Let $a\in V^{\mathcal{S}}$ and $R$ be a class. We want to show the following statement holds:
	\begin{equation*}
	    \llbr (R\colon a\rrarrows V) \to \exists b (R\colon a\lrlrarrows b)\rrbr = \top.
	\end{equation*}
	Here $\to$ is translated to a Heyting implication operation, and the focal property of the implication operation is the following: $q\to r=\top$ if and only if for every $p\le q$, we have $p\le r$.
	
	Hence we try the following strategy: take any $p\in\Low(\mathcal{S})_\jmath$ such that $p\subseteq \llbr R \colon a \rrarrows V \rrbr$(that is, $p \leq \llbr R \colon a \rrarrows V \rrbr$ in the inclusion ordering). Then we claim that we can find $b\in V^\mathcal{S}$ such that $p\subseteq \llbr R\colon a\lrlrarrows b\rrbr$.
	
	We will, and must use a form of Strong Collection to prove the validity of Strong Collection over $V^\mathcal{S}$. The role of Strong Collection is to confine the codomain of a class-sized multi-valued function to a set-sized range.
	Next note that, using the third and fourth items in \autoref{Remark:lowersubclassproperties}, the assumption $p\subseteq \llbr R\colon a\rrarrows V\rrbr$ is equivalent to $p\land a(x)\subseteq \llbr \exists y\ R(x,y)\rrbr$ for all $x\in\dom a$. So, using this assumption, we can codify the relation $R$ by using
	\begin{equation*}
	    P = \{\langle x,y,z\rangle \mid x\in\dom a,\ y\in V^\mathcal{S},\ z\in p\land a(x)\land \llbr R(x,y)\rrbr\}.
	\end{equation*}

    Intuitively, $P$ encodes a family of classes
	\begin{equation*}
	    \{p\land a(x)\land \llbr R(x,y)\rrbr \mid x\in\dom a,\ y\in V^\mathcal{S}\}.
	\end{equation*}
	Since $\llbr R(x,y)\rrbr$ could be a proper class in which we cannot form a collection of all $p\land a(x)\land \llbr R(x,y)\rrbr $, we introduce $P$ to code this family into a single class. 
	
	We will construct $b$ by searching for an appropriate subset of $P$ and making use of it. We will cast appropriate lemmas when we need them.
	
	\begin{lemma}\label{Lemma:StrongCollectionValid-MainLemma}
	    Fix $p\in \Low(\mathcal{S})_\jmath$ and $a\in V^\mathcal{S}$ such that $p\subseteq \llbr R\colon a\rrarrows V\rrbr$.
	    Let $P$ be the class we defined before. Then we can find a set $r\subseteq P$ such that $p\land a(x)\subseteq \jmath \{ z\mid \exists y \ \langle x,y,z\rangle\in r\}$ for all $x \in \dom a$.
	\end{lemma}
	
	\begin{proof}
	    Before starting the proof, let us remark that the $b$ we will eventually construct as our witness for Strong Collection will satisfy $\jmath \{ z\mid \exists y \ \langle x,y,z\rangle\in r\}\subseteq \llbr \exists y\in b\ R(x,y)\rrbr$.
	    
	    Observe that $p\subseteq \llbr R\colon a\rrarrows V\rrbr$ is equivalent to $\forall x\in\dom a (p\land a(x)\subseteq \bigvee_{y\in V^\mathcal{S}} \llbr R(x,y)\rrbr)$. Hence we have
	    
	    \begin{equation*}
	        p\land a(x)\subseteq \bigvee_{y\in V^\mathcal{S}} p\land a(x)\land \llbr R(x,y)\rrbr
	        = J \left( \bigcup_{y\in V^\mathcal{S}} p\land a(x)\land \llbr R(x,y)\rrbr \right)
	    \end{equation*}
	    
	    for every $x\in\dom a$. We want to have a family of classes $Q_x = \bigcup_{y\in V^\mathcal{S}} p\land a(x)\land \llbr R(x,y)\rrbr$.
	    For this, we define the coding $Q$ of the family $Q_x$ by
	    \begin{equation*}
	        Q = \{\langle x,z\rangle \mid \exists y\in V^\mathcal{S} [\langle x,y,z\rangle \in P]\}
	    \end{equation*}
	    and let $Q_x = \{z\mid\langle x,z\rangle\in Q\}${, which one can easily see satisfies our requirement}.
	    Then $p\land a(x)\subseteq JQ_x$. By the definition of $J$, we have the following sublemma, which is Lemma 2.8 of \cite{Gambino2006}:
	   
	   \renewcommand{\qedsymbol}{$\dashv$}
	    \begin{lemma} 
	        Let $P$ be a lower subclass of $\mathcal{S}$ and $u\subseteq JP$. Then we can find $v\subseteq P$ such that $u\subseteq \jmath v$. \qed
	    \end{lemma}
	    
	    In sum, we have that for each $x\in \dom a$ there is a $v\subseteq Q_x$ such that $p\land a(x)\subseteq \jmath v$. The following lemma shows we can find $v$ in a uniform way:
	    
	    \begin{lemma} \label{Lemma:TargetSubsetMVChoice-uniform-overV}
	        Let $a$ be a set, $S \colon a \rrarrows V$ a multi-valued class function, $Q\subseteq a\times V$ a class. 
	        For each $x \in a$, let $Q_x := \{ z \mid \lag x, z \rag \in Q\}$.
	        Moreover, assume that
		\begin{enumerate}
			\item For each $x\in a$ there is $u\subseteq Q_x$ such that $S(x,u)$ holds and
			\item (Monotone Closure) {If} $S(x,u)$ holds and $u\subseteq v\subseteq Q_x$ then $S(x,v)$, 
		\end{enumerate}
		then there is $f\colon a\to V$ such that $f(x)\subseteq Q_x$ and $S(x,f(x))$ for all $x\in a$.
	    \end{lemma}
	    
	    \begin{proof}
	        Consider the multi-valued function $S'$ of domain $a$ defined by
	        \begin{equation*}
	            S'(x,u) \qquad \text{if and only if} \qquad S(x,u)\text{ and } u\subseteq Q_x.
	        \end{equation*}
	        
	        By mimicking the proof of \autoref{Lemma:SetMV}, we can find a set $g\colon a\rrarrows V$ such that $g\subseteq S'$.
	        Now take $f(x)=\bigcup g_x = \bigcup \{u \mid \lag x,u\rag \in g\}$, which is a set by Union and Replacement.
	        We can see that $f(x) \subseteq Q_x$ for each $x\in a$. By the first clause of our assumptions and monotone closure of $S$, we have $S(x,f(x))$ for all $x\in a$.
	    \end{proof}
	    
	    Let us return to the proof of \autoref{Lemma:StrongCollectionValid-MainLemma}. Consider the relation $S$ defined by
	    \begin{equation*}
	        S(x,u)\quad\text{if and only if}\quad(u\subseteq Q_x\text{ and }p\land a(x)\subseteq \jmath u).
	    \end{equation*}
	    It is easy to see that $S$ satisfies the hypotheses of \autoref{Lemma:TargetSubsetMVChoice-uniform-overV}. Therefore, by \autoref{Lemma:TargetSubsetMVChoice-uniform-overV} applied to $S$, we can find $f\colon \dom a\to \Low(\mathcal{S})_\jmath$ such that $f(x)\subseteq Q_x$ and $p\land a(x)\subseteq \jmath f(x)$. 
	    
	    Recall that $f(x)\subseteq Q_x = \bigcup_{y\in V^\mathcal{S}} p\land a(x)\land\llbr R(x,y)\rrbr$. Hence for each $x\in\dom a$ and $z\in f(x)$ we can find $y\in V^\mathcal{S}$ such that $z\in p\land a(x)\land \llbr R(x,y)\rrbr$. However, the class of such $y$ may not be a set, so we will apply Strong Collection to find a subset of the class uniformly.
	    
	    Now let
	    \begin{equation*}
	        q=\{\lag x,z\rag\mid x\in\dom a,\ z\in f(x)\}.
	    \end{equation*}
	    Then for each $\lag x,z\rag \in q$ there is $y$ such that $\lag x,y,z\rag\in P$. Now consider the multi-valued function $P'\colon q\rrarrows V^\mathcal{S}$ defined by
	    \begin{equation*}
	        P'(\lag \lag x,z\rag,y\rag) \quad\text{if and only if} \quad P(x,y,z),
	    \end{equation*}
	    By Strong Collection applied to $\mathcal{A}(P')\colon q \rrarrows q \times V^\mathcal{S}$, we have $d$ such that $\mathcal{A}(P')\colon q \lrlrarrows d$. By \autoref{Lemma:PrelimAdjuectmentFtn}, $d\subseteq P'$ and $d\colon q \rrarrows V^\mathcal{S}$. Define $g(x,z) = \{y \mid \lag \lag x,z\rag, y\rag \in d\}$. Then we can see that $y\in g(x,z)$ implies $\lag x,y,z\rag \in P$.
	    
	    Finally, let
	    \begin{equation*}
	        r = \{\lag x,y,z\rag \mid \lag x,z\rag \in q,\ y\in g(x,z)\}.
	    \end{equation*}
	    
	    It is clear that $r$ is a set. We know that $p\land a(x)\subseteq \jmath f(x)$, and $z\in f(x)$ implies $\exists y[y\in g(x,z)]$, so $\exists y [\lag x,y,z\rag \in r]$. By combining these facts, we have $p\land a(x)\subseteq \jmath \{z\mid \exists y\ \lag x,y,z\rag \in r\}$. 
	    Thus concluding the proof of \autoref{Lemma:StrongCollectionValid-MainLemma}.
	    % Restoring the QED symbol 
	    \renewcommand{\qedsymbol}{$\square$}
	    \qedhere 
	\end{proof}
	
	\begin{proposition}\label{Proposition:StrongCollectionValid}
	    Working over $\CZFminus$, $V^\mathcal{S}$ validates Strong Collection.
	\end{proposition}
	
	\begin{proof}
	    Let $a\in V^\mathcal{S}$ and let $R$ be a class relation.
	    Let $p\in\Low(\mathcal{S})_\jmath$ and assume that $p\subseteq \llbr R\colon a\rrarrows V\rrbr$.
	    By \autoref{Lemma:StrongCollectionValid-MainLemma} (and using the previously defined notation), we can fix some $r\subseteq P$ such that 
	    \begin{equation*}
	        p\land a(x)\subseteq \jmath \{ z\mid \exists y \ \langle x,y,z\rangle\in r\}
	    \end{equation*}
	    holds for all $x\in\dom a$. Now define $b\in V^\mathcal{S}$ as follows: $\dom b = \{y \mid \exists x\exists z[\lag x,y,z\rag \in r]\}$, and \mbox{$b(y) = \jmath \{ z \mid \exists x \in \dom a [\lag x,y,z\rag \in r]\}$.}
	    
	    We claim that $p\subseteq \llbr R:a\lrlrarrows b \rrbr$.
	    \begin{itemize}
	        \item $p\subseteq \llbr R:a\rrarrows b\rrbr$: 
	        if $\lag x,y,z\rag \in r$, then $y\in \dom b$, hence
	        \begin{align*}
	            p\land a(x)& \subseteq \jmath \{ z\mid \exists y \ \langle x,y,z\rangle\in r\} \subseteq \bigvee_{y\in\dom b} \jmath \{z \mid \lag x,y,z\rag\in r\} 
	            \\ &\subseteq \bigvee_{y\in\dom b} \jmath\{z\mid \lag x,y,z\rag \in r\}\land p\land a(x)\land \llbr R(x,y)\rrbr 
	            \\ &\subseteq \bigvee_{y\in\dom b} b(y) \land \llbr R(x,y)\rrbr = \llbr \exists y\in b \ R(x,y)\rrbr.
	        \end{align*}
	        
	        \item $p\subseteq \llbr R:b\rrarrows a\rrbr$: note that $\lag x,y,z\rag\in r$ implies $x\in\dom a$.
	        Hence 
	        \begin{align*}
	            b(y)& = \jmath \{z \mid \exists x \in \dom a \lag x,y,z\rag \in r\} \subseteq \bigvee_{x\in\dom a} a(x)\land \llbr R(x,y) \rrbr = \llbr\exists x\in a\ R(x,y)\rrbr. \qedhere 
	        \end{align*}
	    \end{itemize}
	\end{proof}
	
	Hence we have
	
	\begin{theorem}\label{Theorem:CZFPersistence}
		Working over $\CZFminus$, the Heyting-valued model $V^\mathcal{S}$ also satisfies $\CZFminus$. If $\mathcal{S}$ is set-presented and Subset Collection holds, then $V^\mathcal{S}\models \CZF$.
		Furthermore, if our background theory satisfies Full Separation or Power Set, then so does $V^\mathcal{S}$ respectively.
	\end{theorem}
	
	\begin{proof}
		The first part of the theorem is shown by Gambino \cite{Gambino2006}, and we have already replicated the proof for the validity of Strong Collection. Hence we omit this part of the proof, and we concentrate on the preservation of Full Separation and Power Set.
		
		For Full Separation, it suffices to see that the proof for Bounded Separation over $V^\mathcal{S}$ also works for Full Separation. It actually works since Full Separation ensures $\llbr\phi\rrbr$ is a set for every formula $\phi$ because $\llbr\forall x \phi(x)\rrbr = \{s\in S \mid \forall x (x\in V^\mathcal{S} \to s\in \llbr\phi(x)\rrbr)\}$ and $\llbr\exists x\phi(x)\rrbr = \jmath(\{s\in S\mid \exists x (x\in V^\mathcal{S} \land s\in \llbr\phi(x)\rrbr)\})$.
		
		For Power Set, let $a\in V^\mathcal{S}$. We can show that $\Low_\jmath(S)$ is a set due to Power Set. Thus we have the name $b \in V^{\mathcal{S}}$ defined by $\dom b = {{}^{\dom a}}(\Low_\jmath(\mathcal{S}))$ and $b(c)=\top$. We claim that $b$ witnesses Power Set.
		
		Let $c\in V^\mathcal{S}$. We will find $d\in\dom b$ such that
		\begin{equation*}
			\llbr c\subseteq a\rrbr \le \bigvee_{d\in\dom b}\llbr c=d\rrbr.
		\end{equation*}
		
		Let $d$ be the name such that $\dom d=\dom a$ and $d(y)=\llbr y\in c\rrbr$.
		Then $d\in \dom b$. Furthermore, we have
		
		\begin{align*}
			\llbr c\subseteq a\rrbr =\bigwedge_{x\in\dom c} c(x)\to \left(\bigvee_{y\in\dom d} a(y)\land\llbr x=y\rrbr \right)
			\le \bigwedge_{x\in\dom c} c(x)\to \left(\bigvee_{y\in\dom d} a(y)\land\llbr x=y\rrbr\land c(x) \right)
			\\ \le \bigwedge_{x\in\dom c} c(x)\to \left(\bigvee_{y\in\dom d} \llbr x=y\rrbr\land \llbr x\in c\rrbr \right) \le \bigwedge_{x\in\dom c} c(x)\to \left(\bigvee_{y\in\dom d} \llbr x=y\rrbr\land \llbr y\in c\rrbr \right) 
			= \llbr c\subseteq d\rrbr
		\end{align*}
	 	and
	 	\begin{align*}
	 		\llbr d\subseteq c\rrbr = \bigwedge_{x\in \dom a}\llbr x\in c\rrbr \to \llbr x\in c\rrbr = \top.
	 	\end{align*}
 	Hence $\llbr c\subseteq a\rrbr \le \llbr c=d\rrbr$.
	\end{proof}
	
	Let us finish this subsection with some constructors, which we will need later.
	\begin{definition}
		For $\mathcal{S}$-names $a$ and $b$, $\up(a,b)$ is defined by $\dom(\up(a,b))=\{a,b\}$ and $(\up(a,b))(x)=\top$. $\op(a,b)$ is the name defined by $\op(a,b)=\up(\up(a,a),\up(a,b))$.
	\end{definition}
	
	$\up(a,b)$ is a canonical name for the unordered pair $\{a,b\}$ over $V^\mathcal{S}$. That is, we can prove that
	\begin{equation*}
		\llbr \forall x\forall y \forall z [z=\up(x,y) \lr z=\{x,y\}]\rrbr = \top.
	\end{equation*}
	Hence the name $\op(a,b)$ is the canonical name for the ordered pair given by $a$ and $b$ over $V^\mathcal{S}$.
	
	\subsection{Double-negation formal topology} \label{Subsection:DoubleNegationFormalTopology}
	Our main tool to determine the consistency strength of the theories in this paper is the Heyting-valued interpretation with the \emph{double-negation formal topology}.
	
	\begin{definition} \label{Definition: double-negation topology}
		The \emph{double-negation formal topology}, $\Omega$, is the formal topology $(1,=,\vartriangleleft)$, where $x\vartriangleleft p$ if and only if $\lnot\lnot(x\in p)$.
	\end{definition}
    
    Unlike set-sized realizability or set-represented formal topology, the double-negation topology and the resulting Heyting-valued interpretation need not be absolute between BCST-regular sets or transitive models of $\CZFminus$. 
    For example, even for a transitive model $M$ of $\CZFminus$, it need not be the case that $\Omega=\Omega^M$ holds. This is because $\vartriangleleft \subseteq 1 \times \mathcal{P}(1)$ and $\mathcal{P}(1) = \mathcal{P}^M(1)$ may not in general hold between transitive sets. For instance, if $\varphi^M$ is not logically equivalent to $\varphi$ then it is unclear why $\{ 0 \in 1 \mid \varphi\}$ should be in $M \cap \mathcal{P}(1)$.
    One example of the failure of $\Omega=\Omega^M$ happens when $M$ is the set $\mathsf{HF}$ of all hereditarily finite sets. \cite{Jeon2022Ackermann} proved that $\mathsf{HF}$ is a model of $\mathsf{CZF}$ without Infinity, and it additionally satisfies $\mathsf{LEM}$ for atomic formulas. This shows $\Omega^{\mathsf{HF}}=2$. However, there is no reason to believe that $\Omega=2$ in general unless we have $\Delta_0\mhyphen\mathsf{LEM}$.
    Hence we need a careful analysis of the double-negation formal topology, which is the aim of this subsection.
	
	We can see that the class of lower sets, $\Low(\Omega)=\{p\subseteq 1 \mid p= \downwards p\}$, is just the powerclass of 1, $\mathcal{P}(1)$, and the nucleus of $\Omega$ is given by the double complement $p^\dnot = (p^\lnot)^\lnot$, where
	\begin{equation*}
		p^\lnot=\{0\mid \lnot(0\in p)\},
	\end{equation*}
	so $p^\dnot = \{0\mid \lnot\lnot (0\in p)\}$. Hence $\Low(\Omega)_\jmath$ is the collection of all \emph{stable} subsets of 1, that is, those sets $p\subseteq 1$ such that $p=p^\dnot$.
	
	The main feature of $\Omega$ is that Heyting-valued interpretation over $\Omega$ forces a law of excluded middle for some class of formulas:
	
	\begin{proposition}\label{Proposition:DoubleNegationTopologyForcesLEM}
	    Let $p\in \Low(\Omega) = \mathcal{P}(1)$. Then $(p\cup p^\lnot)^\dnot =1$. As a corollary, if $\llbr \phi\rrbr$ is a set, then $\llbr\phi\lor\lnot\phi\rrbr=1$. Especially, 
	    
	    \begin{enumerate}
	        \item $\llbr\phi\lor\lnot\phi\rrbr=1$ holds for every bounded formula $\phi$, and
	        \item If Full Separation holds, then $\llbr\phi\lor\lnot\phi\rrbr=1$ holds for every $\phi$.
	    \end{enumerate}
	\end{proposition}
	
	\begin{proof}
	    $\lnot\lnot(p\cup \lnot p)=1$ follows from the fact that $\lnot\lnot(\phi\lor\lnot\phi)$ is derivable in intuitionistic logic. Moreover, $\llbr\phi\lor\lnot\phi\rrbr = (\llbr\phi\rrbr \cup \lnot \llbr\phi\rrbr)^\dnot$ if $\llbr\phi\rrbr$ is a set. 
	    Lastly, under $\CZFminus$, if $\phi$ is bounded then $\llbr \phi \rrbr$ is a set and if Full Separation holds then $\llbr \phi \rrbr$ is a set for every $\phi$.
	\end{proof}
	
	The following corollary is immediate from the previous proposition and \autoref{Theorem:CZFPersistence}:
	
	\begin{corollary} \label{Corollary:VOmegamodelsZFminus} \pushQED{\qed}
        If $V$ satisfies $\CZFminus$, then $V^\Omega\models \CZFminus+\Delta_0\mhyphen\mathsf{ LEM}$. Furthermore, if $V$ satisfies Full Separation, then $V^\Omega\models \ZFminus$. \qedhere 
	\end{corollary}

	We will frequently mention the relativized Heyting-valued interpretation. For a transitive model $A$ of $\CZFminus$, we can consider the construction of $V^\Omega$ internal to $A$. We define the following notions to distinguish relativized interpretation from the usual one.
	
	\begin{definition}\label{Definition:AOmega}
		Let $A$ be a transitive model of $\CZFminus$. Then $A^\Omega:=(V^\Omega)^A$ is the $\Omega$-valued universe relativized to $A$. If $A$ is a set, then $\tilde{A}$ denotes the $\Omega$-name defined by $\dom \tilde{A}:=A^\Omega$ and $\tilde{A}(x)=\top$ for all $x\in\dom \tilde{A}$.
	\end{definition}
	
    We shall see in \autoref{Lemma:HeytingUniversePreserving} that whenever $A$ is a set, so is $A^\Omega$ which means that the definition of $\tilde{A}$ is well defined. Also, we will
	often confuse $\tilde{A}$ and $A^\Omega$ if the context is clear. Finally, it is worth mentioning that if $j$ is an elementary embedding, then $j(\tilde{K})=\widetilde{j(K)}$, so we may write $j^n(\tilde{K})$ instead of $\widetilde{j^n(K)}$.
	
    As before with $\Omega$, we do not know whether $\Low(\Omega)_\jmath$ is equal to its relativization $(\Low(\Omega)_\jmath)^M$ in $M$.
    As a result, we do not know whether its Heyting-valued universe, $V^\Omega$, and Heyting-valued interpretation, $\llbr\cdot\rrbr$, are absolute. Fortunately, the formula $p\in\Low(\Omega)_\jmath$, which is $p\subseteq 1\land p=p^\dnot$, is $\Delta_0$. Hence $p\in\Low(\Omega)_\jmath$ is absolute between transitive models of $\CZFminus$. As a result, we have the following absoluteness result on the Heyting-valued universe:	
	
	\bigskip 
	
	\begin{lemma}\label{Lemma:HeytingUniversePreserving}
		Let $A$ be a transitive model of $\CZFminus$ without Infinity. Then we have $A^\Omega=V^\Omega\cap A$. Moreover, if $A$ is a set, then $A^\Omega$ is also a set.
	\end{lemma}
	
	\begin{proof}
		We will follow the proof of \cite[Lemma 6.1]{Rathjen2003Realizability}.
		Let $\Phi$ be the inductive definition given by
		\begin{equation*}
			\lag X,a \rag\in\Phi\iff \text{$a$ is a function such that $\dom a\subseteq X$, $a(x)\subseteq 1$ and $a(x)^\dnot=a(x)$ for all $x\in\dom a$}.
		\end{equation*}
		We can see that the $\Phi$ defines the class $V^\Omega$. Furthermore, $\Phi$ is $\Delta_0$, so it is absolute between transitive models of $\CZFminus$.
		By \autoref{Lemma:ItreationClass}, we have a class $J$ such that $V^\Omega = \bigcup_{a\in V} J^a$, and for each $s\in V$, $J^s=\Gamma_\Phi(\bigcup_{t\in s} J^t)$.
		Now consider the operation $\Upsilon$ given by
		\begin{equation*}
			\Upsilon(X):=\{a\in A\mid\exists Y\in A(Y\subseteq X\land \lag Y,a\rag\in\Phi)\}.
		\end{equation*}
		By \autoref{Lemma:ItreationClass} again, there is a class $Y$ such that $Y^s=\Upsilon(\bigcup_{t\in s}Y^t)$ for all $s\in V$.
		Furthermore, we can see that $Y^s\subseteq V^\Omega$ by induction on $a$.
		
		Let $Y=\bigcup_{s\in A}Y^s$. We claim by induction on $s$ that $J^s\cap A\subseteq Y$.
		Assume that $J^t\cap A\subseteq Y$ holds for all $t\in s$. If $a\in J^s\cap A$, then the domain of $a$ is a subset of $A\cap \left(\bigcup_{t\in s} J^t\right)$, which is a subclass of $Y$ by the inductive assumption and the transitivity of $A$.
		Moreover, for each $x\in\dom a$ there is $u \in A$ such that $x\in Y^u$. By Strong Collection over $A$, there is $v\in A$ such that for each $x\in \dom a$ there is $u\in v$ such that $x\in Y^u$. Hence $\dom a\subseteq \bigcup_{u\in v} Y^u$, which implies $a\in Y^v\subseteq Y$. 
		
		Hence $V^\Omega\cap A\subseteq Y$, which gives us that $Y=V^\Omega\cap A$. We can see that the construction of $Y$ is the relativized construction of $V^\Omega$ to $A$, so $Y=A^\Omega$. Hence $A^\Omega=V^\Omega\cap A$.
		If $A$ is a set, then $\Upsilon(X)$ is a set for each set $X$, so we can see by induction on $a$ that $Y^a$ is also a set for each $a\in A$. Hence $A^\Omega=Y=\bigcup_{a\in A} Y^a$ is also a set.
	\end{proof}
	
	We extended the nucleus $\jmath$ to $J$ for subclasses of $\Low\mathcal{S}$, and used it to define the validity of formulas of the forcing language.
	Now, we are working with the specific formal topology $\mathcal{S}=\Omega$, and in this case, for a class $P\subseteq 1$, $JP = \bigcup\{q^\dnot \mid q\subseteq P\}$. It is easy to see that $P\subseteq JP\subseteq P^\dnot$. 
	
	We also define the following relativized notion for any transitive class $A$ such that $1\in A$:
	\begin{equation*}
		J^AP = \bigcup\{q^\dnot \mid q\subseteq P\text{ and }q\in A\}.
	\end{equation*}
	If $P\in A$, then $J^AP=P^\dnot$, and in general, we have $P\subseteq J^AP\subseteq JP\subseteq P^\dnot$. Moreover, we can prove the following facts by straight forward computations.
	
	\begin{lemma}\label{Lemma:PropertiesOfRelativizedJ} \pushQED{\qed}
		Let $A$ and $B$ be transitive classes such that $1\in A,B$ and let $P\subseteq 1$ be a class.
		\begin{enumerate}
			\item $A\subseteq B$ implies $J^AP\subseteq J^BP$.
			\item If $\mathcal{P}(1)\cap A=\mathcal{P}(1)\cap B$, then $J^AP=J^BP$. \qedhere
		\end{enumerate}
	\end{lemma}
	
	However, the following proposition shows that $\CZFminus$ does not prove $J^AP=J^BP$ or $JP = P^\dnot$ in general:
	
	\begin{proposition} \phantom{a}
		\begin{enumerate}
			\item If $A\cap \mathcal{P}(1)=2$ (for example when $A=2$ or $A=V$ and $\Delta_0\mhyphen\mathsf{LEM}$ holds), then $J^AP=P$.
			
			\item If $P^\dnot\subseteq JP$ for every class $P$, then if $\Delta_0\mhyphen \mathsf{LEM}$ holds so does the law of excluded middle for arbitrary formulas. \qed
		\end{enumerate}
	\end{proposition}
	
	$J^A$ has a crucial role in defining Heyting-valued interpretation, but $J^A$ and $J^B$ might have different effects unless $A=B$. This causes absoluteness problems, which appears to be impossible to in general avoid.

    The following lemma states provable facts about relativized Heyting interpretations:
	
	\begin{lemma}\label{Lemma:HeytingInterpretationBetweenModels}
		Let $A\subseteq B$ be transitive models of $\CZFminus$. Assume that $\phi$ is a formula with parameters in $A^\Omega$.
		\begin{enumerate}
			\item If $\phi$ is bounded, then $\llbr\phi\rrbr^A=\llbr\phi\rrbr^B$.
			
			\item If $\phi$ only contains bounded quantifications, logical connectives between bounded formulas, unbounded $\forall$, and $\land$, then $\llbr\phi\rrbr^A=\llbr\phi^{\tilde{A}}\rrbr^B$.
			
			\item If every conditional appearing as a subformula of $\phi$ is of the form $\psi\to\chi$ for a bounded formula $\psi$ and a formula $\chi$, then $\llbr\phi\rrbr^A\subseteq \llbr\phi^{\tilde{A}}\rrbr^B$.
			
			\item If $\mathcal{P}(1)\cap A=\mathcal{P}(1)\cap B$, then $\llbr\phi\rrbr^A=\llbr\phi^{\tilde{A}}\rrbr^B$.
		\end{enumerate}
	\end{lemma}
	
	\bigskip
	
	\begin{proof} \phantom{a}
	    \begin{enumerate}
	        \item If $\phi$ is bounded, then $\llbr\phi\rrbr$ is defined in terms of double complement, Heyting connectives between subsets of 1, and set-sized union and intersection. These notions are absolute between transitive sets, so we can prove $\llbr\phi\rrbr$ is also absolute by induction on $\phi$. In the case of atomic formulas as an initial stage, we apply the induction on $A^\Omega$-names.
	        
	        \item We apply induction on formulas. For the unbounded $\forall$, we have
	        \begin{equation*}
    			\llbr\forall x\phi(x)\rrbr^A = \bigwedge_{x\in A^\Omega}\llbr \phi(x)\rrbr^A
    			= \bigwedge_{x\in A^\Omega} \llbr\phi^{\tilde{A}}(x)\rrbr^B
    			= \llbr\forall x\in\tilde{A} \phi^{\tilde{A}}(x)\rrbr^B.
    		\end{equation*}
    		The remaining clauses follows from the absoluteness argument we used in the proof of the previous item, so we omit them.
    		
    		\item The proof again uses induction on formulas. By the previous argument, if $\llbr\phi(x)\rrbr^A\subseteq\llbr\phi^{\tilde{A}}(x)\rrbr^B$ for all $x\in A^\Omega$, then $\llbr \forall x \phi(x)\rrbr^A\subseteq \llbr\forall x\in\tilde{A} \phi^{\tilde{A}}(x)\rrbr^B$. For an unbounded $\exists$, we have 
    		\begin{equation*}
    			\llbr\exists x\phi(x)\rrbr^A = J^A \left(\bigcup \{\llbr \phi(x)\rrbr^A\mid x\in A^\Omega\}\right)
    			\subseteq J^B \left(\bigcup \{\llbr \phi^{\tilde{A}}(x)\rrbr^B\mid x\in A^\Omega\}\right) = \llbr\exists x\in \tilde{A} \phi^{\tilde{A}}(x) \rrbr^B.
    		\end{equation*}
    		The remaining cases are straightforward except for $\to$, which requires some inspection to see how the inclusion works. By the assumption, our conditional is of the form $\psi\to\chi$ for some bounded formula $\psi$ and a (possibly unbounded) formula $\chi$. Since $\psi$ is bounded, we have $\llbr\psi\rrbr^A=\llbr\psi^{\tilde{A}}\rrbr^B$. Furthermore, we can see that if $a, b, c\in\Low(\Omega)_\jmath$ satisfies $a\subseteq b$, then $c\to a\subseteq c\to b$. Hence
    		\begin{equation*}
    		    \llbr \psi\to\chi\rrbr^A = \left( \llbr \psi\rrbr^A \to \llbr\chi\rrbr^A \right)\subseteq
    		    \left( \llbr \psi^{\tilde{A}}\rrbr^B \to \llbr\chi^{\tilde{A}}\rrbr^B \right) = \llbr(\psi\to\chi)^{\tilde{A}}\rrbr^B.
    		\end{equation*}
		    
		    \item We can see that $\llbr\phi(x)\rrbr^A = \llbr\phi^{\tilde{A}}(x)\rrbr^B$ holds by induction on $\phi$. The calculations we did before are helpful to see the inductive argument works for unbounded quantifications. In particular, we observe that if $\mathcal{P}(1) \cap A = \mathcal{P}(1) \cap B$ then $J^A = J^B$. \qedhere
            \end{enumerate}
	\end{proof}
	
   	\begin{remark}
        In the third clause of \autoref{Lemma:HeytingInterpretationBetweenModels}, if $A\in B$, then $\bigcup_{x\in A^\Omega} \llbr \phi^{\tilde{A}}(x)\rrbr^B\in \mathcal{P}(1)\cap B$. 
	    Thus, in this case, we have $J^B\left(\bigcup_{x\in A^\Omega} \llbr \phi^{\tilde{A}}(x)\rrbr^B\right) = \left(\bigcup_{x\in A^\Omega} \llbr \phi^{\tilde{A}}(x)\rrbr^B\right)^\dnot$.
	\end{remark}
	
	\begin{remark}
	    We can see that \autoref{Lemma:HeytingInterpretationBetweenModels} also holds for arbitrary formal topologies $\mathcal{S}\in A$. The proof is identical, and its verification is left to the reader. 
	\end{remark}

	\subsection{Preservation of small large sets} \label{Section:PreservationSmallLargeSets}
	Heyting-valued models do not necessarily preserve large sets unless we impose some additional restrictions.
	For example, on the one hand, Ziegler proved in \cite{ZieglerPhD} that large set properties are preserved under `small' pcas and formal topologies. 
	
	\begin{proposition}[Ziegler, \cite{ZieglerPhD}, Chapter 4] \label{Proposition:negationtranslationpreservesinaccessibles} \pushQED{\qed}
		Let $K$ be either a regular set, inaccessible set, critical set or Reinhardt set. 
		\begin{enumerate}
			\item Let $\mathcal{A}$ be a pca and $\mathcal{A}\in K$. Then $K[\mathcal{A}]:=V[\mathcal{A}]\cap K$ is a set and the realizability model $V[\mathcal{A}]$ thinks $K[\mathcal{A}]$ possesses the same large set property $K$ does.
			\item Let $\mathcal{S}$ be a formal topology that is set-represented by $R$. Assume that $\mathcal{S},R\in K$. Then $K^\mathcal{S} := V^\mathcal{S}\cap K$ is a set and $V^\mathcal{S}$ thinks the canonical name $\tilde{K}$ given by $\dom\tilde{K}=K^\mathcal{S}$, $\tilde{K}(a)=\top$ possesses the same large set property $K$ does.
		\end{enumerate}
	\end{proposition}
	
	On the other hand, however, our double-negation formal topology will usually lose large set properties. For example, $\CZFminus$ cannot prove that Heyting-valued interpretations under $\Omega$ preserve regular sets regardless of how many regular sets exist in $V$: If it were possible, then $\mathsf{CZF+REA}$ would interpret $\mathsf{ZF}$, but the latter theory proves the consistency of the former.
	
	Hence our preservation results under the double-negation topology are also quite limited, which is a great obstacle for deriving the consistency strength of large sets over constructive set theories. Fortunately, almost all lower bounds we derived in \autoref{Section:StructuralLowerbdd} are models of $\IZF$. Furthermore, we mostly work with power inaccessible sets instead of mere inaccessible sets.
	This will make deriving the lower bounds easier, and we will examine the detailed account of the aforementioned statements in this subsection.
	
	We mostly follow the proof of \autoref{Lemma:StrongCollectionValid-MainLemma} and \autoref{Proposition:StrongCollectionValid}. However, for the sake of verification, we will provide most of the details of relevant lemmas and their proofs. Throughout this section, $A$ and $B$ are classes such that
	\begin{itemize}
	    \item $A\in B$ (thus $A$ is a set),
	    \item $\mathcal{P}(1)\cap A=\mathcal{P}(1)\cap B$, 
	    \item $A$ is BCST-regular and $B$ is a transitive (set or class) model of $\CZFminus$.
	\end{itemize}
	Finally, $R\in B$ will denote a multi-valued function over the Heyting-valued universe $B^\Omega$ unless specified. \\
	
	The main goal of this subsection is proving the following preservation theorem: 
	\begin{theorem}\label{Theorem:BCSTset-preserving}
		Let $A$ be a BCST-regular set and $B \supseteq A$ a transitive model of $\CZFminus$ such that $\mathcal{P}(1)\cap A = \mathcal{P}(1)\cap B$. Then $B^\Omega$ thinks $\tilde{A}$ is BCST-regular.
	\end{theorem}
	
	Its proof requires a sequence of lemmas in a similar way to how we proved \autoref{Proposition:StrongCollectionValid}, the validity of Strong Collection over $V^\mathcal{S}$. The following lemma is an analogue of \autoref{Lemma:TargetSubsetMVChoice-uniform-overV}:
	
	\begin{lemma}\label{Lemma:TargetSubsetMVChoice}
		Let $a\in A$, $S \colon a\rrarrows A$ a multi-valued function and $Q\subseteq a\times A$ a class. Moreover, assume that
		\begin{enumerate}
			\item For each $x\in a$ there is $u\subseteq Q_x=\{z \mid \lag x,z\rag\in Q\}$ such that $\lag x,u\rag\in S$, and
			\item (Monotone Closure) If $\lag x,u\rag \in S$ and $u \subseteq v \subseteq Q_x$ then $\langle x,v\rangle\in S$.
		\end{enumerate}
		Then there is an $f\in A\cap \prescript{a}{}A$ such that $f(x)\subseteq Q_x$ and $\lag x,f(x)\rag \in S$ for all $x\in a$.
	\end{lemma}
	
	\begin{proof}
	    As before, consider the multi-valued function $S'$ with domain $a$ defined by
	        \begin{equation*}
	            S'(x,u) \qquad \text{if and only if} \qquad S(x,u)\text{ and } u\subseteq Q_x.
	        \end{equation*}
	   By \autoref{Lemma:SetMV}, there is $g\in A$ such that $g \colon a\rrarrows A$ and $g\subseteq S'$. Let $g_x=\{y\mid \lag x,y\rag \in g\}$, then $\bigcup g_x\subseteq Q_x$. Now take $f(x) = \bigcup g_x$, then $f\in A$ since $A$ satisfies Union and second-order Replacement. Moreover, since $S'$ is monotone closed, we have $\lag x,f(x)\rag\in S'$ for all $x\in a$.
	\end{proof}
	
	The following lemma is an analogue of \autoref{Lemma:StrongCollectionValid-MainLemma}. The reader is reminded that the proof of the following lemma necessarily uses the assumption $\mathcal{P}(1)\cap A=\mathcal{P}(1)\cap B$.
	
	\begin{lemma}\label{Lemma:Preserving-MainLemma}
        Let $a\in A^\Omega$ and $R\in B$.
	    Fix $p\in A$ such that $p=p^\dnot$ and $p\subseteq\llbr R \colon a\rrarrows \tilde{A}\rrbr^B$. If we define $P$ as 
	    \begin{equation*}
			P = \{\lag x,y,z\rag\in \dom a\times A^\Omega\times 1 \mid z\in (p\land a(x)\land \llbr \op(x,y) \in R\rrbr^B)\},
		\end{equation*}
		then there exists $r\in A$ such that $r\subseteq P$ and \mbox{$p\land a(x)\subseteq \{z\mid \exists y \in A^\Omega\ \lag x,y,z\rag \in r\}^\dnot$.}
	\end{lemma}
	
	\begin{proof}
	    Again, observe that $p\subseteq \llbr R \colon a\rrarrows \tilde{A}\rrbr^B$ is equivalent to
	    \begin{equation*}
		    p\land a(x)\subseteq \sideset{}{^B}\bigvee_{y\in A^\Omega} \llbr \op(x,y)\in R\rrbr^B = J^B \left(\bigcup_{y\in A^\Omega} \llbr \op(x,y)\in R\rrbr^B\right)
		\end{equation*}
		for all $x\in \dom a$.
		Now let us take 
		\begin{equation*}
		    Q=\{\lag x,z\rag\mid \exists y\in A^\Omega (\lag x,y,z\rag\in P)\},
		\end{equation*}
		then take $Q_x=\{z\mid \lag x,z\rag \in Q\}\subseteq 1$.
		We can easily see that \mbox{$Q_x=\bigcup_{y\in A^\Omega} \llbr\op(x,y)\in R\rrbr^B$} holds. 
		Thus we have $p\land a(x)\subseteq J^BQ_x$. Furthermore, it is also true that $Q_x\in \mathcal{P}(1)\cap B=\mathcal{P}(1)\cap A$. So we have $Q_x\in \mathcal{P}(1)\cap B${, which implies that} $J^B Q_x=Q_x^\dnot$.
		
		Now consider the relation $S\subseteq \dom a\times (\mathcal{P}(1)\cap A)$ defined by
		\begin{equation*}
			\lag x,u\rag \in S \quad\text{ if and only if }\quad u\subseteq Q_x\text{ and } p\land a(x)\subseteq u^\dnot.
		\end{equation*}
		We want to apply \autoref{Lemma:TargetSubsetMVChoice} to $S$, so we will check the hypotheses of \autoref{Lemma:TargetSubsetMVChoice} hold.
		
		The first condition holds because $\lag x,Q_x\rag \in S$ for each $x\in\dom a$. Furthermore, this shows that $S$ is a multi-valued function with domain $\dom a$. For the second condition, the relation is monotone closed in $A$ because if $v \subseteq w$ are in $\mathcal{P}(1) \cap A$ then $v^\dnot \subseteq w^\dnot$. 
		
		Therefore, by \autoref{Lemma:TargetSubsetMVChoice} applied to $S$, we have a function $f\in {^{\dom a}}A\cap A$ such that $p\land a(x)\subseteq f(x)^\dnot$ and $f(x)\subseteq Q_x$ for all $x\in\dom a$.
		Now let
		\begin{equation*}
			q = \{\lag x,z\rag \mid x\in\dom a\text{ and } z\in f(x)\}.
		\end{equation*}
		Then for each $\lag x,z\rag\in q$ there is $y\in A^\Omega$ such that $\lag x,y,z\rag \in P$ holds. By \autoref{Lemma:SetMV} applied to $P'\colon q\rrarrows A^\Omega$, defined by
		\begin{equation*}
		    P'(\lag x,z\rag, y) \qquad \text{if and only if} \qquad P(x,y,z),
		\end{equation*}
		there is $r\in A$ such that $r\subseteq P$ and $r \colon q\rrarrows A^\Omega$. It is easy to see that $r$ satisfies our desired property.
	\end{proof}
	
	\begin{remark}
	    There is a technical note for the proof of \autoref{Lemma:Preserving-MainLemma}, which is that there is no need for $P$, $Q$, and $Q_x$ to be definable over $A$ in general. The reason is that we do not know if either $R$ or $\llbr\cdot \rrbr^B$ are accessible from $A$. However, we do not need to {worry} about this since we are relying on the second-order Strong Collection over $A$.
	    
	    Furthermore, the proof of \autoref{Lemma:Preserving-MainLemma} also uses an assumption that $B$ is a transitive model of $\CZFminus$ implicitly. The reason is that we made use of Heyting operations relative to $B$, which we formulate over $\CZFminus$. Thus relativized Heyting operations are viable when $B$ satisfies $\CZFminus$. 
	\end{remark}

    We are now ready to prove the preservation theorem, \autoref{Theorem:BCSTset-preserving}. Its proof is parallel to that of \autoref{Proposition:StrongCollectionValid}.
    
	\begin{proof}[Proof of \autoref{Theorem:BCSTset-preserving}]
		First, observe that since $B$ is a transitive model of $\CZFminus$ both $B^\Omega$ and the Heyting operations over $B$ are well defined. Furthermore, it
		is easy to see that $B^\Omega$ thinks $\tilde{A}$ is transitive, and closed under Pairing, Union and Binary Intersection.
		Hence it remains to show that $B^\Omega$ thinks $\tilde{A}$ satisfies second-order Strong Collection, that is, 
		
		\begin{equation*}
			\llbr \forall a\in \tilde{A} \, \forall R [R \colon a\rrarrows \tilde{A}\to \exists b\in \tilde{A}(R \colon a\lrlrarrows b)] \rrbr^B = 1.
		\end{equation*}
		Take $a\in \dom\tilde{A}$, $R\in B^\Omega$ and $p\in B$ such that $p\subseteq 1$ and $p=p^\dnot$.
		
		We claim that if $p\subseteq \llbr R \colon a\rrarrows \tilde{A}\rrbr^B$, then there is $b\in \dom\tilde{A}$ such that $p\subseteq \llbr R \colon a\lrlrarrows b\rrbr^B$. Taking $P$ as in the statement of \autoref{Lemma:Preserving-MainLemma}, by \autoref{Lemma:Preserving-MainLemma}, we can find some $r\in A$ such that $r\subseteq P$ and
		
		\begin{equation*}
		    p\land a(x)\subseteq \{z\mid \exists y \lag x,y,z\rag\in r\}^\dnot.
		\end{equation*}
		
		Define $b$ such that $\dom b=\{y\mid\exists x,z(\lag x,y, {z} \rag\in r)\}$ and 
		\begin{equation*}
			b(y) =\{z\mid \exists x \lag x,y,z\rag\in r\}^\dnot.
		\end{equation*}
		
		for $y\in\dom b$. Note that $b(y)\in A$ since $r$ is. Then we can show that $p\subseteq \llbr R \colon a\lrlrarrows b\rrbr^B$ by following the computation given in the proof of \autoref{Proposition:StrongCollectionValid}.
	\end{proof}
	
	As a corollary, we have
	
	\begin{corollary}\label{Corollary:PreservationofPowerInaccessibles}
	    Let $K$ be a power inaccessible set. Then $V^\Omega$ thinks $\tilde{K}$ is power inaccessible.
	\end{corollary}
	
	\begin{proof}
	    Since $\mathcal{P}(1)\in K$, we have $\mathcal{P}(1)\cap K=\mathcal{P}(1)$. Thus \autoref{Theorem:BCSTset-preserving} shows $V^\Omega$ thinks $\tilde{K}$ is BCST-regular. It remains to show that $V^\Omega$ thinks $\tilde{K}$ is closed under the true power set.
	    
    	We claim that the argument for the preservation of Power Set given in the proof of \autoref{Theorem:CZFPersistence} relativizes to $K$. Work in $V$, and take $a\in K^\Omega$. Since $K$ is power inaccessible, the name defined by $\dom b =$ \mbox{$\prescript{\dom a}{} \{p^\dnot \mid p\subseteq 1\}$}, $b(c)=1$ is a member of $K$.
    	By the same calculation as in the proof of \autoref{Theorem:CZFPersistence}, we have \mbox{$\llbr\forall c (c\subseteq a\to c\in b)\rrbr = 1$.} Hence the desired result holds.
	\end{proof}
	
	\subsection{Interpreting an elementary embedding}
	In this subsection, we work over $\mathsf{CZF + BTEE}_M$ unless otherwise specified.
	
	We will show that elementary embeddings are persistent under Heyting-valued interpretation over $\Omega$.
	We will mostly follow Subsection 4.1.6 of Ziegler \cite{ZieglerPhD}, but we need to check his proof works in our setting since his applicative topology does not include Heyting algebras generated by formal topologies that are not set-presentable. Furthermore, we are working with a weaker theory than Ziegler assumed. Especially, we do not assume Set Induction, $\Delta_0$-Separation and Strong Collection for $(j,M)$-formulas, which calls for additional care.
	
	In most of our results in the rest of the paper, we will consider the case $M=V$, so considering the target universe, $M$, of $j$ will not be needed for our main results. Nevertheless, to begin with we will work in the more general setting and not assume that $M = V$.
	
	We need to define Heyting-valued interpretations for $M$ and $j$ in the forcing language.
	Since $j$ preserves names, we can interpret $j$ as $j$ itself. We will interpret $M$ by $M^\Omega$, as defined in \autoref{Definition:AOmega}, which was given by following the construction of $V^\Omega$ inside $M$.
	Thus defining $M^\Omega$ does not require Strong Collection or Set Induction for the language extended by $j$. The existence of $M^\Omega$ follows from the assumption that $M$ satisfied $\CZFminus$.
	
	One possible obstacle to defining the interpretation for the extended language is a non-absoluteness of $\Omega$ and the resulting double-negation formal topology between transitive models of $\CZFminus$.
	We discussed in \autoref{Subsection:DoubleNegationFormalTopology} that the Heyting interpretation $\llbr\phi\rrbr$ need not be absolute between transitive sets.
	
	It would be convenient if we have $\llbr\phi\rrbr=\llbr\phi\rrbr^M$, which follows from $\mathcal{P}(1)=\mathcal{P}(1)\cap M${, which thankfully we have.}
	
	\begin{lemma}
	    For any formula $\phi$ with parameters in $M^\Omega$, we have $\llbr \phi\rrbr^M = \llbr \phi^{M^\Omega}\rrbr$.
	\end{lemma}
	
	\begin{proof}
	   By \autoref{Lemma:PowerSetPreserving}, $\mathcal{P}(1)=\mathcal{P}(1)\cap M$. Hence the conclusion follows from the last clause of \autoref{Lemma:HeytingInterpretationBetweenModels}.
	\end{proof}
	
	Thus we do not need to worry about the absoluteness issue on the Heyting interpretation.
	Now we are ready to extend our forcing language to $\{\in, j, M\}$.
	
	\begin{definition}
	    Define $\llbr\phi\rrbr$ for the extended language as follows: 
	    \begin{itemize}
	        \item $\llbr j^m(a)\in j^n(b)\rrbr$ and $\llbr j^m(a)=j^n(b)\rrbr$ are defined in the same way as $\llbr x\in y\rrbr$ and $\llbr x=y\rrbr$ were.
	        \item $\llbr a\in M \rrbr := \bigvee_{x\in M^\Omega}\llbr a=x\rrbr$,
	        \item $\llbr\forall x\in M \phi(x)\rrbr := \bigwedge_{x\in M^\Omega} \llbr \phi(x)\rrbr$, and
	        \item $\llbr\exists x\in M \phi(x)\rrbr := \bigvee_{x\in M^\Omega} \llbr \phi(x)\rrbr$.
	    \end{itemize}
	\end{definition}
	
	\begin{remark}\label{Remark:InterpretationJFormulas}
	    The reader should be careful that $\llbr\phi(\vec{x})\rrbr$ is not in general a set. Provably over $\mathsf{CZF + BTEE}_M$, $\llbr j^m(a)\in j^n(b)\rrbr$ and $\llbr j^m(a)=j^n(b)\rrbr$ are always sets. In general, we can see that if $\phi(\vec{x})$ is a Heyting combination of atomic formulas, then $\llbr \phi(\vec{x})\rrbr$ is a set regardless of what the background theory is.
	    
	    However, taking a bounded quantification could make $\llbr\phi\rrbr$ a class that is not provably a set. For example, consider taking bounded universal quantification over $\llbr \phi(x.y)\rrbr$ for a Heyting combination of atomic $(j,M)$-formulas $\phi$. Then $\llbr\forall x\in a \phi(x,y) \rrbr = \bigwedge_{x\in \dom a} a(x)\land \llbr \phi(x,y)\rrbr$.
	    Each of $\llbr \phi(x,y)\rrbr$ can be a set, but the family $\{\llbr \phi(x,y)\rrbr \mid x\in\dom a\}$ need not be a set unless we have Strong Collection and $\Delta_0$-Separation for $(j,M)$-formulas.
	    
	    The situation is even worse for bounded existential quantifiers and disjunctions: we took a nucleus $\jmath$ when defining the interpretation of these, so $\llbr\exists x\in a\phi(x)\rrbr = \bigvee_{x\in\dom a}\llbr \phi(x)\rrbr =\jmath\left( \bigcup_{x\in\dom a}\llbr \phi(x)\rrbr\right)$.
	    This is ill-defined when the union $\bigcup_{x\in\dom a}\llbr \phi(x)\rrbr$ is a class rather than a set. Hence we have to use the join operator for classes instead of sets. The trade-off for the new definition of $\llbr\phi\rrbr$, in this case, is that we do not know if $\llbr\phi\lor\lnot\phi\rrbr=1$ for a $j$-formula $\phi$, unless we can ensure $\llbr\phi\rrbr$ is a set.
	\end{remark}
   
   	From this definition, we have an analogue of Lemma 4.26 of \cite{ZieglerPhD}, which is useful to check that $j$ is still elementary over $V^\Omega$:
	
	\begin{lemma}\label{Lemma:ValuePreserving}
		For any bounded formula $\phi(\vec{x})$ with all free variables displayed in the language $\in$ (that is, without $j$ and $M$), we have
		\begin{equation*}
			\llbr \phi(\vec{a}) \rrbr = \llbr \phi^{M^\Omega} (j(\vec{a})) \rrbr = \llbr \phi (j(\vec{a})) \rrbr
		\end{equation*}
		for every $\vec{a}\in V^\Omega$.
	\end{lemma}

	\begin{proof}
		For the first equality, note that $\llbr\phi(\vec{a})\rrbr\subseteq 1$. Hence we have
		\begin{equation*}
			\llbr \phi(\vec{a}) \rrbr = j(\llbr \phi(\vec{a}) \rrbr) = \llbr\phi(j(\vec{a}))\rrbr^M
			= \llbr \phi^{M^\Omega}(j(\vec{a})) \rrbr.
		\end{equation*}
		by \autoref{Lemma:PowerSetPreserving}. 
		
		The second equality will follow from the claim that for any bounded formula $\phi$ and $\vec{b} \in M^\Omega$, $\llbr \phi(\vec{b})\rrbr = \llbr\phi^{M^\Omega}(\vec{b})\rrbr$. The proof proceeds by induction on $\phi$: the atomic case and cases for $\land$, $\lor$, and $\to$ are trivial. For bounded $\forall$, observe that if $c,\vec{b}\in M^\Omega$ then $\llbr c = c \cap M^\Omega\rrbr = 1$, so $\llbr \forall x\in c \ \phi(x,\vec{b}) \lr (\forall x\in c \ \phi(x,\vec{b}) )^{M^\Omega}\rrbr=1$. This proves $\llbr \forall x\in c \ \phi(x,\vec{b})\rrbr = \llbr \forall x\in c \ \phi^{M^\Omega}(x,\vec{b})\rrbr$. The case for bounded existential quantifiers is similar.
	\end{proof}
	
	Moreover, we can check the following equalities easily:
	
	\bigskip
	
	\vbox{
	\begin{proposition} \phantom{a}
		\begin{enumerate}
			\item $\llbr \forall x,y (x=y\to j(x)=j(y))\rrbr=1$,
			\item $\llbr \forall x (j(x)\in M)\rrbr=1$,
			\item $\llbr \forall x (x\in M\to \forall y\in x(y\in M)) \rrbr=1$.
		\end{enumerate}
	\end{proposition}
	}
	
	\begin{proof}
		The first equality follows from $\llbr x=y\rrbr = \llbr j(x)=j(y)\rrbr$, which holds by the previous lemma, and the remaining two follow from direct calculations.
	\end{proof}
	
	\begin{lemma}\label{Lemma:ElementarityPersistentUnderDoubleNegationTranslation}
		For every $\vec{a}\in V^\Omega$ and formula $\phi$ that does not contain $j$ or $M$, we have \mbox{$\llbr \phi(\vec{a})\lr \phi^{M^\Omega}(j(\vec{a})) \rrbr=1$.}
	\end{lemma}

	\begin{proof}
		\autoref{Lemma:ValuePreserving} proved that this lemma holds for bounded formulas $\phi$.
		We will use full induction on $\phi$ to prove $\llbr\phi(\vec{a})\rrbr=\llbr\phi^M(j(\vec{a}))\rrbr$ for all $\vec{a}\in V^\Omega$.
		If $\phi$ is $\forall x\psi(x,\vec{a})$, we have $\llbr \forall x\psi(x,\vec{a})\rrbr=\bigwedge_{x\in V^\Omega}\llbr\psi(x,\vec{a})\rrbr$. Now
		
		\begin{align*}
			0\in \bigwedge_{x\in V^\Omega}\llbr \psi(x,\vec{a})\rrbr &\iff \forall x\in V^\Omega (0\in \llbr\psi(x,\vec{a})\rrbr)\\
			&\iff \forall x \in (V^\Omega)^M (0\in \llbr\psi(x,j(\vec{a}))\rrbr),
		\end{align*}
		
		where the last equivalence follows from applying $j$ to the above formula. Since the last formula is equivalent to $0\in \bigwedge_{x\in M^\Omega} \llbr\psi(x, j(\vec{a}))\rrbr$, we have 
        \begin{equation*}
            \bigwedge_{x\in V^\Omega} \llbr\psi(x, \vec{a})\rrbr=\bigwedge_{x\in M^\Omega} \llbr\psi(x, j(\vec{a}))\rrbr=\llbr \forall x\in M \psi^M(x, j(\vec{a}))\rrbr.
        \end{equation*}
		
		If $\phi$ is $\exists x\psi(x,\vec{a})$, we have $\llbr \exists x\psi(x,\vec{a})\rrbr=\bigvee_{x\in V^\Omega}\llbr\psi(x,\vec{a})\rrbr$. Moreover,
		
		\begin{align*}
			0\in \bigvee_{x\in V^\Omega}\llbr \psi(x,\vec{a})\rrbr 
			&\iff \exists p\subseteq 1 \left[ p\subseteq \bigcup_{x\in V^\Omega} \llbr\psi(x,\vec{a})\rrbr \text{ and } 0\in p^\dnot \right]
			\\&\iff \exists p\subseteq 1 \left[ p\subseteq \bigcup_{x\in M^\Omega} \llbr\psi^{M^\Omega}(x,j(\vec{a}))\rrbr \text{ and } 0\in p^\dnot \right]
		\end{align*}
		
		Hence $0\in \bigvee_{x\in V^\Omega}\llbr \psi(x,\vec{a})\rrbr$ if and only if $0\in \bigvee_{x\in M^\Omega}\llbr \psi(x,j(\vec{a}))\rrbr$
	\end{proof}
	
	We now work in the extended language $\CZF_{j,M}$. We need to check that $\Delta_0$-Separation and Strong Collection under the extended language are also persistent under the double-negation interpretation. We can see that the proof given by \cite{Gambino2006} and \autoref{Theorem:CZFPersistence} carries over, so we have the following claim:
	
	\begin{proposition}\pushQED{\qed} \label{Proposition:InterpretingJFormulaAxioms}
	    If $V$ satisfies any of Set Induction, Strong Collection, $\Delta_0$-Separation or Full Separation for the extended language, then the corresponding axiom for the extended language is valid in $V^\Omega$. \qedhere
	\end{proposition}
	
	The essential property of a critical set is that it is inaccessible. However, inaccessibility is not preserved under Heyting-valued interpretations in general. Fortunately, being a critical point is preserved provided it is regular:
	
	\begin{lemma}\label{Lemma:BeingACriticalPointIsPreserved}
		Let $K$ be a regular set such that $K\in j(K)$ and $j(x)=x$ for all $x\in K$.
		Then \mbox{$\llbr\tilde{K}\in j(\tilde{K})\land \forall x\in \tilde{K}(j(x)=x)\rrbr=1$.}
	\end{lemma}
	
	\begin{proof}
		Since $(\text{$j(K)$ is inaccessible})^M$, we have $j(K)\models\CZFminus$.
		Thus, by \autoref{Lemma:HeytingUniversePreserving}, $j(K)^\Omega = j(K)\cap V^\Omega$. By applying the same argument internal to $M$, we have $(j(K)^\Omega)^M = (j(K)\cap V^\Omega)^M = j(K)\cap M^\Omega$.
		Since $j(K)\cap V^\Omega \subseteq j(K)\subseteq M$, we have $j(K)\cap V^\Omega\subseteq M$. This implies $j(K)\cap V^\Omega = j(K)\cap V^\Omega\cap M = j(K)\cap M^\Omega$. In sum, we have
		\begin{equation*}
		    j(K)^\Omega = j(K)\cap V^\Omega = j(K)\cap M^\Omega = (j(K)^\Omega)^M,
		\end{equation*}
		and these are sets by \autoref{Lemma:HeytingUniversePreserving}.
		
		Also, $K\in j(K)$ implies $\tilde{K}\in j(K)$. Since the domain of $j(\tilde{K})=\widetilde{j(K)}$ is $j(K)\cap V^\Omega$, we have $\tilde{K}\in \dom j(\tilde{K})$, which implies $\llbr\tilde{K}\in j(\tilde{K})\rrbr=1$.
		For the assertion $\llbr\forall x\in \tilde{K} (j(x)=x)\rrbr=1$, observe that if $x\in \dom\tilde{K}$ then $j(x)=x$, so we have the desired conclusion.
	\end{proof}
	
	\subsection{Consistency strength: intermediate results}
	By using the results from the previous sections and subsections, we have the following:
	
	\begin{theorem}\label{Theorem:InterpretationResults} \pushQED{\qed} \phantom{a} 
	    \begin{enumerate}
	        \item $\IZF+\BTEE$ interprets $\ZF+\BTEE$.
	        \item $\IZF+\BTEE+\text{Set Induction}_j$ interprets $\ZF+\BTEE+\TIj$.
	        \item $\IZF+\WA$ interprets $\ZF+\WA$. \qedhere 
	    \end{enumerate}
	\end{theorem}
	
	Hence
	
	\begin{corollary} \label{Corollary:MainConsistencyResults} \phantom{a} 
	    \begin{enumerate}
	        \item $(\IKP_{j,M})$ $\IKP$ with a critical point implies $\mathsf{Con(ZF + BTEE + \TIj)}$.
	        
	        \item $\CZF$ with a Reinhardt set implies $\mathsf{Con(ZF+WA)}$.
	    \end{enumerate}
	\end{corollary}
	
	\begin{proof}
	    Using \autoref{Theorem:IKPSigmaOrdimpliesIZF+BTEE}, if $K$ is a critical point of an embedding $j \colon V \rightarrow M$ then we can find some ordinal $\lambda$ for which $\lag L_\lambda, j \restricts L_\lambda \rag$ is a model of $\mathsf{IZF + BTEE}+\text{Set Induction}_j$. Therefore the first claim follows from the first claim of \autoref{Theorem:InterpretationResults}. The second claim follows directly from \autoref{Theorem:ReinhardtCriticalPtModelsWA}.
	\end{proof}
	
	By mimicking Bagaria-Koellner-Woodin's argument in \cite{BagariaKoellnerWoodin2019} over $V^\Omega$, we have the following consistency result:
	
	\begin{theorem}[$\CGB_\infty$]\label{Theorem:LowerBound-superReinhardt}
	    Let $K$ be a super Reinhardt set. Then $V^\Omega$ satisfies $\mathsf{ZF}$ plus there is a proper class of inaccessible cardinals $\gamma$ such that $(V_\gamma,V_{\gamma+1})\models \mathsf{ZF_2} + \text{there is a Reinhardt cardinal}$.
	\end{theorem}

	\begin{proof}
	    The proof will proceed as follows: First, we will show that the background theory interprets some moderate semi-intuitionistic theory. Then we will derive that this semi-intuitionistic theory proves that there is a proper class of inaccessible cardinals $\gamma$ such that $V_\gamma$ is a model of $\mathsf{ZF}$ with a Reinhardt cardinal.
	
	    So, let $K$ be a super Reinhardt set. For any $a\in V^\Omega$, we can find an amenable elementary embedding $j\colon V\to V$ with critical set $K$ such that $a\in j(K)$. As per usual, let $\Lambda = j^{\omega}(K)$. We shall
	    restrict our background theory to its first-order part to facilitate the proof.
	    
	    By \autoref{Prop:SuperReinhardtPowerInaccessible} and \autoref{Corollary:SuperReinhardtmodelsIZFpIEA}, $K$ is power inaccessible and $V$ satisfies $\mathsf{IZF+pIEA}$.
	    Thus by \autoref{Corollary:VOmegamodelsZFminus} and \autoref{Corollary:PreservationofPowerInaccessibles}, $V^\Omega$ interprets the following statements:
	    \begin{itemize}
	        \item Axioms of $\mathsf{ZF}$, and
	        \item There is a proper class of inaccessible cardinals.
	    \end{itemize}
	    Especially, $V^\Omega$ interprets the law of excluded middle for formulas without $j$. 
	    However, we do not know if $V^\Omega$ satisfies $\mathsf{ZF}_j$ because there is no reason why $V$ should satisfy the Separation scheme for $j$-formulas. As a result, the double-negation translation does not force the law of excluded middle for $j$-formulas.
	    Despite this, $V^\Omega$ still believes the following statements are valid due to \autoref{Lemma:ElementarityPersistentUnderDoubleNegationTranslation} and \autoref{Proposition:InterpretingJFormulaAxioms}:
	    
	    \begin{itemize}
	        \item $j$ is amenable and elementary, and 
	        \item Collection and Set Induction for $j$-formulas.
	    \end{itemize}
	    
	    Amenability of $j$ needs some justification. Working over $V$, $j$ is amenable by \autoref{Lemma:FunctionalClass-amenable}. As we observed in \autoref{Remark:SuperReinhardtsareReinhardts}, if $j$ is amenable then Separation for $\Delta_0^j$-formulas holds. By \autoref{Proposition:InterpretingJFormulaAxioms}, $\Delta_0$-Separation for $j$-formulas is valid in $V^\Omega$ from which it follows that $V^\Omega$ thinks $j$ is amenable.
	    
	    Now work in $V^\Omega$. Since $V^\Omega$ validates $\mathsf{ZF}$ with a proper class of inaccessible cardinals, we can find the least inaccessible cardinal $\gamma$ such that $\gamma > j^\omega(\kappa)$, where $\kappa=\rank \tilde{K}$.
	    Here $j^\omega(\kappa)$ is well-defined since the sequence $\langle j^n(\kappa)\mid n\in\omega\rangle$ required Set Induction for $j$-formulas for its definition and Collection for its supremum to exist.
	    Furthermore, one can see that $V^\Omega$ thinks $\tilde{K}$ is a critical point of $j$ and $\tilde{K}=V_\kappa$. (The latter equality follows from the fact that $\tilde{K}$ is power inaccessible.) From which it follows that $\kappa=\crit j$.
	    
	    Since $\gamma$ is definable from the parameter $j^{\omega}(\kappa)$, which is fixed by $j$, by elementarity we have $\gamma=j(\gamma)$.
	    Moreover, $V^\Omega$ also believes $j\restricts V_\gamma\in V_{\gamma+1}$ and $\crit (j\restricts V_\gamma)=\kappa$. (Amenability of $j$ ensures $j\restricts V_\gamma$ exists.) Hence $V^\Omega$ thinks $(V_\gamma,V_{\gamma+1})$ is a model of $\mathsf{ZF_2}$ with a Reinhardt cardinal.
	    
	    Finally, recall that $a\in j(K)$ and $a$ were arbitrary. Hence we have proved that for every $a\in V^\Omega$, $V^\Omega$ thinks there is an inaccessible cardinal $\gamma$ such that $(V_\gamma,V_{\gamma+1})$ believes that there is a Reinhardt cardinal $\kappa$ for which $a\in V_{j(\kappa)}$. Hence $V^\Omega$ thinks there is a proper class of such $\gamma$. 
	\end{proof}
	
	\begin{remark} \label{Remark:DoubleNegationvsTopology}
	    Most results in this section can be obtained by using Friedman's double-negation translation as defined in \cite{Friedman1973} and \cite{FriedmanScedrov1984}. In some cases, applying Friedman's argument would be simpler: for example, verifying the axioms of $\IZF$ under Friedman's double-negation translation does not involve any lengthy proof, unlike the verification of Strong Collection over Gambino's Heyting universe $V^\mathcal{S}$.
	    However, Friedman's proof of the validity of Collection under this translation heavily relies on Full Separation, an axiom scheme that $\CZF$ does not enjoy.
	    
	    Thus we may ask what the advantages are to using Gambino's Heyting-valued interpretation over Friedman's double-negation interpretation. The main reason is that Gambino's presentation is closer to forcing, which is a technique more familiar to set theorists. Secondly, as an intermediate step, Friedman first interprets a non-extensional set theory. This can be very notationally heavy and difficult to follow when one first tries to understand the arguments. Finally, most of this work has been done over the weak system of $\CZF$ and it is unknown how to achieve Friedman's interpretation in this system.
	    
	    On the other hand, there may be no advantage to our main results concerning critical sets and Reinhardt sets. This is because in Theorems \ref{Theorem:IKPSigmaOrdimpliesIZF+BTEE} and \ref{Theorem:ReinhardtCriticalPtModelsWA} we obtained a lower bound for their consistency strength in terms of $\IZF$ plus some large set axioms. One could then take Freidman's double-negation translation of this theory to obtain the consistency of $\ZF + \BTEE + \TIj$ and $\mathsf{ZF + WA}$ for critical sets and Reinhardt sets respectively. However, using Gambino's method we can strengthen this result to preserve our background theory while containing a set model of this classical theory. Namely, we will see that $V^\Omega$ preserves the theory $\CZFminus_j$ and contains, as a set, $\tilde{\Lambda}$ which is a model of the above theory. Since Friedman's interpretation does not work over weak theories, it does not seem to be possible to obtain a similar result with that translation.
    \end{remark}
	
	\begin{theorem} \label{Theorem:CZF+ReinhardtinterpretsZF+WA}
	    Working with the theory $\CZF_{j,M}$ with a critical set, $V^\Omega$ validates $\CZFminus_j + \Delta_0\mhyphen \mathsf{LEM}$ with a critical point $\tilde{K}$ of $j$. Furthermore, $\tilde{\Lambda} = j^\Omega(\tilde{K})$ satisfies $\ZF + \BTEE$ plus Set Induction for $j$-formulas.
	    If we strengthen a critical set to a Reinhardt set, then $V^\Omega$ validates $\tilde{\Lambda}\models \ZF+\WA$.
	    
	    If we add Full Separation into the background theory, then $V^\Omega$ satisfies not only $\CZFminus_j + \Delta_0\mhyphen \mathsf{LEM}$, but also $\ZFminus$.
	\end{theorem}
	
	\begin{proof}
	    The results in the previous subsections show $V^\Omega$ validates $\CZFminus_j + \Delta_0\mhyphen \mathsf{LEM}$, and that $j$ is an elementary embedding from $V^\Omega$ to $M^\Omega$ with a critical point $\tilde{K}$.
	    Since $K$ is BCST-regular, \mbox{$K\in j(K)$}, and $\mathcal{P}(1)\cap K = \mathcal{P}(1)\cap j(K)$, we can apply \autoref{Theorem:BCSTset-preserving} to $K$ and $j(K)$. Hence we have that $\llbr \tilde{K}\text{ is BCST-regular}\rrbr^{j(K)^\Omega}=1$, which implies that $\llbr \tilde{K}\models \CZFminus\rrbr^{j(K)^\Omega}=1$.
	    Observe that the claim $\tilde{K}\models \CZFminus$ is $\Delta_0$, so applying the second clause of \autoref{Lemma:HeytingInterpretationBetweenModels} proves $\llbr\tilde{K}\models \CZFminus\rrbr=1$.
	    
	    Now working in $V^\Omega$, the excluded middle for bounded formulas gives us that every transitive set satisfies the full excluded middle. Especially, both $\tilde{K}$ and $\tilde{\Lambda}$ satisfy the full excluded middle. In addition, $\tilde{\Lambda}$ satisfies $\IZF+\BTEE$ plus Set Induction for $j$-formulas by \autoref{Corollary:LambdaModelsIZFBTEEInd}. Thus we have the desired result.
	    
	    The case for a Reinhardt set is analogous to the previous one, except that we apply \autoref{Theorem:ReinhardtCriticalPtModelsWA} to show that $\tilde{\Lambda}$ validates $\IZF+\WA$.
        
        Finally, if we have Full Separation, then $V^\Omega$ also satisfies Full Separation which completes the proof since the combination of Full Separation and $\Delta_0\mhyphen \mathsf{LEM}$ implies the full excluded middle.
	\end{proof}
	
    We will end this section by observing that the translation preserves the cofinality of an elementary embedding over a moderate extension of $\CZF$.
    None of the previous analysis has required either Full Separation or Subset Collection in the background universe. On the other hand, the following proof, which we include for completeness, requires either Full Separation or $\mathsf{REA}$.
    
	\begin{lemma}[$\CZF_{j}$]\label{Lemma:CofinalityPreserved}
	    Assume either Full Separation or $\mathsf{REA}$. 
		If $j \colon V\prec V$ is a cofinal elementary embedding, then $V^\Omega$ thinks $j$ is cofinal.
	\end{lemma}
	
	\begin{proof}
		Let $a\in V^\Omega$. Then there is a set $X$ such that $a\in j(X)$. If we assume Full Separation, then $X\cap V^\Omega$ is a set, and $j(X\cap V^\Omega)=j(X)\cap V^\Omega$.
		Let $b$ be a name such that $\dom b= X\cap V^\Omega$ and $b(y)=1$ for all $y\in\dom b$.
		Then $\llbr a\in j(b) \rrbr=\top$.
		
		An additional step is required if we assume $\mathsf{REA}$ instead: Take a set $X$ such that $a\in j(X)$. By $\mathsf{REA}$, we can find a regular set $Y$ such that $X\in Y$. By \autoref{Lemma:HeytingUniversePreserving}, $Y^\Omega=Y\cap V^\Omega$ is a set. The remaining argument is then identical to the previous one.
	\end{proof}
	
	While \autoref{Lemma:CofinalityPreserved} does not directly suggest anything about consistency strength, when combined with \autoref{Theorem:CZF+ReinhardtinterpretsZF+WA} it does tell us that $\mathsf{CZF+Sep}$ with a Reinhardt set interprets $\ZFminus$ with the existence of a non-trivial cofinal embedding \mbox{$j\colon V\to V$.}
	
	As we pointed out after \autoref{Proposition:CZFReinhardtEmbeddingCofinal}, over $\ZFCminus$ with the $\mathsf{DC}_\mu$-schemes for all cardinals $\mu$ there cannot be a non-trivial cofinal embedding $j\colon V\to V$.
	The $\mathsf{DC}_\mu$-schemes are a variant of the Axiom of Choice and adding these schemes to $\ZFminus$ does not bolster its consistency strength. In fact, $\ZFminus$ proves $L$ satisfies the $\mathsf{DC}_\mu$-scheme for every cardinal $\mu$ in $L$. Thus the absence of a cofinal embedding over $\ZFCminus$ with the $\mathsf{DC}_\mu$-schemes for any cardinal $\mu$ can be seen as a variant of the Kunen inconsistency phenomenon. 
	
	However, it is unclear how we can use the above results to obtain a stronger bound for the consistency strength of $\CZF$ with a Reinhardt set. It is known that if one extends $\CZF$ by the \emph{Relation Reflection Scheme} ($\mathsf{RRS}$), as defined by Aczel \cite{Aczel2008}, then this is also persistent under Gambino's Heyting-values interpretation and $\ZFCminus$ proves that $\mathsf{RRS}$ is equivalent to the $\mathsf{DC}$-Scheme. It is further possible that one might be able to generalize such a scheme to $\mathsf{DC}_\mu$ for larger cardinals, however, this does not appear to be sufficient to derive an inconsistency due to the heavy use of Well-Ordering in the proof of inconsistency in \cite{Matthews2020}.
	
	Alternatively, it is proven in \cite{Matthews2020} that we can remove the assumption of the Dependent Choice Schemes if we instead require that $V_{\crit j} \in V$ (in fact, the main reason one assumes the $\mathsf{DC}_\mu${-scheme} for every cardinal $\mu$ is to prove this). The issue is that we do not know whether $V^\Omega$ believes that $V_{\crit j} \in V$, even if we assume $V$ satisfies Full Separation or $\mathsf{REA}$. If we were to assume that there was a Reinhardt set which was Power Inaccessible then this would be the case, however, we can only obtain that $\Lambda$ believes that $K$ is Power Inaccessible which is insufficient to derive the required result.

	\section{Double-negation translation of second-order set theories}\label{Section:DNT-SOST}

	In this section, we will provide a double-negation translation of $\IGB$ and $\TR$. One technical issue is that the statement of $\TR$ requires an infinite conjunction. However, we know that the full elementarity of $j \colon V\to V$ is definable in a classical context, namely $\Sigma_1$-elementarity, and that this can be codified into a single formula by using the partial truth predicate for $\Sigma_1$-formulas. Thus we will not try to interpret infinite connectives.
	
	\subsection{Interpreting \texorpdfstring{$\mathsf{GB}$}{GB} from \texorpdfstring{$\IGB$}{IGB}}
	
	\begin{definition}
	    A class $A$ is an \emph{$\Omega$-class name}, or simply a \emph{class name} if:
	    \begin{itemize}
	        \item The elements of $A$ are of the form $\langle a, p\rangle$, where $a\in V^\Omega$, $p\subseteq 1$ and $p^\dnot = p$,
	        \item It is \emph{functional} in the sense that $\langle a,p\rangle, \langle a,q\rangle\in A$ implies $p=q$.
	    \end{itemize}
	    
	    If $A$ is a class name, define $\dom A := \{a \mid \exists p \ \langle a,p\rangle\in A\}$ and $A(a)=p$ for the unique $p$ such that $\langle a,p\rangle\in A$.
	\end{definition}
	
	Then we can extend $\llbr\phi \rrbr$ to atomic second-order formulas as follows:
	\begin{definition}
	    Let $a\in V^\Omega$ and $A$, $B$ be class names. 
	    Define $\llbr a\in A\rrbr := \bigvee_{x\in\dom A}A(x)\land \llbr x=a\rrbr$ and 
	    \begin{equation*}
	        \llbr A=B\rrbr :=\left( \bigwedge_{x\in\dom A} A(x) \to\bigvee_{y\in\dom B} B(y)\land \llbr x=y\rrbr \right) \land 
	        \left( \bigwedge_{y\in\dom B} B(y) \to\bigvee_{x\in\dom A} A(x)\land \llbr x=y\rrbr \right).
	    \end{equation*}
	\end{definition}
	
	Based on the above definition, we can extend our interpretation to any formula $\phi$ with no second-order quantifiers. 
	
	Now we want to extend the double-negation interpretation to formulas which have second-order quantifiers, but this calls into question how this should be formulated in our second-order context.
	To illustrate the problem, observe that we take a class meet and join to interpret unbounded first-order quantifiers, and the resulting $\llbr\phi\rrbr$ is a class rather than a set.
	Thus we may suspect that interpreting second-order quantifiers results in $\llbr\phi\rrbr$ being a \emph{hyperclass}, a collection of classes, which is clearly not an object of $\CGB$ or $\IGB$.
	
	This situation is analogous to the one encountered in \autoref{Remark:InterpretationJFormulas}, where $\llbr\phi\rrbr$ for a $j$-formula $\phi$ may not be a set. This was because the definitions of $\lor$ and $\exists$ for $j$-formulas used $J$ instead of $\jmath$ and we only know that $J(\llbr\phi\rrbr\cup\lnot\llbr\phi\rrbr)=1$ is true when $\llbr\phi\rrbr$ is a set.
	
	We may suspect the same issue happens for second-order formulas: we may need to extend $\jmath$ to some function $\mathcal{J}$ which is defined for hyperclasses, as $J$ serves as an extension of $\jmath$, to define $\lor$ and $\exists$ for second-order formulas. However, even if we have a proper definition of $\mathcal{J}$, there is no reason to believe that $\mathcal{J}(\llbr\phi\rrbr\cup\lnot\llbr\phi\rrbr)=1$ should hold for a second-{order} formula $\phi$.
	
	The situation would be better if we can ensure $\llbr\phi\rrbr$ is a set. We know that if we additionally assume Full Separation then, for $\phi$ a first-order formula, $\llbr\phi\rrbr$ is a set. Analogously, we may suspect that we can prove $\llbr\phi\rrbr$ is a set when $\phi$ is a second-order formula, if we have Full Separation and the \emph{Full Class Comprehension scheme}, which is the assertion that $\{x\mid\phi(x)\}$ is a class for any class formula $\phi$. However, adding this axiom to $\mathsf{GBC}$ results in Kelley-Morse theory, $\mathsf{KM}$, which would considerably increase the strength of the underlying theory we wish to work with.

	Because we do not have a good framework for working with hyperclasses we will circumvent this issue in another way, by combining Gambino's interpretation with Friedman's double-negation translation. We will take Gambino's interpretations for formulas with no second-order quantifiers, and extend it using the same translation as in Friedman's double-negation interpretation.
	
	\begin{definition}
	    For a second-order formula $\phi$ with set parameters from $V^\Omega$ and $\Omega$-class name parameters, define $\phi^-$ inductively as follows:
	    \begin{itemize}
	        \item If $\phi$ is an atomic formula, then $\phi^- \equiv (\llbr\phi\rrbr = 1)$,
	        \item $(\phi\land\psi)^- \equiv \phi^-\land \psi^-$,
	        \item $(\phi\lor\psi)^- \equiv \lnot\lnot(\phi^-\lor \psi^-)$,
	        \item $(\phi\to\psi)^- \equiv \phi^-\to \psi^-$,
	        \item $(\forall^0 x \phi(x))^- \equiv \forall^0 x\in V^\Omega\ \phi^-(x)$,
	        \item $(\exists^0 x \phi(x))^- \equiv \lnot\lnot\exists^0 x\in V^\Omega\ \phi^-(x)$,
	        \item $(\forall^1 X \phi(X))^- \equiv \forall^1 X ((\text{$X$ is a class name})\to \phi^-(X))$, 
	        \item $(\exists^1 X \phi(X))^- \equiv \lnot\lnot \exists^1 X ((\text{$X$ is a class name})\land \phi^-(X))$.
	    \end{itemize}
	\end{definition}
	
	\begin{lemma}[$\IGB$]\label{Lemma:TwoInterpretationsCoindice}
	    If $\phi$ is a formula with no second-order quantifiers, then $\phi^-\equiv (\llbr\phi\rrbr=1)$.
	\end{lemma}
	
	\begin{proof}
	    The proof proceeds by induction on $\phi$. Atomic cases follows by definition. For conjunction, 
	    \begin{align*}
	        (\phi\land \psi)^- \equiv \phi^- \land \psi^- \iff \llbr\phi\rrbr = 1\text{ and } \llbr\psi\rrbr = 1 \iff \llbr\phi \rrbr\cap \llbr\psi\rrbr = 1
	    \end{align*}
	    and we can see that $\llbr\phi\rrbr\cap\llbr\psi\rrbr = \llbr\phi\cap\psi\rrbr$ by the definition of our translation. The cases for bounded and unbounded $\forall$ and implications are analogous.
	    
	    We need some care for the case for disjunction and $\exists$ and in particular we will need to make use of both Full Separation and Powerset. We examine the case for unbounded $\exists$, to see why we need these two axioms.
	    We know that $\llbr\exists x \phi(x)\rrbr = \bigvee_{x\in V^\Omega} \llbr \phi(x)\rrbr$, and this is $J\left(\bigcup_{x\in V^\Omega} \llbr \phi(x)\rrbr\right)$. By Full Separation, $\bigcup_{x\in V^\Omega} \llbr \phi(x)\rrbr = \{0\mid \exists x\in V^\Omega \llbr \phi(x)\rrbr=1\}$ is a set. Furthermore, Powerset proves $Jp=p^\dnot$ for all $p\subseteq 1$.
	    Hence $\llbr\exists x\phi(x)\rrbr = \left(\bigcup_{x\in V^\Omega} \llbr \phi(x)\rrbr\right)^\dnot$, so
	    
	    \begin{align*}
	        (\exists x\phi(x))^- & \equiv \dnot\exists x\in V^\Omega \phi^-(x) \iff
	        \dnot\exists x\in V^\Omega(\llbr\phi(x)\rrbr =1) \\ &
	        \iff 0 \in \{ 0 \mid \dnot ( \exists x \in V^\Omega ( \llbr \phi(x) \rrbr = 1)) \} 
	        \iff \left(\bigcup_{x\in V^\Omega} \llbr \phi(x)\rrbr\right)^\dnot = 1.
	    \end{align*}
	    The cases for bounded existential quantifiers and disjunctions are analogous, so we omit them.
	\end{proof}
	
	\begin{lemma}\label{Lemma:ClassEqualityRespectsLogicalFormulas}
	    Let $A$ and $B$ be class names and $\phi(X)$ be a formula of second-order set theory. Then
	    \begin{equation*}
	        ( \llbr A=B\rrbr=1 \land \phi^-(A) ) \to \phi^-(B).
	    \end{equation*}
	\end{lemma}
	
	\begin{proof}
	    The proof proceeds by induction on $\phi$, and one can see that the only non-trivial part of this is the atomic case where $\phi(X)$ is $x\in X$.
	    
	    To do this, we show that $\llbr A=B\rrbr \land \llbr a\in A\rrbr \le \llbr a\in B\rrbr$:
	    \begin{align*}
	        \llbr A=B\rrbr \land \llbr a\in A\rrbr &\le \left( \bigwedge_{x\in\dom A} A(x) \to\bigvee_{y\in\dom B} B(y)\land \llbr x=y\rrbr \right) \land\left(\bigvee_{x\in\dom A}A(x)\land \llbr x=a\rrbr\right)\\
	        & \le \bigvee_{x\in\dom A} \left( \left(\bigvee_{y\in\dom B} B(y)\land \llbr x=y\rrbr \right) \land \llbr x=a\rrbr \right)\\
	        & \le \bigvee_{y\in\dom B} B(y)\land \llbr y=a\rrbr = \llbr a\in B\rrbr.
	    \end{align*}
	    Therefore, if both $\llbr A=B\rrbr=1$ and $\llbr a\in A\rrbr=1$ then $\llbr a \in B \rrbr = 1$.
	\end{proof}
	
	\begin{theorem}[$\IGB$]
	    Working over $\IGB$, every axiom of $\mathsf{GB}$ is valid in $V^\Omega$.
    \end{theorem}

	\begin{proof}
	By \autoref{Theorem:CZFPersistence} and \autoref{Lemma:TwoInterpretationsCoindice}, we can see that the first-order part of $\IGB$ is valid. Moreover, $(\phi\lor\lnot\phi)^-$ is valid since it is equivalent to $\lnot\lnot(\phi^-\lor\lnot\phi^-)$, which is constructively valid.
	    
	    It remains to show that the second-order part of $\IGB$ is valid under the interpretation.
	    \begin{itemize}
	        \item Class Extensionality: follows from \autoref{Lemma:ClassEqualityRespectsLogicalFormulas}.
	        \item Elementary Comprehension: Let $a\in V^\Omega$, $A$ be a class name, and $\phi$ be a formula without class quantifiers. Then $\llbr\phi(x,a,A)\rrbr$ is well-defined for $x\in V^\Omega$. Now consider $B$ to be the class name
	        \begin{equation*}
	            B = \{\langle x,\llbr\phi(x,a,A)\rrbr\rangle \mid x\in V^\Omega\}.
	        \end{equation*}
	        
	        Furthermore, we can easily see that $\llbr x\in B\lr \phi(x,a,A)\rrbr=1$. Hence $B$ witnesses Elementary Comprehension for $\phi(x,a,A)$.
	        
	        \item Class Set Induction: the usual argument for first-order Set Induction carries over. 
	        \item Class Strong Collection: we can see that the proof of \autoref{Proposition:StrongCollectionValid} works if we replace $R(x,y)$ with $\lag x,y\rag \in R$.
	    \end{itemize}
	
	    We do not need to check that the double-negation translation interpretation also validates Class Separation since $\mathsf{GB}$ without Class Separation proves Class Separation: it follows from Class Replacement, and the proof is similar to the derivation of Separation from Replacement over classical set theory.
	\end{proof}
	
	\subsection{Interpreting \texorpdfstring{$\TR$}{TR}}
	This subsection is devoted to the following result:
	
	\begin{theorem}\label{Theorem:IGBTRinterprets}
	   If $\mathsf{GB+TR}\vdash\phi$, then $\mathsf{IGB_\infty+TR}\vdash\phi^-$.
	\end{theorem}
	
	As mentioned before, we will dismiss infinite connectives from the interpretation. The main reason is that in $\mathsf{GB}$ every $\Sigma^A_1$-elementary embedding is fully $A$-elementary embedding by essentially \autoref{Lemma:DeltaA0elementarityandfullelementarity} (2).

	\begin{proof}[Proof of \autoref{Theorem:IGBTRinterprets}]
	    It suffices to show that $\mathsf{IGB_{\infty}+TR}$ proves $(\TR)^-$.
	    Let $A$ be a class name and $a\in V^\Omega$. We claim that there is a class name $\tilde{\jmath}$ such that 
	    \begin{equation*}
	        \llbr \tilde{\jmath}\text{ is $A$-elementary}\land a\in \tilde{\jmath}({\tilde{K}})\land {\tilde{K}} \text{ is a critical point of } {\tilde{\jmath}} \rrbr = 1.
	    \end{equation*}
	    
	    By $\TR$, we can find some elementary embedding $j \colon V\to V$ and an inaccessible set $K$ such that $j$ is $A$-elementary, $a\in j(K)$ and $K$ is a critical point of $j$.
	    Define $\tilde{\jmath} = \{\langle \op(a,j(a)),1\rangle \mid a\in V^\Omega\}$. It is easy to see that $\llbr\tilde{\jmath}\text{ is a function}\rrbr=1$, and, by \autoref{Lemma:BeingACriticalPointIsPreserved}, $\llbr \tilde{K}\text{ is a critical point of }\tilde{\jmath}\rrbr=1$. Furthermore, $\llbr a\in\tilde{\jmath}(\tilde{K})\rrbr=1$ is equivalent to $\llbr a\in j(\tilde{K})\rrbr=1$, the latter of which follows from $a\in j(K)$.
	    
	    It remains to show that $V^\Omega$ thinks $\tilde{\jmath}$ is an $A$-elementary embedding. For this, it suffices to show that $\tilde{\jmath}$ is $\Sigma^A$-elementary by \autoref{Lemma:DeltaA0elementarityandfullelementarity}. Let $\psi(x)$ be a $\Sigma^A$-formula. Then
	    \begin{equation*}
	        {\bigllbr}\forall \vec{x} {\big(} \psi(\vec{x}) \lr\psi(j(\vec{x})) {\big)} {\bigrrbr} =1
	    \end{equation*}
	    is equivalent to
	    \begin{equation*}
	        \forall \vec{x} \in V^\Omega \big( \llbr \psi(\vec{x})\rrbr = 1 \lr \llbr \psi(j(\vec{x}))\rrbr = 1 \big).
	    \end{equation*}
	    
	    Finally, since $\llbr \psi(\vec{x}) \rrbr = 1$ is expressible in $V$ using a $\Sigma^A$-formula (without any other class parameters), the above formula is immediate from the $A$-elementarity of $j$.
    \end{proof}
	
	Combining the above analysis with \autoref{Theorem:CGBTRinterpretsIGBTR}, the following corollary is immediate.
	
	\begin{corollary}\pushQED{\qed} \label{Corollary:CGBTRinterpretsGBTR}
        $\mathsf{CGB_\infty+TR}$ interprets $\mathsf{GB+TR}$.
	    \qedhere
	\end{corollary}
	
    \begin{remark}
	    Some readers may wonder about whether we need full elementarity in the formulation of $\TR$, because the proof of \autoref{Theorem:IGBTRinterprets} would work if $j$ preserves formulas of some bounded complexity, probably $\Sigma^A$-formulas.
	    It is correct that what the proof of \autoref{Theorem:IGBTRinterprets} actually shows is $\mathsf{IGB_{\infty}}$ with a statement weaker than $\TR$ interprets $\mathsf{GB+TR}$.
	    However, the full power of $\TR$ was necessary in \autoref{Theorem:CGBTRinterpretsIGBTR} to derive $\mathsf{IGB_{\infty}+TR}$ from $\mathsf{CGB_{\infty}+TR}$. 
	\end{remark}

	\section{Consistency strength: final results} \label{Section:ConsistencyStrength:Final}
    We have proven in \autoref{Corollary:MainConsistencyResults} that $\IKP$ with a critical point implies the consistency of $\ZF + \BTEE + \TIj$ while $\CZF$ with a Reinhardt set implies the consistency of $\mathsf{ZF + WA}$. However, one may ask how strong these notions are with respect to the traditional large cardinal hierarchy over $\ZFC$, which is a question we address here. 
	
	% ZF + BTEE + TI_j
	Let us examine the theory $\ZF + \BTEE + \TIj$ first. We can see that $\ZF + \BTEE + \TIj$ proves $L$ satisfies the same theory. Hence we have the consistency of $\ZFC + \BTEE + \TIj$ from that of $\ZF + \BTEE + \TIj$.

	Let us compare the consistency strength of $\ZFC + \BTEE + \TIj$ with that of other large cardinal axioms to illustrate how strong it is.
	We can see that \emph{virtually rank-into-rank cardinals}, defined by Gitman and Schindler \cite{GitmanSchindler2018}, provide an upper bound:
	
	\begin{definition}[\cite{GitmanSchindler2018}]
	    A cardinal $\kappa$ is \emph{virtually rank-into-rank} if in some set-forcing extension it is the critical point of an elementary embedding $j\colon V_\lambda\to V_\lambda$ for some $\lambda > \kappa$.
	\end{definition}
	
	\begin{lemma}
	    Let $\kappa$ be a virtually rank-into-rank cardinal and $\lambda$ witness this. If $j\colon V_\lambda\to V_\lambda$ is a elementary embedding over a set-generic extension with $\crit j=\kappa$, then, in this extension, $(V_{j^\omega(\kappa)},\in j)$ satisfies $\mathsf{ZFC+BTEE} + \TIj$.
	\end{lemma}
	
	\begin{proof}
	    First, $\kappa$ is inaccessible in $V$. This follows from Theorem 4.20 of \cite{GitmanSchindler2018} and known facts about $\omega$-iterable cardinals, but we shall also give a direct proof for it. Suppose, for a contradiction, that $V$ thinks $\kappa$ is singular. Then there is a cofinal sequence $\langle \alpha_\xi\mid\xi<\cf\kappa \rangle \in V_\lambda$ that converges to $\kappa$. Since $j(\cf\kappa)=\cf\kappa$ and $j(\alpha_\xi)=\alpha_\xi$ for all $\xi<\cf\kappa$, $j(\kappa)=\kappa$, which gives our contradiction. Similarly, if $\kappa$ is not a strong limit cardinal in $V$, then there is $\xi<\kappa$ and a surjection $f:\mathcal{P}^{V}(\xi)\to \kappa$ in $V$. Then $f\in V_\lambda$ since $\rank f\le \kappa+3<j^3(\kappa)<\lambda$. (This follows from the fact that the $j^n(\kappa)$ for $n<\omega$ form a strictly increasing sequence.) Now we can derive a contradiction in the usual way by considering $j(f)$.
	    
	    Hence $V_\kappa$ is a model of $\ZFC$. Also, we can see in the extension that $V_\kappa\prec V_{j^n(\kappa)}$ for all $n<\omega$ which shows that $V_{j^\omega(\kappa)}=\bigcup_{n<\omega} V_{j^n(\kappa)}$ is a model of $\ZFC$.
	    Finally, $j\restricts V_{j^\omega(\kappa)}:V_{j^\omega(\kappa)}\to V_{j^\omega(\kappa)}$ and, by the transitivity of $V_{j^\omega(\kappa)}$, $V_{j^\omega(\kappa)}$ satisfies \TIj.
	\end{proof}
	
	As a lower bound, Corazza \cite{Corazza2006} proved that $\mathsf{ZFC+BTEE}$ proves there is an $n$-ineffable cardinal for each (meta-)natural number $n$.
	
	% ZF + WA
	The authors do not know the exact consistency strength of $\mathsf{ZF+WA}$ (the assumption of $\Sigma^j$-Induction is useful to ensure that the critical sequence is total) in the $\ZFC$-context but we can still find a lower bound for it. Suppose that $\kappa$ was the critical point of such an embedding. Then we have that the critical sequence $\lag j^n(\kappa) | n \in \omega \rag$ is definable, although it may not be a set (see Proposition 3.2 of \cite{Corazza2000} for details). From this, we have

	\begin{lemma}[$\mathsf{ZF+WA_0}$]
		If the critical sequence is cofinal over the class of all ordinals, then $\kappa$ is extendible.
	\end{lemma} 
	
	\begin{proof}
		Let $\eta$ be an ordinal. Take $n$ such that $\eta<j^n(\kappa)$, then $j^n \colon V_{\kappa+\eta}\prec V_{j^n(\kappa+\eta)}$ and $\crit j^n=\kappa$. Hence $\kappa$ satisfies $\eta$-extendibility.
	\end{proof}
	
	However, it should be noted that it is also possible that the critical sequence $\lag j^n(\kappa)\mid n<\omega \rag$ is bounded.
	In this case $V_\lambda$, for $\lambda = \sup_{n<\omega} j^n(\kappa)$, is a model of $\mathsf{ZF+WA}$ in which the critical sequence is cofinal.
	Thus we can proceed with the argument by cutting off the universe at $\lambda$.
	
	By an easy reflection argument, we can see also see that $\mathsf{ZF+WA_0}$, with the critical sequence cofinal, proves not only that there is an extendible cardinal, but also the consistency of $\mathsf{ZF}$ with a proper class of extendible cardinals, an extendible limit of extendible cardinals, and much more. 
	Since extendible cardinals are preserved by Woodin's forcing \cite[Theorem 226]{Woodin2010}, we have a lower bound for the consistency strength of $\mathsf{ZF+WA_0}$, e.g., $\ZFC$ plus there is a proper class of extendible cardinals. Clearly, $\mathsf{ZF + WA_0}$ should be much stronger than this but to find a better lower bound that involves even more sophisticated machinery than is currently available.
	
	To summarize our consistency bound derived in this section, we have the following:
	
	\vbox{
	\begin{corollary} \phantom{a}
	    \begin{itemize}
	        \item $\IKP_{j,M}$ with a $\Sigma$-$\Ord$-inary elementary embedding or $\CZF$ with a critical set implies the consistency of $\mathsf{ZFC+BTEE} + \TIj$.
	        Furthermore, $\mathsf{ZFC+BTEE} + \TIj$ proves that the critical point, $\kappa$, of $j$ is $n$-ineffable for every (meta-)natural $n$.
	        
	        \item $\CZF$ with a Reinhardt set implies the consistency of $\mathsf{ZF+WA}$. Furthermore, the consistency $\ZF+\WA$ implies the consistency of $\ZFC$ with a proper class of extendible cardinals.
	    \end{itemize}
	\end{corollary}
	}

	\section{Future works and Questions}\label{Section:RemarkQuestions}
	
	% Deriving the upper bounds
	\subsection{Upper bounds for large large set axioms}
	We may wonder how to find an upper bound for the consistency strength of $\CZF$ with very large set axioms in terms of classical set theories. The authors believe that the currently known methods do not suffice to provide non-trivial upper bounds for the proof-theoretic strength of $\CZF$ with very large set axioms. The known methods for analyzing the strength of $\CZF$ and its extensions are the followings:

	\begin{enumerate}
	    \item Reducing $\CZF$ or its extension to Martin-L\"of type theory or its extension, and constructing a model of the type theory into a classical theory such as $\KP$ or its extensions. This is how Rathjen (cf., \cite{Rathjen1993}, \cite{Rathjen2005Brouwer}, \cite{Rathjen2014Omniscience}, \cite{Rathjen2017}) provides a relative proof-theoretic strength for extensions of $\CZF$.
 
	    \item More generally, sets-as-trees interpretation or its variants: for example, Lubarsky \cite{Lubarsky2006SOA} proved that we can reduce $\mathsf{CZF+Sep}$ into Second-order Arithmetic by combining realizability with a sets-as-trees interpretation. We may associate Lubarsky's construction with Rathjen's interpretations because the sets-as-types interpretation is a special case of the sets-as-trees interpretation. (Types have tree-like structures.) 
	    Another construction on that line is functional realizability, which we define below for the sake of completeness: 
	    
	    \begin{definition}
            Let $\mathcal{A}$ be a pca and $V(\mathcal{A})$ be the realizability universe. (See Definition 3.1 of \cite{Rathjen2003Realizability} or Definition 2.3.1 of \cite{McCartyPhD}.) A name $a\in V(\mathcal{A})$ is \emph{functional} if for every $\lag e,b\rag, \lag e, c\rag \in a$, $b=c$. $V^f(\mathcal{A})$ is the class of all functional names $a\in V(\mathcal{A})$.
            
            The realizability relation $\Vdash^f$ over $V^f(\mathcal{A})$ is identical with the usual realizablity relation $\Vdash$ over $V(\mathcal{A})$. (See Definition 4.1 of \cite{Rathjen2003Realizability}), except that we restrict quantifiers to $V^f(\mathcal{A})$ instead of $V(\mathcal{A})$.
            We say that $V^f(\mathcal{A})\models \phi(\vec{a})$ if there is an $e\in\mathcal{A}$ such that $e\Vdash^f \phi(\vec{a})$.
        \end{definition}
	    
	    \item Set realizability, which appears in \cite{Rathjen2012Power} and \cite{Rathjen2012Existence}. This exploits the computational nature of sets to construct an interpretation.
	    Unfortunately, set realizablity does not result in models of $\CZF$: it produces an interpretation of $\CZF^-$, $\CZF^-+\mathsf{Exp}$, or $\mathsf{CZF+Pow}$ to $\IKP$, $\IKP(\mathcal{E})$, or $\IKP(\mathcal{P})$ respectively. 
	    Furthermore, set realizability can be used to prove the existence property of $\CZF^-$, $\CZF^-+\mathsf{Exp}$, or $\mathsf{CZF+Pow}$.
	    However, Swan proved in \cite{Swan2014} that $\CZF$ does not have the existence property, which suggests set realizability cannot model $\CZF$.
	\end{enumerate}
	
	% Ziegler's subsection 9.4.4
	Thus the only currently known way to analyze the consistency strength of $\CZF$ and its extensions is by combining realizability with sets-as-trees interpretations. However, this has significant issues when we try and generalize the method to large large cardinals. The first point is that almost all of the currently known methods (possibly except for \cite{Rathjen2005Brouwer} and \cite{Swan2014}) rely on Kleene's first pca to construct an interpretation of variants of Martin-L\"of type theory into classical theories.
	The upshot is that the resulting interpretation of $\CZF$ also validates \emph{the Axiom of Subcountability}, which claims that every set is subcountable.
	However, Ziegler \cite{ZieglerPhD} observed that the Axiom of Subcountability is incompatible with critical sets.
    This feature may simply be due to the fact that the pca we are using in the construction of the interpretation is countable, so we may avoid this issue by using a larger pca.

	However, it seems that there is no obvious way of constructing a realizer for $\forall\vec{x}\phi^M(j(\vec{x}))\to \phi(\vec{x})$ regardless of how large the pca $\mathcal{A}$ is.
	Since we can view sets-as-types interpretations as a special case of sets-as-trees interpretations, we can expect there would be a similar difficulty when we construct a model of $\CZF$ with large large set axioms by using the combination of realizability models and type-theoretic interpretations of $\CZF$.
	
	\begin{question}
	    Can we provide any non-trivial upper bound for the consistency strength of $\CZF$ with a critical set or a Reinhardt set?
	\end{question}
	
	\subsection{Improving lower bounds}
	Our current lower bound could also be improved. For example, the proof of \autoref{Corollary:LambdaModelsIZFBTEEInd}, \autoref{Theorem:IKPSigmaOrdimpliesIZF+BTEE}, and \autoref{Theorem:ReinhardtCriticalPtModelsWA} produces a set model of some theory. Hence the resulting lower bound for the consistency strength is strict. This brings into question whether we can provide a better lower bound for the given theories, possibly by constructing a class model of $\IZF$ with a large set axiom.
	
	Since obtaining a lower bound heavily relies on double-negation translations, it would be important to develop the relationship between double-negation translations and large set axioms. For example, Avigad \cite{Avigad2000} provided a way to interpret $\KP$ from $\IKP$ by combining Friedman's double-negation interpretation \cite{Friedman1973} and a proof-theoretic forcing. It may be possible to extend Avigad's interpretation for $\IKP$ with a $\Sigma$-$\Ord$-inary elementary embedding with a critical point, which could result in a better lower bound (or possibly an equiconsistency result). 
	
	The `classical' side of the lower bounds should also be improved. For example, we claimed that $\CZF$ with a Reinhardt set implies the consistency of $\mathsf{ZF+WA}$, and the latter implies the consistency of $\ZFC$ with a proper class of extendible cardinals.
	The authors do not know if it is possible to prove the consistency of $\mathsf{ZFC+WA}$ from that of $\mathsf{ZF+WA}$. 
	
	We conclude this subsection with the following obvious question: 

	\begin{question}
		Can we obtain a better lower bound for the consistency of any of the theories analyzed in this paper?
	\end{question}
	
	\subsection{Developing technical tools}
	Some concepts in this paper are of independent interest. For example, we defined second-order constructive set theories $\CGB$ and $\IGB$ to handle super Reinhardt sets and $\TR$.
	However, constructive second-order set theories bring their own questions, associated with their classical counterparts. The following questions are untouched in this paper:

	\begin{question}\phantom{a}
	    \begin{enumerate}
	        \item Williams \cite{WilliamsPhD} defined and analyzed second-order set-theoretic principles that bolster second-order set theory, including \emph{Class Collection schema}, \emph{Elementary Transfinite Recursion schema}, $\mathsf{ETR}$, and its restrictions $\mathsf{ETR}_\Gamma$.
	        Can we define constructive analogues of these principles? If so, what are their consistency strength? For example, does $\mathsf{CGB+ETR_\omega}$ prove the existence of the truth predicate of first-order set theory (and hence proves $\mathsf{Con(CZF})$)?
	        \item Do constructive second-order set theories admit the unrolling construction that was introduced by Williams \cite{WilliamsPhD}?
	        \item Can we develop a realizability model or Heyting-valued model for constructive second-order theories? It is known that not every class forcing preserves Collection and Powerset over the classical $\mathsf{GB}$. We may expect that we need some restriction on a class realizability or a class formal topology to ensure they preserve $\CGB$.
	\end{enumerate}
	\end{question}
	
	The relativization of Heyting-valued models is also an interesting topic. We only focused on the double-negation topology, and it seems that the formal topologies appearing in the current literature are either set-presentable or the double-negation topology. However, it is plausible that another formal topology might appear in the future, and its interaction with different inner models could be non-trivial. In that case, absoluteness and relativization issues become important.
	
	\vbox{
	\subsection{Other large set axioms}
	In this paper, we analyzed critical sets, Reinhardt sets, and some constructive analogue of choiceless large cardinals. We did not provide the definition and analysis for analogues of other large cardinal axioms, such as supercompactness or hugeness. These concepts were first defined over $\IZF$ by Friedman and \v{S}\v{c}edrov} in \cite{FriedmanScedrov1984}, in such a way as to be equivalent to their classical counterparts. However, for the sake of completeness, we include variations below that work better in our weaker context of $\CZF$.
	
	\begin{convention}
	    Assume that $K$ is a critical point of an elementary embedding $j\colon V\to M$. Over $\ZFC$, many large cardinals above measurable cardinals are defined as critical points of some elementary embedding with additional properties, for example, closure properties of $M$. We shall give constructive analogues of the main methods used to define closure where, here, $\mathcal{MP}$ is Ziegler's modified powerclass operator and $\hat{V}_a$ is Ziegler's modified hierarchy. (See Section 5.3 of \cite{ZieglerPhD} for the details).
	\end{convention}
	
	\begin{itemize}
	    \item Replace the \emph{closure under $<\gamma$-sequences}, $\prescript{< \gamma}{} M\subseteq M$, to \emph{closure under multi-valued functions whose domain is an element of $\hat{V}_a$}: if $b\in \hat{V}_a$ and $R:b\rrarrows M$ is a multi-valued function\footnote{It is unclear whether we need to restrict $R$ to set-sized multi-valued functions. Future research should analyze the difference between these two.}, then there is $c\in M$ such that $R:b\lrlrarrows c$.
	    
	    \item Replace the \emph{closure under $\gamma$-sequences}, $\prescript{\gamma}{} M\subseteq M$, to \emph{closure under multi-valued functions whose domain is an element of $\hat{V}_{a\cup\{a\}}=\hat{V}_a\cup \mathcal{MP}(\hat{V}_a)$}.
    
        \item Replace $V_\alpha \subseteq M$ with $\hat{V}_\alpha\subseteq M$.
	\end{itemize}
	
	The reader might wonder why we use multi-valued functions of the domain in $\hat{V}_a$ and $\hat{V}_{a\cup\{a\}}$. For example, we may formulate a constructive analogue of the closure under $\gamma$-sequences $\prescript{\gamma}{} M\subseteq M$ as the closure under multi-valued functions of domain $a$. There is a reason why we should allow multi-valued functions of the domain in $a\cup\{a\}$: for a transitive set $a$, the closure under multi-valued functions of domain in $a\cup\{a\}$ proves $a\in M$. However, it is unclear if this can be achieved from just those multi-valued functions whose domain is in $a$.
	
	It still remains a question why we use $\hat{V}_{a\cup\{a\}}$ instead of $a\cup\{a\}$. The main reason is that $\hat{V}_{a\cup\{a\}}$ includes $a\cup\{a\}$ but has a much richer set-theoretic structure, which is the same reason we have worked with large sets rather than large cardinals. This will mean that our definitions will be equivalent to those presented in \cite{FriedmanScedrov1984} over $\IZF$.

    Observe that, over $\ZFC$, these modified definitions are equivalent to the standard ones. Thus, for example, we can define supercompact sets or strong sets\footnote{This is an analogue of strong cardinals. Unfortunately, this terminology overlaps with 2-strong sets (See \cite{Rathjen1998} or \cite{ZieglerPhD} for the definition of 2-strongness.)} as follows:
	
	\begin{definition}[$\CGB_\infty$] \label{Definition:Supercompact,Huge,Strong} \phantom{a} 
	    \begin{enumerate}
	        \item Let $a$ be a transitive set. A set $K$ is \emph{$a$-supercompact} if there is an elementary embedding $j\colon V\to M$ such that $K$ is a critical point of $j$ and $M$ is closed under multi-valued functions whose domain is in $\hat{V}_{a\cup\{a\}}$ A set $K$ is \emph{supercompact} if $K$ is $a$-supercompact for all transitive sets $a$.
	        \item A set $K$ is \emph{$n$-huge} if there is an elementary embedding $j\colon V\to M$ such that $K$ is a critical point of $j$ and $M$ is closed under multi-valued functions whose domain is in $\hat{V}_{j^n(K)\cup \{j^n(K)\}}$.
	        \item A set $K$ is \emph{$\alpha$-strong} if there is an elementary embedding $j\colon V\to M$ such that $K$ is a critical point of $j$ and $\hat{V}_\alpha\subseteq M$. A set $K$ is \emph{strong} if $K$ is $\alpha$-strong for all $\alpha$. 
	        
	        Here it suffices to restrict our attention to ordinals rather than defining $a$-strongness for arbitrary sets since $\hat{V}_a = \hat{V}_{\rank(a)}$.
	    \end{enumerate}
	\end{definition}
	
	The reader is reminded that we may formulate $a$-supercompactness or $\alpha$-strongness over $\CZF_{j,M}$. However, formulating the full supercompactness and strongness would require us to quantify over $j$ and $M${, so these should be stated over $\CGB_\infty$}.
	
	As a remark, let us mention that Friedman and \v{S}\v{c}edrov \cite{FriedmanScedrov1984} also defined hugeness and supercompactness over $\IZF$. A set $K$ is \emph{huge in the sense of Friedman-\v{S}\v{c}edrov} if $K$ is a critical point of an elementary embedding $j\colon V\to M$ which is power inaccessible and satisfies the following statement: for any subset $u$ of $j(K)$, if $t\colon u\rrarrows M$ then we can find some $v\in M$ such that $t\colon u\lrlrarrows v$.\footnote{Their definition is adjusted for the intensional $\IZF$, but we work with Extensionality. Also, they required $v$ only to satisfy $t\colon u\rrarrows v$, but their definition is `incorrect' in the sense that their formulation does not imply the closure of $M$ under sequences over $\ZFC$. The reason is that we do not know whether $t$ is amenable over $M$ in general.} The following proposition shows our definition of hugeness and that of Friedman-\v{S}\v{c}edrov coincides in some sense:
	
	\begin{proposition}[$\IGB_\infty$]
	    A transitive set $K$ is huge in the sense of Friedman-\v{S}\v{c}edrov if and only if $K$ is power inaccessible and huge in the sense of \autoref{Definition:Supercompact,Huge,Strong}.
	\end{proposition}
	
	\vbox{
	\begin{proof}
        Assume that $K$ is huge in the sense of Friedman-\v{S}\v{c}edrov. It is known that $\hat{V}_K=K$ if $K$ is inaccessible (See \cite[Theorem 5.18]{ZieglerPhD}). Furthermore, power inaccessibility implies $\mathcal{P}(1)\in K$, so $\mathcal{MP}(K)=\mathcal{P}(K)$. Since transitivity implies $K\subseteq \mathcal{P}(K)$, we have $\hat{V}_K\cup \mathcal{MP}(K)=\mathcal{P}(K)$. Thus we have
        \begin{equation*}
            \forall u\in \hat{V}_K\cup \mathcal{MP}(K) [\forall t\colon u\rrarrows M \to \exists v\in M (t\colon u\lrlrarrows v)].
        \end{equation*}
        Hence $K$ is huge. The remaining direction is trivial.
	\end{proof}
	}
	
	The case for supercompactness is tricky because Friedman and \v{S}\v{c}edrov employed a family of elementary embeddings instead of working over a second-order set theory to formulate it. We leave examining the difference between these two definitions of supercompact sets to possible future work.
	
	These new notions of large set axioms bring the following question:

	\begin{question}
	    Can we provide any consistency result for the large set axioms we can define by following the above schemes? Can we define an $\IKP$-analogue of such very large cardinals? 
	\end{question}
	
\overfullrule=0pt

\printbibliography

\appendix
\clearpage
\section{Tables for notions appearing in this paper}

    \begin{figure}[h]
        \begin{center}
        \renewcommand{\arraystretch}{1.2}
        \begin{tabular}{|c|m{80mm}|c|c|c|}
            \hline
            Name & \centering{Rough definition} & Def.& Synonyms & Note \\
            \hline \hline
            Regular & $M$ {is} transitive and satisfies 2\textsuperscript{nd}-order Strong \hbox{Collection} & \ref{Definition:RegularityandInaccessibility} &  & \\
            \hline
            Weakly regular & $M$ {is} transitive and satisfies 2\textsuperscript{nd}-order Collection & \ref{Definition:RegularityandInaccessibility} &  & \\
            \hline
            Functionally regular & $M$ {is} transitive and satisfies 2\textsuperscript{nd}-order Replacement & \ref{Definition:RegularityandInaccessibility} &  & \\
            \hline
            $\bigcup$-regular & $M$ {is} regular and satisfies Union & \ref{Definition:RegularityandInaccessibility} & & \\
            \hline
            Strongly regular & $M$ {is} $\bigcup$-regular and satisfies Exponentiation & \ref{Definition:RegularityandInaccessibility} &  & \\
            \hline
            BCST-regular & $M$ {is} regular and $M\models \mathsf{BCST}$ & \ref{Definition: BCST-regular} & & \\
            \hline
            Inaccessible & $M$ {is} regular and $M\models \CZF_2$ & \ref{Definition:RegularityandInaccessibility} & & \\
            \hline
            REA-Inaccessible & $M$ {is} inaccessible and $M\models\mathsf{REA}$ & \ref{Definition:REA-inaccessible} & \cite{Rathjen1999Realm}: \emph{Inaccessible} &  \\ 
            \hline
            Power Inaccessible & $M$ {is} inaccessible and $b\subseteq a\in M\to b\in M$ & \ref{Definition:PowerInaccessible} & \cite{FriedmanScedrov1984}: \emph{Inaccessible} & \eqref{Table1item-1} \\
            \hline\hline
            $\mathsf{REA}$ & Cofinally many regular sets & \ref{Definition:RegularityandInaccessibility} & & \\
            \hline
            $\mathsf{IEA}$ & Cofinally many inaccessible sets & \ref{Definition:RegularityandInaccessibility} & & \\
            \hline
            $\mathsf{pIEA}$ & Cofinally many power inaccessible sets & \ref{Definition:PowerInaccessible} & & \\
            
            \hline\hline
            $T_{j,M}$ & $j \colon V \to M$ is an elementary embedding and the axiom schemes of $T$ (e.g. Collection) allow $j$ and $M$ as parameters & \ref{Convention:T_j}& & \eqref{Table1item-2} \\ 
            \hline
            $\Delta_0\mhyphen\BTEE_M$ & $j\colon V\to M$ and $j$ is $\Delta_0$-elementary & \ref{Definition:BTEE} & &\eqref{Table1item-2} \\
            \hline
            $\Sigma\mhyphen\BTEE_M$ & $j\colon V\to M$ and $j$ is $\Sigma$-elementary & \ref{Definition:BTEE} & &\eqref{Table1item-2} \\
            \hline
            $\BTEE_M$ & $j\colon V\to M$ and $j$ is fully elementary & \ref{Definition:BTEE}  & & \eqref{Table1item-2} \\
            \hline
            $\WA$ & $\BTEE$ + Separation for $j$-formulas & \ref{Definition:WA} & & \\
            \hline \hline
            Critical & A set $K$ such that $j\colon V\to M$, $K\in j(K)$, $j\restricts K = \mathrm{Id}$ & \ref{Definition:CriticalReinhardt} & \makecell{\cite{Schlutzenberg2020Extenders}: \emph{$V$-critical}\\ \cite{ZieglerPhD}: \emph{Measurable}} & \eqref{Table1item-3} \\
            \hline
            Reinhardt & Same {as} $V$-critical but {with} $V=M$  & \ref{Definition:CriticalReinhardt} & & \\
            \hline
            Super Reinhardt & A set $K$ such that for every set $a$ there is an elementary embedding $j\colon V\to V$ such that $a, K\in j(K)$ and  $j\restricts K=\mathrm{Id}$& \ref{Definition:SuperReinhardt} & & \eqref{Table1item-4} \\
            \hline
            $A$-Super Reinhardt & Same as above but $j$ is $A$-elementary & \ref{Definition:SuperReinhardt} & & \eqref{Table1item-4} \\
            \hline
            $\TR$ & $\forall^1 A \forall^0 a\exists j\colon V\to V$ ($j$ is $A$-elementary and $a\in j(K)$) & \ref{Definition:TR} & & \eqref{Table1item-4} \\ 
            \hline \hline
            
            $V^\mathcal{S}$ & The Heyting universe over the formal topology $\mathcal{S}$ & \ref{Definition: Heyting interpretation} & & \\
            \hline
            $\Omega$ & The double-negation formal topology $(1, =, \vartriangleleft)$, where $x \vartriangleleft p$ iff $\neg \neg (x \in p)$ & \ref{Definition: double-negation topology} & & \\
            \hline
            $a^\dnot$ & The class $\{x\mid \lnot\lnot(x\in a)\}$ & \makecell{Below\\ \ref{Definition: double-negation topology}} & & \\
            \hline
            $A^{\Omega}$ & The Heyting-valued universe over $\Omega$ relativized to $A$ & \ref{Definition:AOmega} & & \\
            \hline
            $\tilde{A}$ & The $\Omega$-name satisfying $\dom \tilde{A}=A^\Omega$ and $\tilde{A}(x)=\top$ & \ref{Definition:AOmega} & & \eqref{Table1item-5} \\
            \hline
        \end{tabular}
        \renewcommand{\arraystretch}{1.2}
        \end{center}
        
        { \footnotesize 
        \begin{enumerate}
            \item \label{Table1item-1} The existence of a Power Inaccessible set implies Powerset over $\CZF$.
            \item \label{Table1item-2} We omit $M$ when $V=M$.
            \item \label{Table1item-3} We call such $K$ a critical point of $j\colon V\to M$. Also, our critical set is different from the notion of critical cardinals defined in \cite{HayutKaragila2020}.
            \item \label{Table1item-4} This formulation requires a second-order set theory with infinite connectives.
            \item \label{Table1item-5} We use $A^\Omega$ to mean $\tilde{A}$ when the context is clear. 
        \end{enumerate}
        }
        \caption{Notions appearing in this paper}
    \end{figure}

    \clearpage
    \begin{figure}[h]
        \begin{center}
        \renewcommand{\arraystretch}{1.2}
        \begin{tabular}{|c|c|c|c|c|c|c|}
            \hline
            Name & Logic & $\in$-Induction & Infinity & Separation & Collection & Powerset \\
            \hline \hline
            $\mathsf{BCST}$ & IFL & & & $\Delta_0$ & Replacement &\\
            $\IKP$ & IFL & \checkmark & \checkmark & $\Delta_0$ & $\Delta_0$ ($\Sigma$) & \\
            $\KP$ & CFL & \checkmark & \checkmark & $\Delta_0$ & $\Delta_0$ ($\Sigma$) & \\
            $\CZFminus$ & IFL & \checkmark & \checkmark & $\Delta_0$ & Strong Collection & \\
            $\CZF$ & IFL & \checkmark & \checkmark & $\Delta_0$ & Strong Collection& Subset Collection \\
            $\IZF$ & IFL & \checkmark & \checkmark & \checkmark & \checkmark & \checkmark \\
            $\ZFminusrep$ & CFL & \checkmark & \checkmark & \checkmark & Replacement & \\
            $\ZFminus$ & CFL & \checkmark & \checkmark & \checkmark & \checkmark & \\
            $\mathsf{ZF}$ & CFL & \checkmark & \checkmark & \checkmark & \checkmark & \checkmark \\
            $\CGB$ & \eqref{Table2Item-1} & \checkmark & \checkmark & $\Delta_0$ & Class Strong Collection & \\
            $\IGB$ & \eqref{Table2Item-1} & \checkmark & \checkmark & \checkmark & Class Strong Collection & \checkmark \\
            $\CGB_\infty$ & $\mathsf{G3i}_\omega$ & \checkmark & \checkmark & $\Delta_0$ & Class Strong Collection & \\
            $\IGB_\infty$ &  $\mathsf{G3i}_\omega$& \checkmark & \checkmark & \checkmark & Class Strong Collection & \checkmark \\ 
            \hline 
        \end{tabular}
        \end{center}
        { \footnotesize 
        \begin{enumerate}
            \item Every theory on the table has Extensionality, Union, and Pairing as axioms.
            \item IFL stands for Intuitionistic first-order logic, and CFL stands for Classical first-order logic.
            \item $\Delta_0$ means the schema is restricted to $\Delta_0$-formulas. $(\Sigma)$ means we can prove the schema for $\Sigma$-formulas. 
            \item \label{Table2Item-1} We define second-order set theories over intuitionistic first-order logic with two sorts. See \autoref{Definition:CGB}, \ref{Definition:IGB}, and \ref{Definition:CGBIGBinfty} for the definition of $\CGB$, $\IGB$,$\CGB_\infty$, and $\IGB_\infty$ respectively.
        \end{enumerate}
        }
        \caption{Theories appearing in this paper}
    \end{figure}

    % Table ended

\end{document}